\documentclass[12pt]{scrartcl}
\pdfoutput=1
\usepackage{amsmath}
\usepackage{url,color}
\usepackage[pdftex,colorlinks]{hyperref}
\definecolor{darkblue}{cmyk}{1,0,0,0.8}
\definecolor{darkred}{cmyk}{0,1,0,0.7}
\hypersetup{anchorcolor=black,
  citecolor=darkblue, filecolor=darkblue,
  menucolor=darkblue,pagecolor=darkblue,urlcolor=darkblue,linkcolor=darkblue}
\usepackage[scaled=0.8]{helvet}
\usepackage[charter]{mathdesign}
\usepackage{microtype} 
\usepackage{mythmbox}
\typearea{12}

\newtheorem[S]{corollary}{Corollary}[section]
\newtheorem[S]{theorem}[corollary]{Theorem}
\newtheorem[S]{definition}[corollary]{Definition}
\newtheorem[S]{assumption}[corollary]{Assumption}
\newtheorem[S]{lemma}[corollary]{Lemma}
\newtheorem[S]{proposition}[corollary]{Proposition}
\newtheorem[S]{introthm}{Main Result}

\newcommand{\R}{\mathbb{R}}

\newcommand{\T}{\mathbb{T}}
\newcommand{\C}{\mathbb{C}}
\newcommand{\Lint}{\mathbb{L}}

\newcommand{\id}{I}

\newcommand{\rg}{\operatorname{rg}}
\newcommand{\cut}{\operatorname{cut}}
\renewcommand{\mod}{\operatorname{mod}}
\newcommand{\ev}{\operatorname{ev}}

\renewcommand{\d}{\mathop{}\!\mathrm{d}}
\renewcommand{\Re}{\operatorname{Re}}
\renewcommand{\Im}{\operatorname{Im}}
\newcommand{\tpartial}{\partial}
\renewcommand{\phi}{\varphi}
\renewcommand{\epsilon}{\varepsilon}
\newcommand{\pstart}[1]{\paragraph*{#1}}
\newcommand{\eop}{\hfill$\square$}

\title{Finding periodic orbits in state-dependent delay differential
  equations as roots of algebraic equations} \author{Jan Sieber}
\begin{document}
\maketitle
\begin{abstract}
  \noindent 
  In this paper we prove that periodic boundary-value
  problems (BVPs) for delay differential equations are locally
  equivalent to finite-dimensional algebraic systems of equations. 
  We rely only on regularity assumptions that follow those of the
  review by Hartung \emph{et al.}  (2006). Thus, the equivalence
  result can be applied to differential equations with state-dependent
  delays (SD-DDEs), transferring many results of bifurcation theory
  for periodic orbits to this class of systems. We demonstrate this by
  using the equivalence to give an elementary proof of the Hopf
  bifurcation theorem for differential equations with state-dependent
  delays. This is an alternative and extension to the original Hopf
  bifurcation theorem for SD-DDEs by Eichmann (2006).
\end{abstract}
\section{Introduction}
\label{sec:intro}
If a dynamical system is described by a differential equation where
the derivative at the current time may depend on states in the past
one speaks of delay differential or, more generally, functional
differential equations (FDEs).  A reasonably general formulation of an
autonomous dynamical system of this type looks like this:
\begin{equation}\label{eq:ivp}
  \dot x(t)=f(x_t,\mu)
\end{equation}
where $\tau>0$ is an upper bound for the delay. On the right-hand side
$f$ is a functional, mapping $C^0([-\tau,0];\R^n)$ (the space of
continuous functions on the interval $[-\tau,0]$ with values in
$\R^n$) into $\R^n$. The dependent variable $x$ is a function on
$[-\tau,T_{\max})$ for some $T_{\max}>0$, and $x_t$ is the current
\emph{function segment}: $x_t(s)=x(t+s)$ for $s\in[-\tau,0]$ such that
$x_t\in C^0([-\tau,0];\R^n)$. The second argument $\mu\in\R^\nu$ is a
system parameter. For a system of the form \eqref{eq:ivp} one would
have to prescribe a continuous function $x$ on the interval
$[-\tau,0]$ as the initial value and then extend $x$ toward time
$T_{\max}$ (see textbooks on functional differential equations such as
\cite{DGLW95,HL93,S89}). 

A long-standing problem with certain types of FDEs is that they do
not fit well into the general framework of smooth infinite-dimensional
dynamical system theory. The problem occurs whenever the functional
$f$ invokes the evaluation operation in a non-trivial way, that is,
for example, if one has a state-dependent delay. A prototypical caricature
example would be the functional
\begin{align}
  \label{eq:example}
  f:&\ U\times\R\mapsto \R\mbox{,}& f(x,\mu)&=\mu-x(-x(0))\mbox{,}
  \intertext{where $U=\{x\in C^0([-\tau,0];\R): 0<x(0)<\tau\}$ is an
    open set in $C^0([-\tau,0];\R)$. The corresponding FDE is} \dot
  x(t)&=\mu-x(t-x(t))\mbox{.}\label{eq:examplede}
\end{align}
Here, $f$ evaluates its first argument $x$ at a point that itself
depends on $x$. We restrict ourselves to solutions $x$ of
\eqref{eq:examplede} with $x(t)\in(0,\tau)$ for $t\geq0$ to avoid
problems with causality and to limit the maximal delay to $\tau$
(always keeping $x_t$ in $U$).

The difficulty with \eqref{eq:examplede} stems from the fact that $f$
as a map is only as smooth as its argument $x$. Specifically, the
derivative of $f$ with respect to its first argument in this example
exists only for $x\in C^1([-\tau,0];\R)$ (the space of all
continuously differentiable functions on $[-\tau,0]$):
\begin{equation}\label{eq:exdf}
  \begin{split}
    \partial^1 f&:C^1([-\tau,0];\R)\times\R\times
    C^1([-\tau,0];\R)\mapsto\R\mbox{,}\\
    \partial^1f&(x,\mu,y)=x'(-x(0))\,y(0)-y(-x(0))\mbox{.}
  \end{split}
\end{equation}
So, if we choose $U$ as the phase space for initial-value problems
(IVPs) in example~\eqref{eq:examplede} then the functional $f$ is not
differentiable for all elements of $U$. In fact, it is not even
locally Lipschitz continuous in $U$. Indeed, Winston \cite{W70} gave
an example of an initial condition in $U$ for \eqref{eq:examplede}
(with $\mu=0$ and $\tau>1$), for which the IVP did not have a unique
solution. This counterexample is not surprising since the right-hand
side $f$ does not fit into the framework that the textbooks
\cite{DGLW95,HL93,S89} assume to be present. A result of Walther
\cite{W04} rescues IVPs with state-dependent delays (such as
\eqref{eq:examplede}) by restricting the phase space in general to the
closed submanifold $C_c$ of $C^1([-\tau,0];\R^n)$:
\begin{displaymath}
  C_c=\{x\in
  C^1([-\tau,0];\R^n): x'(0)=f(x)\}\mbox{.}
\end{displaymath}
Walther \cite{W04} could prove the existence of a semiflow inside this
manifold that is continuously differentiable with respect to its
initial conditions. However, this result is restricted to a single
degree of differentiability. Results about higher degrees of
smoothness are lacking for the semiflow \cite{HKWW06}.

A typical task one wants to perform for problems of type
\eqref{eq:ivp}, or example \eqref{eq:examplede}, is bifurcation
analysis of \emph{equilibria} and \emph{periodic orbits}. Equilibria
are solutions $x$ of \eqref{eq:ivp} that are constant in time, and
periodic orbits are solutions $x\in C^1(\R;\R^n)$ of \eqref{eq:ivp}
that satisfy $x(t+T)=x(t)$ for some $T>0$ and all $t\in\R$. Equilibria
of the general FDE \eqref{eq:ivp} can be determined by finding the
solutions $(p,\mu)\in\R^n\times\R^\nu$ of the algebraic system of
equations
\begin{equation}
  \label{eq:eqsys}
  0=f(E_0p,\mu)
\end{equation}
where $E_0$ is the trivial embedding
\begin{displaymath}
  E_0:\R^n\mapsto C^0([-\tau,0];\R^n)\mbox{,\quad}
  [E_0p](s)=p\mbox{\quad for all $s\in[-\tau,0]$.}
\end{displaymath}
We observe that, even though the FDE \eqref{eq:ivp} is an
infinite-dimensional system, its equilibria can be found as roots of
the finite-dimensional system \eqref{eq:eqsys} of algebraic
equations. Moreover, the regularity problems of the semiflow do not
affect \eqref{eq:eqsys}: in the example \eqref{eq:examplede}, the
algebraic equation \eqref{eq:eqsys} reads $0=\mu-p$, which is smooth
to arbitrary degree, and can be solved even for negative $\mu$ (near
equilibria with $\mu=p<0$ the semiflow does not exist).

In this paper we establish a system similar to \eqref{eq:eqsys}, but
for periodic orbits: we find a finite-dimensional algebraic system of
equations that does not suffer from the regularity problems affecting
the semiflow, and an equivalence between solutions of this algebraic
system and periodic orbits of \eqref{eq:ivp}. In comparison, for
ordinary differential equations (ODEs) of the form $\dot
x(t)=f(x(t),\mu)$ with a smooth $f:\R^n\times\R^\nu\mapsto\R^n$, the
fact that the problem of finding periodic orbits can be reduced to
algebraic root-finding is well known \cite{GH83}. For example, in ODEs
one can use the algebraic system $0=X(T;p,\mu)-p$ where $t\mapsto
X(t;p,\mu)$ is the trajectory defined by the IVP starting from
$p\in\R^n$ and using parameter $\mu\in\R^\nu$.

A central notion in the construction of the equivalent algebraic
system for periodic orbits of FDEs are \emph{periodic boundary-value
  problems} (BVPs) for FDEs on the interval $[-\pi,\pi]$ with periodic
boundary conditions (which we identify with the unit circle
$\T$). Periodic orbits of \eqref{eq:ivp} can then be found as
solutions of periodic BVPs. If one wants to make the equivalence
result useful in practical applications, one has to find a regularity
(smoothness) condition on the right-hand side $f$ that includes the
class of state-dependent delay equations reviewed in \cite{HKWW06},
while still ensuring that it is possible to prove the existence of an
equivalent algebraic system. We use exactly the same condition as used
by Walther in \cite{W04} to prove the existence of a continuously
differentiable semiflow, the so-called \emph{extendable continuous
  differentiability} (originally introduced as ``almost Frech{\'e}t
differentiability'' in \cite{MNP94}), which implies a restricted form
of local Lipschitz continuity. We generalize restricted continuous
differentiability to higher degrees of restricted smoothness (which we
call $EC^k$ smoothness) in a similar fashion as Krisztin \cite{K03}
did for the proof of the existence and smoothness of local unstable
manifolds of equilibria. Our definition of $EC^k$ smoothness is
comparatively simple to state and check, and lends itself easily to
inductive proofs.

After introducing the notation for periodic BVPs and $EC^k$ smoothness
we state the main result, an equivalence theorem between periodic BVPs
and algebraic systems of equations in Section~\ref{sec:theorem}. The
equivalence theorem reduces statements about existence and smooth
dependence of periodic orbits of FDEs to root-finding problems of
smooth algebraic equations. The result is weaker than the
corresponding results for equilibria of FDEs and for periodic orbits
of ODEs because the equivalence is only valid locally. For any given
periodic function $x_0$ with Lipschitz continuous time derivative we
construct an algebraic system that is equivalent to the periodic BVP
in a sufficiently small open neighborhood of $x_0$. However, the
result is still useful, as we then demonstrate in
Section~\ref{sec:po}. We apply the equivalence theorem in the vicinity
of equilibria for which the linearization of \eqref{eq:ivp} has
eigenvalues on the imaginary axis (for example, near $x_0=\mu=\pi/2$
in example \eqref{eq:examplede}) to prove the Hopf Bifurcation
Theorem. The equivalence theorem reduces the
proof of the Hopf Bifurcation Theorem to an application of the
Algebraic Branching Lemma \cite{AG79}. This provides a complete proof
for the Hopf Bifurcation Theorem for FDEs with state-dependent delays,
including the regularity of the emerging periodic orbits. We discuss
differences to the first version of the proof by Eichmann \cite{E06}
and the approach of Hu and Wu \cite{HW10} in Section~\ref{sec:po}.
The equivalence is applicable in other scenarios where one would
expect branching of periodic solutions. Examples are period doublings,
the branching from periodic orbits with resonant Floquet multipliers
on the unit circle in Arnol'd tongues, and branching scenarios in FDEs
with symmetries. We give a tentative list of straightforward
applications and generalizations of the equivalence theorem in the
conclusion (Section~\ref{sec:conc}).

We note that the theorem stated in Section~\ref{sec:theorem} differs
from statements about numerical approximations. As part of the theorem
we also provide a map $X$ that maps the root of the algebraic system
back into a function space to give the \emph{exact} solution of the
periodic BVP, and a projection $P$ that maps functions to
finite-dimensional vectors (and, hence, periodic orbits to roots of
the algebraic system).  In numerical methods one typically has to
increase the dimension of the algebraic system in order to get more
and more accurate \emph{approximations} of the true solution whereas
the dimension of the algebraic system constructed in
Section~\ref{sec:theorem} is finite.
\section{The Equivalence Theorem}
\label{sec:theorem}
This section states the assumptions and conclusions of the main result
of the paper, the Equivalence Theorem stated in
Theorem~\ref{thm:main}. Before doing so, we introduce some basic
notation (function spaces on intervals with periodic boundary
conditions and projections onto the leading Fourier modes).
\subsubsection*{Periodic BVPs}
We first state precisely what we mean by periodic BVP and introduce
the usual hierarchy of continuous, continuously differentiable and
Lipschitz continuous functions on the compact interval $[-\pi,\pi]$
with periodic boundary conditions.  For $j\geq0$ we will use the
notation $C^j(\T;\R^n)$ for the spaces of all functions $x$ on the
interval $[-\pi,\pi]$ with continuous derivatives up to order $j$
(including order $0$ and $j$) satisfying the periodic boundary conditions
$x^{(l)}(-\pi)=x^{(l)}(\pi)$ for $l=0\ldots j$. Elements of
$C^0(\T;\R^n)$ are continuous and satisfy $x(-\pi)=x(\pi)$. For
derivatives of order $j>0$, $x^{(j)}(-\pi)$ is the right-sided $j$th
derivative of $x$ in $-\pi$, and $x^{(j)}(\pi)$ is the left-sided $j$th
derivative of $x$ in $\pi$. The norm in $C^j(\T;\R^n)$ is
\begin{displaymath}
  \|x\|_j=\max_{t\in[-\pi,\pi]}\left\{|x(t)|,|x'(t)|,\ldots,|x^{(j)}(t)|\right\}\mbox{.}
\end{displaymath}
We can extend any function $x$ in $C^j(\T;\R^n)$ to arguments in $\R$
by defining $x(t)=x(t-2k\pi)$ where $k$ is an integer chosen such that
$-\pi\leq t-2k\pi<\pi$ (we will write $t_{\mod[-\pi,\pi)}$
later). Thus, every element of $C^j(\T;\R^n)$ is also an element of
$BC^j(\R;\R^n)$, the space of functions with bounded continuous
derivatives up to order $j$ on the real line. We use the notation
$t\in\T$ for arguments $t$ of $x$, and also call $\T$ the unit
circle. This make sense because the parametrization of the unit circle
by angle provides a cover, identifying $\T$ with $\R$ where we use
$[-\pi,\pi)$ as the fundamental interval.

Additional useful function spaces are the space of Lipschitz
continuous functions and, correspondingly, spaces with Lipschitz
continuous derivatives, denoted by $C^{j,1}(\T;\R^n)$, which are
equipped with the norm
\begin{equation}\label{eq:cj1def}
  \|x\|_{j,1}=\max\left\{\|x\|_j,
    \sup_{
      \begin{subarray}{c}
        t,s\in\R\\[0.2ex]
        t\neq s
      \end{subarray}
}\frac{|x^{(j)}(s)-x^{(j)}(t)|}{|s-t|} \right\}
\end{equation}
($x^{(0)}(t)$ refers to $x(t)$). Note that we used the notation
$t,s\in\R$ in the index of the supremum, as we can apply arbitrary
arguments in $\R$ to a function $x\in C^0(\T;\R^n)$ by considering it
as an element of $BC^0(\R;\R^n)$, as explained above. We use the same
notation ($C^j(J;\R^n)$ and $C^{j,1}(J;\R^n)$) also for functions on
an arbitrary compact interval $J\subset\R$ without periodic boundary
conditions (and one-sided derivatives at the boundaries). As any
function $x\in C^j(\T;\R^n)$ is also an element of $BC^j(\R;\R^n)$, it
is also an element of $C^j(J;\R^n)$ for any compact interval $J$ (and
the norm of the embedding operator equals unity). On the function
spaces $C^j(\T;\R^n)$ we define the time shift operator
\begin{displaymath}
  \Delta_t:C^j(\T;\R^n)\mapsto C^j(\T;\R^n)\mbox{,\qquad}
  [\Delta_tx](s)=x(t+s)\mbox{.}
\end{displaymath}
The operator $\Delta_t$ is linear and has norm $1$ in all spaces
$C^j(\T;\R^n)$. Similarly, $\Delta_t$ maps also
$C^{j,1}(\T;\R^n)\mapsto C^{j,1}(\T;\R^n)$, and has unit norm there as
well.

Let $f$ be a continuous functional on the space of continuous periodic
functions, that is,
\begin{displaymath}
f:C^0(\T;\R^n)\mapsto\R^n\mbox{.}
\end{displaymath}
The right-hand side $f$, together with the shift $\Delta_t$, creates
an operator in $C^0(\T;\R^n)$, defined as
\begin{align} 
  \label{eq:fdef}
  F&:C^0(\T;\R^n)\mapsto C^0(\T;\R^n) &
  [F(x)](t)&=f(\Delta_tx)\mbox{.}
\end{align}
The operator $F$ is invariant with respect to time shift by construction:
$F(\Delta_tx)=\Delta_tF(x)$.  We consider autonomous periodic
boundary-value problems for differential equations where $f$ is the
right-hand side:
\begin{equation}
  \label{eq:perbvp}
  \begin{split}
    \dot x(t)&=f(\Delta_tx)=F(x)(t)\mbox{.}
  \end{split}
\end{equation}
A function $x\in \C^1(\T;\R^n)$ is a solution of \eqref{eq:perbvp} if
$x$ satisfies equation~\eqref{eq:perbvp} for all $t\in\T$ (for each
$t\in\T$ equation \eqref{eq:perbvp} is an equation in $\R^n$).  In
contrast to the introduction we do not expressly include a parameter
$\mu$ as an argument of $f$. This does not reduce generality as we
will explain in Section~\ref{sec:po}. The main result, the Equivalence
Theorem~\ref{thm:main}, will be concerned with equivalence of the
periodic BVP \eqref{eq:perbvp} to an algebraic system of
equations. The notion of the shift $\Delta_t$ on the unit circle and
the operator $F$, combining $f$ with the shift, is specific to
periodic BVPs such that the BVP \eqref{eq:perbvp} looks different from
the IVP~\eqref{eq:ivp} in the introduction. Several results stating
how regularity of $f$ transfers to regularity of $F$ are collected in
Appendix~\ref{sec:basicprop}.

\subsubsection*{Definition of $EC^k$ smoothness and local (restricted)
  $EC$ Lipschitz continuity}
Continuity of the functional $f$ is not strong enough as a condition
to prove the Equivalence Theorem. Rather, we need a notion of
smoothness for $f$. However, as explained in the introduction, we
cannot assume that $f$ is continuously differentiable with degree
$k\geq1$, if we want to include examples such as $f(x)=-x(-x(0))$ (see
FDE~\eqref{eq:examplede} for $\mu=0$) into the class under
consideration.

The review by Hartung \emph{et al.} \cite{HKWW06} observed the
following typical property of functionals $f$ appearing in equations
of type \eqref{eq:perbvp}: the derivative $\tpartial^1f(x)$ of $f$ in
$x$ as a linear map from $C^1(\T;\R^n)$ into $\R^n$ can be extended to
a bounded linear map from $C^0(\T;\R^n)$ into $\R^n$, and the mapping
\begin{displaymath}
  \tpartial^1f: C^1(\T;\R^n)\times C^0(\T;\R^n) \mapsto \R^n
  \mbox{\quad defined by\quad}
  (x,y)\mapsto\tpartial^1 f(x,y)
\end{displaymath}
is continuous as a function of both arguments. In other words, the
derivative of $f$ may depend on $x'$ but not on $y'$. For the example
$f(x)=-x(-x(0))$ this is true (see \eqref{eq:exdf}). Most of the
fundamental results establishing basic dynamical systems properties
for FDEs with state-dependent delay in \cite{HKWW06} rest on this
extendability of $\partial^1f$.


We also rely strongly on this notion of \emph{extendable}
continuous differentiability. The precise definition is given below in
Definition~\ref{def:extdiff}. In this definition we permit the
argument range $J$ to be any compact interval or $\T$. We use the
notation of a subspace of higher-order continuous differentiability
not only for $C^j(J;\R^n)$ but also for products of such spaces in a
natural way. Say, if
\begin{equation}\label{eq:Dspacedef}
  D=C^{k_1}(J;\R^{m_1})\times\ldots\times C^{k_\ell}(J;\R^{m_\ell})
  \mbox{,}
\end{equation}
where $\ell\geq 1$, and $k_j\geq0$ and $m_j\geq1$ are integers, and
denoting the natural maximum norm on the product $D$ by
\begin{displaymath}
\|x\|_D=\|(x_1,\ldots,x_\ell)\|_D=\max_{j\in\{1,\ldots,\ell\}}\|x_j\|_{k_j}\mbox{,}
\end{displaymath}
then for integers $r\geq0$ the space $D^r$ is defined in the natural
way as
\begin{align*}
  D^r&=C^{k_1+r}(J;\R^{m_1})\times\ldots\times C^{k_\ell+r}(J;\R^{m_\ell})
  \mbox{,\quad with \ }\\
  \|x\|_{D,r}&=\max_{
    \begin{subarray}{c}
      0\leq j\leq r\\[0.2ex]
      1\leq i\leq \ell
    \end{subarray}
} \|x_i^{(j)}\|_{k_i}\mbox{.}
\end{align*}
For the simplest example, $D=C^0(J;\R^n)$, $D^k$ is $C^k(J;\R^n)$. If
$J=\T$ then the time shift $\Delta_t$ extends naturally to products of spaces:
\begin{displaymath}
  \Delta_tx=(\Delta_tx_1,\ldots,\Delta_tx_\ell)
  \mbox{\quad for $x=(x_1,\ldots,x_\ell)\in D$.}
\end{displaymath}

\begin{definition}[Extendable continuous differentiability
  $EC^k$]\label{def:extdiff}
  Let $D$ be a product space of the type \eqref{eq:Dspacedef}, and let
  $f:D\mapsto\R^n$ be continuous. We say that $f$ has an extendable
  continuous derivative if there exists a map $\partial^1f$
  \begin{displaymath}
    \partial^1f:D^1\times  D\mapsto \R^n
  \end{displaymath}
  that is continuous in both arguments $(u,v)\in D^1\times D$ and
  linear in its second argument $v\in D$, such that for all $u\in D^1$
  \begin{equation}\label{eq:ass:contdiff:j}
    \lim_{
      \begin{subarray}{c}
        v\in D^1\\[0.2ex]
        \|v\|_{D,1}\to 0          
      \end{subarray}
      } \frac{|f(u+v)-f(u)-\partial_1f(u,v)|}{\|v\|_{D,1}}=0\mbox{.}
  \end{equation}
  We say that $f$ is $k$ times continuously differentiable in this
  extendable sense if the map $\partial^kf$, recursively defined as
  $\partial^kf=\partial^1[\partial^{k-1}f]$, exists and satisfies the
  limit condition \eqref{eq:ass:contdiff:j} for $\partial^{k-1}f$. We
  abbreviate this notion by saying that $f$ is $EC^k$ smooth in $D$.
\end{definition}
The limit in \eqref{eq:ass:contdiff:j} is a limit in $\R$. For $k=1$
the definition is identical to property (S) in the review
\cite{HKWW06}, one of the central assumptions for fundamental results
on the semiflow.  Extendable continuous differentiability requires the
derivative to exist only in points in $D^1$ and with respect to
deviations in $D^1$, but it demands that the derivative must extend in
its second argument to $D$ ($\partial^1f$ is linear in its second
argument). This is the motivation for calling this property
\emph{extendable} continuous differentiability.  

The definition of $EC^k$ smoothness for $k>1$ uses the notation that a
functional (say, $\partial^1f$) of two arguments (say, $u\in D^1$ and
$v\in D$) for which one would write $\partial^1f(u,v)$, is also a
functional of a single argument $w=(u,v)\in D^1\times D$, such that
one can also write $\partial^1f(w)$.  When using this notation we observe that
the space $D^1\times D$ is again a product of type
\eqref{eq:Dspacedef} such that $\partial^1f$ is again a functional of
the same structure as $f$. For example, let us consider the functional
$f:x\mapsto -x(-x(0))$ from example \eqref{eq:example} (setting
$\mu=0$). The functional is well defined and continuous also on
$D=C^0(\T;\R)$. Moreover, $f$ is $EC^k$ smooth in $D$ to arbitrary
degree $k$. Its first two derivatives are:
\begin{align}
  &\partial^1f: C^1(\T;\R)\times C^0(\T;\R)\mbox{,}\nonumber\\
  &\partial^1f(u,v)=u'(-u(0))\,v(0)-v(-u(0))\mbox{, and}\label{eq:exampledf1}\\
  &\partial^2f: \left[C^2(\T;\R)\times C^1(\T;\R)\right]\times
  \left[C^1(\T;\R)\times C^0(\T;\R)\right]\mbox{,}\nonumber\\
  &
   \begin{aligned}
     \partial_2f(u,v,w,x)=&-u''(-u(0))\,w(0)\,v(0)+u'(-u(0))\,x(0)\\
     &+w'(-u(0))\,v(0)- v'(-u(0))\,w(0)-x(-u(0))\mbox{.}
   \end{aligned}\label{eq:exampledf2}
\end{align}
As one can see, the first derivative $\partial^1f$ has the same
structure as $f$ itself if we replace $D=C^0(\T;\R^n)$ by $D^1\times
D$. So, it is natural to apply the definition again to $\partial^1f$
on the space $D^1\times D$.

Assuming that $f$ is $EC^1$ smooth on $C^0(J;\R^n)$ implies classical
continuous differentiability of $f$ as a map from $C^1(J;\R^n)$ into
$\R^n$ and is, thus, strictly stronger than assuming that $f$ is
continuously differentiable on $C^1(\T;\R^n)$.

Since every element of $C^j(\T;\R^n)$ is also an element of
$C^j(J;\R^n)$ for any compact interval $J$ (and the embedding operator
has unit norm), any $EC^k$ smooth functional $f:C^0(J;\R^m)\mapsto\R^n$
is also a $EC^k$ smooth functional from $C^0(\T;\R^m)$ into $\R^n$.

It is worth comparing Definition~\ref{def:extdiff} with the definition
for higher degree of regularity used by Krisztin in \cite{K03}. With
Definition~\ref{def:extdiff} the $k$th derivative has $2^k$
arguments. In contrast to this, the $k$th derivative as defined in
\cite{K03} has only $k+1$ arguments (the first argument is the base
point, and the derivative is a $k$-linear form in the other $k$
arguments). The origin of this difference can be understood by looking
at the example $f(x)=-x(-x(0))$ and its derivatives in
\eqref{eq:exampledf1}--\eqref{eq:exampledf2}. Krisztin's definition
applied to the second derivative does not include the derivative of
$\partial^1f$ with respect to the linear second argument $v$ (as is
often convention, because it is the identity). One would obtain the
second derivative according to Krisztin's definition by setting $x=0$
in \eqref{eq:exampledf2}. Indeed, the terms containing the argument
$x$ in \eqref{eq:exampledf2} are simply $\partial^1f(u,x)$, as one
expects when differentiating $\partial^1f(u,v)$ with respect to $v$,
calling the deviation $x$. While in practical examples it is often
more economical to use the compact notation with $k$-forms, inductive
proofs of higher-order differentiability using the full derivative
only require the notion of at most bi-linear forms, making them less
complex.

If $f$ is $EC^1$ smooth then it automatically
satisfies a restricted form of local Lipschitz continuity \cite{HKWW06},
which we call local $EC$ Lipschitz continuity:
\begin{definition}[Local $EC$ Lipschitz continuity]\label{def:loclip}
  We say that $f:C^0(\T;\R^n)\mapsto \R^n$ is locally $EC$ Lipschitz
  continuous if for every $x_0\in C^1(\T;\R^n)$ there exists a
  neighborhood $U(x_0)\subset C^1(\T;\R^n)$ and a constant $K$ such
  that
  \begin{equation}\label{eq:loclip}
    |f(y)-f(z)|\leq K\|y-z\|_0
  \end{equation}
  holds for all $y$ and $z$ in $U(x_0)$.  
\end{definition}
That $EC^1$ smoothness implies local $EC$ Lipschitz continuity has
been shown, for example, in \cite{W04} (but see also
Lemma~\ref{thm:flip} in Appendix~\ref{sec:basicprop}).  Note that the
estimate \eqref{eq:loclip} uses the $\|\cdot\|_0$-norm for the upper
bound. This is a sharper estimate than one would obtain using the
expected $\|\cdot\|_1$-norm. The constant $K$ may depend on the
derivatives of the elements in $U(x_0)$ though. For example, for
$f(x)=-x(-x(0))$ as in \eqref{eq:example} with $\mu=0$, one would have
the estimate
\begin{displaymath}
  \left|f(x+y)-f(x)\right|\leq \left[1+\|x'\|_0\right]\|y\|_0
\mbox{\quad such that\quad} 
K\leq1+\max_{x\in U(x_0)}\|x\|_1\mbox{.}
\end{displaymath}
This means that in this example, the neighborhood $U(x_0)$ can be
chosen arbitrarily large as long as it is bounded in
$C^1(\T;\R^n)$. 

The following lemma states that we can extend the neighborhood $U(x)$
for $x\in C^1$ in Definition~\ref{def:loclip} into the space of
Lipschitz continuous functions ($C^{0,1}$ instead of $C^1$) and
include time shifts (which possibly increases the bound $K$).
\begin{lemma}[$EC$ Lipschitz continuity uniform in
  time]\label{thm:loclip}
  Let $f$ be locally $EC$ Lipschitz continuous, and let $x_0$ be in
  $C^1(\T;\R^n)$. Then there exists a bounded neighborhood
  $U(x_0)\subset C^{0,1}(\T;\R^n)$ and a constant $K$ such that
  \begin{displaymath}
    |f(\Delta_ty)-f(\Delta_tz)|\leq K\|\Delta_ty-\Delta_tz\|_0=K\|y-z\|_0
  \end{displaymath}
  holds for all $y$ and $z$ in $U(x_0)$, and for all $t\in\T$. Thus,
  \begin{math}
    \|F(y)-F(z)\|_0\leq K\|y-z\|_0
  \end{math}
  for all $y$ and $z$ in $U(x_0)$.
\end{lemma}
Recall that $F(x)(t)=f(\Delta_tx)$. See Lemma~\ref{thm:flip} and
Lemma~\ref{thm:Flipbound} in Appendix~\ref{sec:basicprop} for the
proof of Lemma~\ref{thm:loclip} and note that we require that the center $x$ of the $C^{0,1}$-neighborhood is $C^1$ to ensure continuity of $t\to\Delta_tx$.  
A consequence of Lemma~\ref{thm:loclip} is that the time
derivative of a solution $x_0$ of the periodic BVP is also Lipschitz
continuous (in time): if $\dot x_0(t)=f(\Delta_tx_0)$ then there
exists a constant $K$ such that
\begin{equation}\label{eq:lipx0t}
  \|x_0'(t)-x_0'(s)\|_0\leq K|t-s|
\end{equation}
Thus, $x_0\in C^{1,1}(\T;\R^n)$. This follows from
Lemma~\ref{thm:loclip} by inserting $\Delta_tx_0$ and
$\Delta_sx_0$ for $y$ and $z$ and using that $x_0'(t)=f(\Delta_tx_0)$
(it is enough to show \eqref{eq:lipx0t} for $|t-s|$ small).

\subsubsection*{Projections onto subspaces spanned by Fourier modes}
The variables of the algebraic system in the Equivalence
Theorem will be the coefficients of the first $N$
Fourier modes (where $N$ will be determined as sufficiently large
later) of elements of $C^{0,1}(\T;\R^n)$ (the space of Lipschitz
continuous functions on $\T$). Consider the functions on $\T$
\begin{displaymath}
  b_0=t\mapsto\frac{1}{2}\mbox{,\quad}b_k=t\mapsto\cos(kt)
  \mbox{,\quad}b_{-k}=t\mapsto\sin(kt)
\end{displaymath}
for $k=1,\ldots,\infty$ (which is the classical Fourier basis of
$\Lint^2(\T;\R)$). For any $m\geq1$ we define the projectors and maps
\begin{equation}\label{eq:proj}
  \begin{aligned}
    P_N&:C^j(\T;\R^m)\mapsto C^j(\T;\R^m)\mbox{,}& [P_Nx](t)_i&=\sum_{k=-N}^N
    \left[\frac{1}{\pi}\int_{-\pi}^\pi b_k(s)x_i(s)\d s\right]\,b_k(t)\mbox{,}
    \allowdisplaybreaks\\
    Q_N&=\id-P_N\mbox{,}\\
    E_N&:\R^{m\times (2N+1)}\mapsto C^j(\T;\R^m)\mbox{,}&
    [E_Np](t)_i&=\sum_{k=-N}^Np_{i,k}b_k(t)\mbox{,}\allowdisplaybreaks\\
    R_N&:C^j(\T;\R^m)\mapsto \R^{m\times(2N+1)}\mbox{,} &
    [R_Nx]_{i,k}\ \ &=\frac{1}{\pi}\int_{-\pi}^\pi b_k(s)x_i(s)\d s\mbox{,}
    \allowdisplaybreaks\\
    L_{\phantom{N}}&:C^j(\T;\R^m)\mapsto C^j(\T;\R^m)\mbox{,}& [Lx](t)\ \ \,&=\int_0^t x(s)-
    R_0x\,\d s=\int_0^t Q_0[x](s)\d s\mbox{.}
  \end{aligned}
\end{equation}
The projector $P_N$ projects a periodic function onto the subspace
spanned by the first $2N+1$ Fourier modes, and $Q_N$ is its
complement. The map $E_N$ maps a vector $p$ of $2N+1$ Fourier
coefficients (which are each vectors of length $n$ themselves) to the
periodic function that has these Fourier coefficients. The map $R_N$
extracts the first $2N+1$ Fourier coefficients from a function. The
simple relation $P_N=E_NR_N$ holds. The vector $R_0x$ is the average
of a function $x$, and $Q_0$ subtracts the average from a periodic
function. The operator $L$ takes the anti-derivative of a periodic
function after subtracting its average (to ensure that $L$ maps back
into the space of periodic functions). In all of the definitions the
degree $j$ of smoothness of the vector space $C^j$ can be any
non-negative integer. The operator $L$ not only maps $C^j$ back into
itself, but it maps $C^j(\T;\R^m)$ into $C^{j+1}(\T;\R^m)$.

We do not attach an index $m$ to the various maps to indicate how many
dimensions the argument and, hence, the value has because there is no
room for confusion: for example, if $x\in C^0(\T;\R^2)$ then $P_Nx\in
C^0(\T;\R^2)$ such that we use the same notation $P_Nx$ for
$x:\T\mapsto\R^m$ with arbitrary $m$. Similarly, we apply all maps
also on product spaces $D$ of the type
$C^{k_1}(\T;\R^{m_1})\times\ldots\times C^{k_\ell}(\T;\R^{m_\ell})$
introduced in Equation~\eqref{eq:Dspacedef} by applying the maps
element-wise. For example,
\begin{align*}
    P_Nx&=(P_Nx_1,\ldots,P_Nx_\ell) &&\mbox{for
      $x=(x_1,\ldots,x_\ell)\in D$,}\\
    E_Np&=(E_Np_1,\ldots,E_Np_\ell) &&\mbox{for
      $p=(p_1,\ldots,p_\ell)\in 
      \R^{m_1\times(2N+1)}\times\ldots\times\R^{m_\ell\times(2N+1)}$.}
\end{align*}

\subsubsection*{Equivalent integral equation}
\label{sec:inteq}
We note the fact that a function $x\in C^1(\T;\R^n)$ solves the
periodic BVP $\dot x(t)=f(\Delta_tx)=F(x)(t)$ if and only if it
satisfies the equivalent integral equation
\begin{equation}\label{eq:inteq}
  x(t)=x(0)+\int_0^tF(x)(s)\d s\mbox{\quad for all $t\in\T$.}
\end{equation}
For each $t\in\T$, Equation \eqref{eq:inteq} is an equation in
$\R^n$. In particular, the term $x(0)$ is in $\R^n$.  Thus, the
integral equation~\eqref{eq:inteq} is very similar to the
corresponding integral equation used in the proof of the
Picard-Lindel{\"o}f Theorem for ODEs \cite{CL55}. This is in contrast
to the abstract integral equations used by Diekmann \emph{et al.}
\cite{DGLW95} to construct unique solutions to IVPs, in which equality
at every point in time is an equality in function spaces. It is the
similarity of \eqref{eq:inteq} to its ODE equivalent that makes the
reduction of periodic BVPs to finite dimensional algebraic equations
possible. One minor problem is that the Picard iteration for
\eqref{eq:inteq} cannot be expected to converge. In fact, the integral
term $\int_0^tF(x)(s)\d s$ does not even map back into the space
$C^0(\T;\R^n)$ of periodic functions, even if $x$ is in
$C^0(\T;\R^n)$. However, a simple algebraic manipulation using the
newly introduced maps $L$, $P_N$, $Q_N$, $E_N$ and $R_N$ removes this
problem (remember that $F(x)(t)=f(\Delta_tx)$):
\begin{lemma}[Splitting of BVP]\label{thm:split}
  Let $N\geq0$ be an arbitrary integer. A function $x\in C^0(\T;\R^n)$
  and a vector $p\in\R^{n\times(2\,N+1)}$ satisfy
  \begin{align}
    \dot x(t)&=f(\Delta_tx)\mbox{\quad and\quad}
    p=R_Nx\mbox{,}\label{eq:parbvp}\\
    \intertext{if and only if they satisfy the system}
    x&=E_Np+Q_NLF(x)\mbox{,}\label{eq:fixp:intro}\\
    0&=R_N\left[P_0F(x)+Q_0\left(E_Np-P_NLF(x)\right)\right]
    \label{eq:lowmodes:intro}\mbox{.}
  \end{align}
\end{lemma}
Note that the map $R_N$ extracts the lowest $2N+1$ Fourier
coefficients from a periodic function. Equation~\eqref{eq:fixp:intro}
can be viewed as a fixed-point equation in a small $C^{0,1}(\T;\R^n)$
neighborhood of the $C^1$ function $x$, parametrized by $p$. We will apply the
Picard iteration to this fixed-point equation instead of
\eqref{eq:inteq}. Equation~\eqref{eq:lowmodes:intro} is an equation in
$\R^{n\times(2\,N+1)}$.  If the Picard iteration converges then the
fixed-point equation \eqref{eq:fixp:intro} can be used to construct
(for sufficiently large $N$) a map
$X:U\subset\R^{n\times(2\,N+1)}\mapsto C^{0,1}(\T;\R^n)$, which maps
the parameter $p$ to its corresponding fixed point $x$. Inserting this
fixed point $x=X(p)$ into \eqref{eq:lowmodes:intro} turns
\eqref{eq:lowmodes:intro} into a system of $n\times(2\,N+1)$ algebraic
equations for the $n\times(2\,N+1)$-dimensional variable $p$, making
the periodic BVP for $x$ equivalent to an algebraic system for its
first $2N+1$ Fourier coefficients, $p$. The proof of
Lemma~\ref{thm:split} is simple algebra, see
Section~\ref{sec:splitproof}.
\subsubsection*{Statement of the Equivalence Theorem}
Using the Splitting Lemma~\ref{thm:split} we can now state the central
result of the paper. The intention to treat \eqref{eq:fixp:intro} as a
fixed-point equation motivates the introduction of the map
\begin{align*}
  M_N&:C^{0,1}(\T;\R^n)\times \R^{n\times(2\,N+1)}\mapsto C^{0,1}(\T;\R^n)
  \mbox{\quad given by}\\
  M_N&(x,p)=E_Np+Q_NLF(x)\mbox{.}
\end{align*}
This means that we will look for fixed points of the map
$M_N(\cdot,p)$ for given $p$ and sufficiently large $N$. We will do
this in small closed balls in $C^{0,1}(\T;\R^n)$ (the space of
Lipschitz continuous functions) such that it is useful to introduce
the notation
\begin{displaymath}
  B_\delta^{0,1}(x_0)=\left\{x\in C^{0,1}(\T;\R^n): 
    \left\|x-x_0\right\|_{0,1}\leq\delta\right\}\mbox{,}
\end{displaymath}
for $\delta>0$ and $x_0\in C^{0,1}(\T;\R^n)$. That is,
$B_\delta^{0,1}(x_0)$ is the closed ball of radius $\delta$ around
$x_0\in C^{0,1}(\T;\R^n)$ in the $\|\cdot\|_{0,1}$-norm (the Lipschitz
norm on $\T$).
\begin{theorem}[Equivalence between periodic BVPs and algebraic
  systems of equations]\label{thm:main}
  Let $f$ be $EC^{j_{\max}}$ smooth, and let $x_0$ have a Lipschitz continuous
  derivative, that is, $x_0\in C^{1,1}(\T;\R^n)$.  Then there exist a
  $\delta>0$ and a positive integer $N$ such that the map
  $M_N(\cdot,p)$ has a unique fixed point in $B_{6\delta}^{0,1}(x_0)$ for
  all $p$ in the neighborhood $U$ of $R_Nx_0$ given by
  \begin{displaymath}
    U=\left\{p\in\R^{n\times(2\,N+1)}: 
      \left\|E_N\left[p-R_Nx_0\right]\right\|_{0,1}\leq2\delta\right\}\mbox{.}
  \end{displaymath}
  The maps
  \begin{align*}
    X&:U\mapsto C^0(\T;\R^n)\mbox{,}&
    X(p)&=\mbox{\ fixed point of $M_N(\cdot,p)$ in $B_{6\delta}^{0,1}(x_0)$,}\\
    g&:U\mapsto\R^{n\times(2\,N+1)}\mbox{,} &
    g(p)&=R_N\left[P_0F(X(p))+Q_0\left(E_Np-P_NLF(X(p))\right)\right]\mbox{,}
  \end{align*}
  are $j_{\max}$ times continuously differentiable with respect to
  their argument $p$, and $X(p)$ is an element of
  $C^{j_{\max}+1}(\T;\R^n)$. Moreover, for all $x\in
  B_{\delta}^{0,1}(x_0)$ the following equivalence holds:
  \begin{align*}
    \dot x(t)&=f(\Delta_tx)\\
    \intertext{if and only if $p=R_Nx$ is in $U$ and satisfies}
    g(p)&=0\mbox{\quad and\quad} x=X(p)\mbox{.}
  \end{align*}
\end{theorem}
Theorem~\ref{thm:main} is the central result of the paper. It implies
that, for any $x_0\in C^{1,1}(\T;\R^n)$ all solutions of the periodic
BVP in a sufficiently small neighborhood of $x_0$ lie in the graph
$X(U)$ of a finite-dimensional manifold. Moreover, these solutions can
be determined by finding the roots of $g$ in $U\subset
\R^{n\times(2\,N+1)}$. We note that Theorem~\ref{thm:main} is
different from statements about numerical approximations. Even though
the integer $N$ is finite, solving the algebraic system $g(p)=0$ and
then mapping the solutions with the map $X$ into the function space
$C^0(\T;\R^n)$ gives an exact solution $x=X(p)$ of the periodic BVP
$\dot x(t)=f(\Delta_tx)$.

The size of the radius $\delta$ of the ball in which the equivalence
holds depends on how large one can choose $\delta$ such that a local
$EC$ Lipschitz constant $K$ for $F$ exists for
$B_{6\delta}^{0,1}(x_0)$ (such neighborhoods exist according to
Lemma~\ref{thm:loclip}). In many applications (in particular, in the
example \eqref{eq:examplede}) this can be any closed ball in which the
right-hand side $f$ is well defined (at the expense of increasing $K$
for larger balls). Once the local $EC$ Lipschitz constant $K$ is
determined, one can find a uniform upper bound $R$ for the norm
$\|F(x)\|_{0,1}$ for all $x\in B_{6\delta}^{0,1}(x_0)$ (see
Lemma~\ref{thm:Flipbound}). The integer $N$, which determines the
dimension of the algebraic system, is then chosen depending on $R$,
$K$ and $\|x_0'\|_{0,1}$. 

Section~\ref{sec:remproofs} contains the complete proof of
Theorem~\ref{thm:main}. The first step of the proof of Equivalence
Theorem~\ref{thm:main} is the existence of the fixed point of $M_N$ in
$B_{6\delta}^{0,1}(x_0)$ for $p\in U$. This is achieved by applying
Banach's contraction mapping principle to the map $M_N(\cdot,p)$ in
the closed ball $B_{6\delta}^{0,1}(x_0)$. The only peculiarity in this
step is that we apply the principle to $B_{6\delta}^{0,1}(x_0)$, which
is a closed bounded set of Lipschitz continuous functions, using the
(weaker) maximum norm ($\|\cdot\|_0$). This is possible because closed
balls in $C^{0,1}(\T;\R^n)$ are complete also with respect to the norm
$\|\cdot\|_0$. With respect to the maximum norm the map
$M_N(\cdot,p):x\mapsto E_Np+Q_NLF(x)$ becomes a contraction for
sufficiently large $N$ (because the norm of the operator $Q_NL$ is
bounded by $C\log(N)/N$, and $F$ has a Lipschitz constant $K$ with
respect to $\|\cdot\|_0$ in $B_{6\delta}^{0,1}(x_0)$).

After the existence of the fixed point of $M_N(\cdot,p)$ is
established in Section~\ref{sec:fixproof} the equivalence between the
algebraic system $g(p)=0$ and the periodic BVP $\dot
x(t)=f(\Delta_tx)$ in the smaller ball $B_\delta^{0,1}(x_0)$ follows
from the Splitting Lemma~\ref{thm:split}.

The smoothness (in the classical sense) of the maps $X$ and $g$
follows, colloquially speaking, from implicit differentiation of the
fixed-point problem $x=E_Np+Q_NLF(x)$ with respect to $p$.
Section~\ref{sec:algdiff1} checks the uniform convergence of the
difference quotient in detail, Section~\ref{sec:smooth} uses the
higher degrees of $EC^{j_{\max}}$ smoothness of $f$ to prove higher
degrees of smoothness for $X$ and $g$. For proving higher-order
smoothness one has to check only if the spectral radius of a linear
operator is less than unity, but the inductive argument requires more
elaborate notation than the first-order continuous differentiability.

\section{Application to periodic orbits of autonomous FDEs 
  --- Hopf Bifurcation Theorem}
\label{sec:po}
Let us come back to the original problem, the parameter-dependent FDE
\eqref{eq:ivp} $\dot x = f(x_t,\mu)$, where $\mu\in\R^\nu$ is a system
parameter and the functional $f:C^0(J;\R^n)\times\R^\nu\mapsto\R^n$ is
defined for first arguments that exist on an arbitrary compact
interval $J$.  Periodic orbits are solutions $x$ of $\dot x =
f(x_t,\mu)$ that are defined on $\R$ and satisfy $x(t)=x(t+T)$ for
some $T>0$ and all $t\in\R$.

Let $x$ be a periodic function of period $T=2\pi/\omega$. Then the
function $y(s)=x(s/\omega)$ is a function of period $2\pi$
($s\in\T$). This makes it useful to define the map
\begin{align*}
  S:&BC^0(\R;\R^n)\times\R\mapsto BC^0(\R;\R^n) &
  [S(y,\omega)](s)=y(\omega s)\mbox{,}
\end{align*}
such that $S(y,\omega)(t)=x(t)$ for all $t\in\R$ (remember that
$BC^0(\R;\R^n)$ is the space of bounded continuous functions on the
real line). Then $x\in C^1(\R;\R^n)$ satisfies the differential
equation
\begin{equation}
  \label{eq:ft}
  \dot x(t)=f(x_t,\mu)
\end{equation}
on the real line and has period $2\pi/\omega$ if and only if
$y=S(x,1/\omega)\in C^1(\T;\R^n)$ satisfies the differential equation
\begin{displaymath}
  \dot y(s)=\frac{1}{\omega}f(S(\Delta_sy,\omega),\mu)\mbox{.}
\end{displaymath}
Let us define an extended differential equation
\begin{align}
  \label{eq:ftext}
  \dot x_\mathrm{ext}(s)&=f_\mathrm{ext}(\Delta_sx_\mathrm{ext})\mbox{,}
\end{align}
where $f_\mathrm{ext}$ maps $C^0(\T;\R^{n+1+\nu})$ into
$\R^{n+1+\nu}$ and is defined by
\begin{align*}
  f_\mathrm{ext}
  \begin{pmatrix}
    y \\ \omega\\ \mu
  \end{pmatrix}&=
  \begin{bmatrix}
    f\left(S(y,R_0\omega),R_0\mu\right)/\cut(R_0\omega)\\ 0\\ 0
  \end{bmatrix}\mbox{,\quad where}\\
  \cut(\omega)&=
  \begin{cases}
    \omega &\mbox{if $\omega>\omega_\mathrm{cutoff}>0$}\\
    \mbox{smooth, uniformly non-negative extension} & \mbox{for
      $\omega<\omega_\mathrm{cutoff}$}
  \end{cases}
\end{align*}
for $y\in C^0(\T;\R^n)$, $\omega\in C^0(\T;\R)$ and $\mu\in
C^0(\T;\R^\nu)$ (recall that $R_0$ takes the average of a function on
$\T$).  We have used in our definition that any functional $f$ defined
for $x\in C^0(J;\R^n)$ is also a functional on $C^0(\T;\R^n)$
(periodic functions have a natural extension
$x(t)=x(t_{\mod[-\pi,\pi)})$ if $t\in \R$ is arbitrary). The extended
system has introduced the unknown $\omega$ and the system parameter
$\mu$ as functions of time, and the additional differential equations
$\dot \omega=0$, $\dot\mu=0$, which force the new functions to be
constant for solutions of \eqref{eq:ftext}.  We have also introduced a
cut-off for $\omega$ close to zero to keep $f_\mathrm{ext}$ globally
defined. The extended BVP \eqref{eq:ftext} is in the form of periodic
BVPs covered by  the Equivalence Theorem~\ref{thm:main}. Thus, if
$f_\mathrm{ext}$ is $EC^{j_{\max}}$ smooth then BVP \eqref{eq:ftext}
satisfies the assumptions of Theorem~\ref{thm:main} in the vicinity of
every periodic function $x_{0,\mathrm{ext}}\in
C^{1,1}(\T;\R^{n+\nu+1})$. Any solution
$x_\mathrm{ext}=(y,\omega,\mu)$ that we find for \eqref{eq:ftext}
corresponds to a periodic solution $t\mapsto y(\omega t)$ of period
$2\pi/R_0\omega$ at parameter $R_0\mu$ for \eqref{eq:ft} and vice
versa, as long as $R_0\omega>\omega_\mathrm{cutoff}$. The condition of
$EC^{j_{\max}}$ smoothness has to be checked only for the first $n$
components of the function $f_\mathrm{ext}$ since its final $\nu+1$
components are zero.

Application of the Equivalence Theorem~\ref{thm:main} results in a
system of algebraic equations that has $(n+\nu+1)(2N+1)$ variables and
equations, where $N$ is the positive integer proven to exist in
Theorem~\ref{thm:main}. Let us denote as $F=(F_y,F_\omega,F_\mu)$ the
components of the right-hand side $F_\mathrm{ext}$ (given by
$F_\mathrm{ext}(x_\mathrm{ext})(t)=f_\mathrm{ext}(\Delta_tx_\mathrm{ext})$),
of which $F_\mu$ and $F_\omega$ are identically zero. Let
$p=(p_y,p_\omega,p_\mu)$ be the $2N+1$ leading Fourier coefficients of
$y$, $\omega$ and $\mu$, respectively (these are the variables of the
algebraic system constructed via Theorem~\ref{thm:main}), and
$X(p)=(X_y(p),X_\omega(p),X_\mu(p))$ be the map from
$R^{(n+\nu+1)(2N+1)}$ into $C^{j_{\max}}(\T;\R^{n+\nu+1})$. Then
several of the components of $p$ can be eliminated as variables, and
the equations for $p$ resulting from Theorem~\ref{thm:main}
correspondingly simplified. Since $F$ is identically zero in its last
$\nu+1$ components we have
\begin{align*}
  X_\omega(p)&=E_Np_\omega\mbox{,} &
  X_\mu(p)&=E_Np_\mu\mbox{.}
\end{align*}
Hence, the right-hand side
\begin{displaymath}
  g(p)=R_N\left[P_0F(X(p))+Q_0\left(E_Np-P_NLF(X(p))\right)\right]\mbox{,}
\end{displaymath}
defined in Theorem~\ref{thm:main}, has $\nu+1$ components that are
identical to zero (since $P_0F(X(p))=0$ for the equations
$\dot\omega=0$ and $\dot\mu=0$). Furthermore, $g(p)=0$ contains the
equations $R_NQ_0E_Np_\omega=0$ and $R_NQ_0E_Np_\mu=0$, which require
that all Fourier coefficients (except the averages $R_0\omega$ and
$R_0\mu$) of $\mu$ and $\omega$ are equal to zero. This means
(unsurprisingly) that the algebraic system forces $\omega$ and $\mu$
to be constant. Thus, we can eliminate $R_NQ_0E_Np_\omega$ and
$R_NQ_oE_Np_\mu$ (which are $2N(\nu+1)$ variables), replacing them by
zero, and drop the corresponding equations. Since $\omega$ and $\mu$
must be constant, we can replace the arguments $p_\omega$ and $p_\mu$
of $X$ by the scalar $R_0E_Np_\omega$ (which we re-name back to
$\omega$) and the vector $R_0E_Np_\mu\in\R^\nu$ (which we re-name back
to $\mu$). This leaves the first $n(2N+1)$ algebraic equations
\begin{align}\label{eq:lowmodes:par}
    0=&R_N\left(P_0F_y(X_y(p_y,\omega,\mu),\omega,\mu)+
      Q_0\left[E_Np_y-
        P_NLF_y(X_y(p_y,\omega,\mu),\omega,\mu)\right]\right)\mbox{,}
\end{align}
which depend smoothly (with degree $j_{\max}$) on the $n(2N+1)$
variables $p_y$ and the parameters $\omega\in\R$ and
$\mu\in\R^\nu$. Overall, \eqref{eq:lowmodes:par} is a system of
$n\times(2\,N+1)$ equations.
\subsubsection*{Rotational Invariance}
The original nonlinearity $F$, defined by $[F(x)](t)=f(\Delta_tx)$  is
equivariant with respect to time shift: $\Delta_tF(x)=F(\Delta_tx)$ for
all $t\in\T$ and $x\in C^0(\T;\R^n)$. Furthermore, $\Delta_t$ commutes
with the following operations:
\begin{align*}
  \Delta_tQ_NL&=Q_NL\Delta_t \mbox{\quad (if $N\geq0$) and }
  &\Delta_tP_N&=P_N\Delta_t\mbox{.}
\end{align*}
This property gets passed on to the algebraic equation in the
following sense: let us define the operation $\Delta_t$ for a vector
$p$ in $\R^{n(2N+1)}$, which we consider as a vector of Fourier
coefficients of the function $E_Np\in C^0(\T;\R^n)$, by
\begin{displaymath}
  \Delta_tp=R_N\Delta_tE_Np\mbox{.}
\end{displaymath}
With this definition $\Delta_t$ commutes with $R_N$ and $E_N$. It is a
group of rotation matrices: $\Delta_t$ is regular for all $t$, and
$\Delta_{2k\pi}$ is the identity for all integers $k$.  The definition
of $X(p)$ as a fixed point of $x\mapsto E_Np+Q_NLF(x)$ implies that
$\Delta_tX(p)=X(\Delta_tp)$. From this it follows that the algebraic
system of equations is also equivariant with respect to $\Delta_t$. If
we denote the right-hand-side of the overall system
\eqref{eq:lowmodes:par} by $G(p_y,\omega,\mu)$ then $G$ satisfies
\begin{align*}
  \Delta_t G(p_y,\omega,\mu)=G(\Delta_tp_y,\omega,\mu)\mbox{\quad for
    all $t\in\T$, $p_y\in\R^{n(2N+1)}$, $\omega>0$ and
    $\mu\in\R^\nu$.}
\end{align*}

\subsubsection*{Application to Hopf bifurcation}

One useful aspect of the Equivalence Theorem is that it provides an
alternative approach to proving the Hopf Bifurcation Theorem for
equations with state-dependent delays. The first proof that the Hopf
bifurcation occurs as expected is due to Eichmann \cite{E06}. The
reduction of periodic boundary-value problems to smooth algebraic
equations reduces the Hopf bifurcation problem to an equivariant
algebraic pitchfork bifurcation.

Let us consider the equation
\begin{equation}
  \label{eq:dynsys}
  0=f(E_0x_0,\mu)
\end{equation}
where $f:C^0(J;\R^n)\times\R\mapsto\R^n$, $\mu\in\R$, $J$ is a compact
interval, $x_0\in\R^n$, and the operator $E_0$ (as defined in
\eqref{eq:proj} in Section~\ref{sec:theorem}) extends a constant to a
function on $\T$ (and thus, on $J$).  This means that
\eqref{eq:dynsys} is a system of $n$ algebraic equations for the $n+1$
variables $(x_0,\mu)$. The definition of $EC^k$ smoothness does not
cover functionals that depend on parameters. We avoid the introduction
of a separate definition of $EC^k$ smoothness for parameter-dependent
functionals that distinguishes between parameters and functional
arguments. We rather extend Definition~\ref{def:extdiff}:
\begin{definition}[$EC^k$ smoothness for parameter-dependent functionals]
  \label{def:extdiff:par} Let $J=[a,b]$ be a compact interval (or $J=\T$),
  and $D$ be a product space of the form
  $D=C^{k_1}(J;\R^{m_1})\times\ldots\times C^{k_\ell}(J;\R^{m_\ell})$
  where $\ell\geq 1$, and $k_j\geq0$ and $m_j\geq1$ are integers. We
  say that $f:D\times\R^\nu\mapsto \R^n$ is $EC^k$ smooth if the
  functional
  \begin{equation}\label{eq:ass:pareck}
    (x,y)\in D\times C^0(J;\R^\nu)\mapsto f(x,y(a))\in\R^n
  \end{equation}
  is $EC^k$ smooth (if $J=\T$ we use $a=-\pi$).
\end{definition}
Requiring $EC^k$-smoothness of the parameter-dependent functional $f$
in this sense, implies that the algebraic system $0=f(E_0x_0,\mu)$ is
$k$ times continuously differentiable. Let us assume that the
algebraic system $0=f(E_0x_0,\mu)$ has a regular solution
$x_0(\mu)\in\R^n$ for $\mu$ close to $0$. Without loss of generality
we can assume that $x_0(\mu)=0$, otherwise, we introduce the new
variable $x_\mathrm{new}=x_\mathrm{old}-E_0x_0(\mu)\in
C^0(J;\R^n)$. Hence, $f(0,\mu)=0$ for all $\mu$ close to $0$.

The $EC^1$ derivative of $f$ in $(0,\mu)$ is a linear functional,
mapping $C^0(J;\R^{n+\nu})$ into $\R^n$. Let us denote its first $n$
components (the derivative with respect to the first argument $x$ of
$f$) by $a(\mu)$. The linear operator $a(\mu)$ can easily be
complexified by defining $a(\mu)[x+iy]=a(\mu)[x]+ia(\mu)[y]$ for
$x+iy\in C^0([-\tau,0];\C^n)$. If $f$ is $EC^k$ smooth with $k\geq2$
then the $n\times n$-matrix $K(\lambda,\mu)$ (called the
\emph{characteristic matrix}), defined by
\begin{equation}\label{eq:charmdef}
  K(\lambda,\mu)\,v=\lambda v-a(\mu)[v\exp(\lambda t)]
\end{equation}
is analytic in its complex argument $\lambda$ and $k-1$ times
differentiable in its real argument $\mu$ (since the functions
$t\mapsto v\exp(\lambda t)$ to which $a(\mu)$ is applied are all
elements of $C^k(J;\R^n)$).  The Hopf Bifurcation Theorem states the
following:
\begin{theorem}[Hopf bifurcation]\label{thm:hopf}
  Assume that $f$ is $EC^k$ smooth ($k\geq 2$) in the sense of
  Definition~\ref{def:extdiff:par} and that the characteristic matrix
  $K(\lambda,\mu)$ satisfies the following conditions:
  \begin{enumerate}
  \item \label{ass:hopf:imag}\textbf{(Imaginary eigenvalue)} there
    exists an $\omega_0>0$ such that $\det K(i\omega_0,0)=0$ and
    $i\omega_0$ is an isolated root of $\lambda\mapsto \det
    K(\lambda,0)$. We denote the corresponding null vector by
    $v_1=v_r+iv_i\in\C^n$ (scaling it such that $|v_r|^2+|v_i|^2=1$).
  \item \label{ass:hopf:nonresonant} \textbf{(Non-resonance)} $\det
    K(ik\omega,0)\neq 0$ for all integers $k\neq \pm1$.
  \item \label{ass:hopf:cross} \textbf{(Transversal crossing)} The
    local root curve $\mu\mapsto \lambda(\mu)$ of $\det
    K(\lambda,\mu)$ that corresponds to the isolated root $i\omega_0$
    at $\mu=0$ (that is, $\lambda(0)=i\omega_0$) has a non-vanishing
    derivative of its real part:
    \begin{displaymath}
      0\neq\left.\frac{\partial}{\partial\mu}
      \Re\lambda(\mu)\right\vert_{\textstyle\mu=0}\mbox{.}
    \end{displaymath}
  \end{enumerate}
  Then there exists a $k-1$ times differentiable curve
  \begin{displaymath}
    \beta\in(-\epsilon,\epsilon)\mapsto (x,\omega,\mu)\in 
    C^1(\T;\R^n)\times\R\times\R
  \end{displaymath}
  such that for sufficiently small $\epsilon>0$ the following holds:
  \begin{enumerate}
  \item\label{thm:hopf:isperiodic} $x(\omega\cdot)$ (or $S(x,\omega)$)
    is a periodic orbit of $\dot x(t)=f(x_t,\mu)$ of period
    $2\pi/\omega$, that is, $x\in C^1(\T;\R^n)$ and
    \begin{equation}
      \dot x(t)=
      \frac{1}{\omega}f(S(\Delta_tx,\omega),\mu)\mbox{,}\label{eq:hopf:po}    
  \end{equation}
\item\label{thm:hopf:phaseamp} the first Fourier coefficients of $x$
  are equal to $(0,\beta)$, that is,
    \begin{equation}
      \begin{split}
        0&=\frac{1}{\pi}\int_{-\pi}^\pi
        \Re\left[v_1\exp(it)\right]^T x(t)\d t\mbox{,\quad and}\\
        \beta&=-\frac{1}{\pi}
        \int_{-\pi}^\pi\Im\left[v_1\exp(it)\right]^T x(t)\d t\mbox{,}
      \end{split}
      \label{eq:hopfphase}
    \end{equation}
  \item\label{thm:hopf:pars} $x\vert_{\beta=0}=0$,
    $\mu\vert_{\beta=0}=0$ and $\omega\vert_{\beta=0}=\omega_0$, that
    is, the solution $x$, the system parameter $\mu$ and the frequency
    $\omega$ of $x$, which are differentiable functions of the
    amplitude $\beta$, are equal to $x=0$, $\mu=0$, $\omega=\omega_0$
    for $\beta=0$.
  \end{enumerate}
\end{theorem}
The statement is identical to the classical Hopf Bifurcation Theorem
for ODEs in its assumptions and conclusions apart from the regularity
assumption on $f$ specific to FDEs. Note that the existence of the
one-parameter family (parametrized in $\beta$) automatically implies
the existence of a two-parameter family for $\beta\neq0$ due to the
rotational invariance: if $x$ is a solution of \eqref{eq:hopf:po} then
$\Delta_sx$ is also a solution of \eqref{eq:hopf:po} for every fixed
$s\in\T$. Condition~\eqref{eq:hopfphase} fixes the time shift $s$ of
$x$ such that $x$ is orthogonal to $\Re[v_1\exp(it)]$ using the
$\Lint^2$ scalar product on $\T$.

The proof of the Hopf Bifurcation Theorem is a simple fact-checking
exercise. We have to translate the assumptions on the derivative of
$f:C^0(J;\R^n)\mapsto\R^n$ into properties of the right-hand side of
the nonlinear algebraic system \eqref{eq:lowmodes:par} near
$(x,\omega,\mu)=(0,\omega_0,0)$, and then apply algebraic bifurcation
theory to the algebraic system. The only element of the proof that is
specific to functional differential equations comes in at the linear
level: the fact that the eigenvalue $i\omega_0$ is simple implies that
the right nullvector $v_1\in\C^n$ and any non-trivial left nullvector
$w_1$ satisfy
\begin{align*}
  w_1^H\left[\left.\frac{\partial}{\partial
      \lambda}K(\lambda,0)\right\vert_{
      \textstyle\lambda=i\omega_0}\right]
  v_1\neq
  0\mbox{.}
\end{align*}
This is the generalization of the orthogonality condition
$w_1^Hv_1\neq 0$, known from ordinary matrix eigenvalue problems, to
exponential matrix eigenvalue problems of the type
$K(\lambda,\mu)\,v=0$. The proof of Theorem~\ref{thm:hopf} is
entirely based on the standard calculus arguments for branching of
solutions to algebraic systems as can be found in textbooks
\cite{AG79}. The details of the proof can be found in
Section~\ref{sec:hopf:proof}.

The statements in Theorem~\ref{thm:hopf} should be compared to two
previous works considering the same situation: the branching of
periodic orbits from an equilibrium losing its stability in FDEs with
state-dependent delays. The PhD thesis of Eichmann (2006, \cite{E06})
proves the existence of the curve
$\beta\mapsto(x(\beta),\mu(\beta),\omega(\beta))$ and that it is once
continuously differentiable (assuming only $EC^2$ smoothness of
$f$). Since $\mu'(0)=0$ due to rotational symmetry (see proof in
Section~\ref{sec:hopf:proof}) this is not enough to determine if the
non-trivial periodic solutions exist for $\mu>0$ or for $\mu<0$ for
small $\beta$ (the so-called \emph{criticality} of the Hopf
bifurcation, which is of interest in applications \cite{IBS08,IST05}).
Moreover, the non-resonance condition in \cite{E06} is slightly too
strong, requiring that $i\omega_0$ is the only purely imaginary root
of $\det K(\lambda,0)$ (this assumption is different in the summary
given in the review by Hartung \emph{et al.} \cite{HKWW06}; note that
the publicly available version of \cite{E06} has a typo in the
corresponding assumption L1), and only the pure-delay case (where the
time interval $J$ equals $[-\tau,0]$) was considered. However, the
techniques employed in \cite{E06}, based on the Fredholm alternative,
are likely to yield exactly the same result as stated in
Theorem~\ref{thm:hopf} if one assumes general $EC^k$ smoothness with
$k\geq2$ (the formulation of $EC^2$ smoothness is already rather
convoluted in \cite{E06}).

Hu and Wu \cite{HW10} use $S^1$-degree theory \cite{EWG92,KW97} to
prove the existence of a branch of non-trivial periodic solutions near
$(x,\mu,\omega)=(0,0,\omega_0)$. This type of topological methods
gives generally weaker results concerning the local uniqueness of
branches of periodic orbits or their regularity, but they require only
weaker assumptions (\cite{HW10} still needs to assume $EC^2$
smoothness, though). Degree methods also give global existence results
by placing restrictions on the number of branches that can occur.


\section{Conclusion, applications and generalizations}
\label{sec:conc}
Periodic boundary-value problems for functional differential equations
(FDEs) are equivalent to systems of smooth algebraic equations if the
functional $f$ defining the right-hand side of the boundary-value
problem satisfies natural smoothness assumptions. These assumptions
are identical to those imposed in the review by Hartung \emph{et al.}
\cite{HKWW06} and do not exclude FDEs with state-dependent
delay. There are several immediate extensions of the results presented
in this paper. The list below indicates some of them.

\subsubsection*{Further potential applications of the Equivalence
  Theorem~\ref{thm:main}}
Theorem~\ref{thm:hopf} on the Hopf bifurcation is not the central
result of the paper, even though it is a moderate extension of the
theorem proved in \cite{E06}. Rather, it is a demonstration of the use
of the Equivalence Theorem~\ref{thm:main}. The main strength of the
equivalence result stated in Theorem~\ref{thm:main} is that it permits
the straightforward application of continuation and Lyapunov-Schmidt
reduction techniques to FDE problems involving periodic orbits of
finite period, regardless if the delay is state dependent, or if the
equation is of so-called mixed type (that is, positive and negative
delays are present). A source of complexity, for example, in
\cite{E06,EWG92,G11,HW10,KW97,W98}, is that techniques such as
Lyapunov-Schmidt reduction or $S^1$-degree theory had to be applied in
Banach spaces. Theorem~\ref{thm:main} removes the need for this,
reducing the analysis of periodic orbits to root-finding in
$\R^{n\times(2\,N+1)}$. For example, Humphries \emph{et al.}
\cite{HDMU12} study periodic orbits in FDEs with two state-dependent
delays numerically using DDE-Biftool \cite{ELR02}, alluding to
theoretical results about bifurcations of periodic orbits that have
been proven only for constant delay. Fist of all, Humphries \emph{et
  al.} \cite{HDMU12} continue branches of periodic orbits.
Theorem~\ref{thm:main} makes clear when these branches as curves of
points $(x,\omega,\mu)$ in the extended space
$C^0(\T;\R^n)\times\R\times\R$ are smooth to arbitrary degree: the
Jacobian of \eqref{eq:lowmodes:par} with respect to $(p_y,\omega,\mu)$
along the curve has to have full rank. Along these branches Humphries
\emph{et al.} \cite{HDMU12} encounter degeneracies of the linearization
and conjecture the existence of the corresponding bifurcations (backed
up with numerical evidence) such as: fold bifurcations,
period-doubling bifurcations, or the branching off of resonance
surfaces (Arnol'd tongues) when resonant Floquet multipliers of the
linearized equation cross the unit circle \cite{K04}.  The Equivalence
Theorem~\ref{thm:main} provides a straightforward route to proofs
that these scenarios occur as expected. Similarly,
Theorem~\ref{thm:main} will likely not only simplify the proofs about
bifurcations of symmetric periodic orbits such as those of Wu
\cite{W98}, but also extend them to the case of FDEs with
state-dependent delays. As long as one considers branching of periodic
orbits with finite periods, the problem can be reduced locally (and
often in every ball of finite size) to a finite-dimensional
root-finding problem. This transfers also a list of results of
symmetric bifurcation theory found in textbooks \cite{GSS88} to FDEs
with state-dependent delays.

\subsubsection*{Globally valid algebraic system} The main result was
formulated locally in the neighborhood of a given $x_0\in
C^{1,1}(\T;\R^n)$ and required only local Lipschitz continuity. The
proof makes obvious that the domain of definition for the map $X$,
which maps between the function space and the finite-dimensional space
is limited by the size of the neighborhood of $x_0$ for which one can
find a uniform ($EC$) Lipschitz constant of the right-hand side
$F$. In problems with delay the right-hand side is typically a
combination of Nemytskii operators generated by smooth functions and
the evaluation operator $\ev:C^0(\T;\R^n)\times\T\mapsto\R^n$, given by
$\ev(x,t)=x(t)$. These typically satisfy a (semi-)global Lipschitz
condition (see also condition (Lb) in \cite{HKWW06}): for all $R$
there exists a constant $K$ such that
\begin{displaymath}
  |f(x)-f(y)|\leq K\|x-y\|_0
\end{displaymath}
for all $x$ and $y$ satisfying $\|x\|_1\leq R$ and $\|y\|_1\leq R$.
Under this condition one can choose for any bounded ball an algebraic
system that is equivalent to the periodic boundary-value problem in
this bounded ball. For periodic orbits of autonomous systems this
means that one can find an algebraic system such that all periodic
orbits of amplitude less than $R$ and of period and frequency at most
$R$ are given by the roots of the algebraic system.

\subsubsection*{Implicitly given delays}
In practical applications the delay is sometimes given implicitly, for
example, the position control problem considered in \cite{W02} and the
cutting problem in \cite{IBS08,IST05} contain a separate algebraic
equation, which defines the delay implicitly. In simple cases these
problems can be reduced to explicit differential equations using the
standard Baumgarte reduction \cite{B72} for index-1 differential
algebraic equations. For example, in the cutting problem the delay
$\tau$ depends on the current position $x$ via the implicit linear
equation
\begin{equation}\label{eq:cutdelay}
  \tau(t)=a-bx(t)-bx(t-\tau(t))\mbox{,}
\end{equation}
which can be transformed into a differential equation by
differentiation with respect to time:
\begin{equation}\label{eq:cutdelaydiff}
  \dot\tau(t)=
  \frac{-bv(t)-\left[\tau(t)-a+bx(t)+bx(t-\tau(t))\right]}{1+bv(t-\tau(t))}
\end{equation}
($v(t)=\dot x(t)$ is explicitly present as a variable in the cutting
model, which is a second-order differential equation). The original
model accompanied with the differential equation
\eqref{eq:cutdelaydiff} instead of the algebraic equation
\eqref{eq:cutdelay} fits into the conditions of the Equivalence
Theorem~\ref{thm:main}. The regularity statement of the Equivalence
Theorem guarantees that the resulting periodic solutions have
Lipschitz continuous derivatives with respect to time. This implies
that the defect $d=\tau-(a-bx-bx(t-\tau))$ occurring in the algebraic
condition~\eqref{eq:cutdelay} satisfies $\dot d(t)=-d(t)$ along
solutions. Since the solutions are periodic the defect $d$ is
periodic, too, and, hence, $d$ is identically zero. The denominator
appearing in Equation~\eqref{eq:cutdelaydiff} becomes zero exactly in
those points in which the implicit condition \eqref{eq:cutdelay}
cannot be solved for the delay $\tau$ with a regular derivative.

The same argument can be applied to the position control problem as
long as the object, at position $x$, does not hit the base at position
$-w$ (the model contains the term $|x+w|$ in the right-hand side).

\subsubsection*{Neutral equations} The index reduction works only if the delay
$\tau$, which is itself a function of time, is not evaluated at
different time points. For example, changing $bx(t-\tau(t))$ to
$bx(t-\tau(t-1))$ on the right-hand-side of \eqref{eq:cutdelay} would
make the index reduction impossible. However, certain simple neutral
equations permit a similar reduction directly on the function space
level. Consider
\begin{equation}
  \label{eq:neutral}
  \frac{\d}{\d t} \left[\Delta_t(x+g(x))\right]=f(\Delta_t x)
\end{equation}
where the functional $f$ satisfies the local $EC$ Lipschitz condition,
defined in Definition~\ref{def:loclip}, in a neighborhood $U$ of a
point $x_0\in\C^{1,1}(\T;\R^n)$, and $g:C^0(\T;\R^n)\mapsto\R^n$ has a
global (classical) Lipschitz constant less than unity (this excludes
state-dependent delays in the essential part of the neutral
equation). Then one can define the map $X_g(y)$ as the unique solution
$x$ of the fixed point problem
\begin{displaymath}
  x(t)=y(t)-g(\Delta_tx)\mbox{\quad for $y$ near $y_0=x_0+g(x_0)$,}
\end{displaymath}
which reduces \eqref{eq:neutral} to the equation
\begin{equation}\label{eq:neutralred}
  \dot y(t)=f(\Delta_tX_g(y))=f(X_g(\Delta_ty))\mbox{.}
\end{equation}
Equation~\eqref{eq:neutralred} satisfies the conditions of the
Equivalence Theorem~\ref{thm:main}. One implication of this reduction
is that periodic solutions of \eqref{eq:neutral} are $k$ times
continuously differentiable if the functional $f$ is $EC^k$ smooth in
the sense of Definition~\ref{def:extdiff} and $g$ is $k$ times
continuously differentiable as a map from $C^0(\T;\R^n)$ into $\R^n$.

\section{Proof of the Equivalence Theorem~\ref{thm:main}}
\label{sec:remproofs}
Theorem~\ref{thm:main} is proved in three steps. First, we establish
the existence of a locally unique fixed point of the map
$M_N(\cdot,p)$ using Banach's contraction mapping principle. This step
requires only local $EC$ Lipschitz continuity in the sense of
Definition~\ref{def:loclip}. In the second step we prove continuous
differentiability of the map $X$ and the right-hand side $g$ of the
algebraic system assuming that $f$ is $EC^1$ smooth. In the final step
we prove higher-order differentiability, assuming that $f$ is $EC^k$
smooth for degrees $k$ up to $j_{\max}$.

\subsection{Decay of Fourier coefficients for integrals and smooth functions}
\label{sec:qnl}
The following preparatory lemma states the well-known fact that,
colloquially speaking, integrating a function makes its high-frequency
Fourier coefficients smaller. In the fixed-point
equation~\eqref{eq:fixp:intro} of Theorem~\ref{thm:main} the term
$Q_NL$ occurs, and we need this term to be small for large $N$. Recall
that $Q_N$ removes the first $N$ Fourier modes from a periodic
function and $Lx$ is the anti-derivative of $x$ (after subtracting the
average of $x$), see Equation~\eqref{eq:proj} for the precise
definitions.
\begin{lemma}[Decay of Fourier coefficients of integrals]\label{thm:qnl}
  The norm of the linear operator $Q_NL$, mapping the space
  $C^j(\T;\R^n)$ back into itself, is bounded by
  \begin{displaymath}
    \|Q_NL\|_j\leq C\frac{\log N}{N}
  \end{displaymath}
  where $C$ is a constant. The same holds in the Lipschitz norm (with
  the same constant $C$):
  \begin{displaymath}
    \|Q_NL\|_{0,1}\leq C\frac{\log N}{N}\mbox{.}
  \end{displaymath}
\end{lemma}
\pstart{Proof} 
We find the norm $\|Q_NL\|_0$ first, and start out with the well-known
estimate for interpolating trigonometric polynomials for continuous
functions on $\T$. Let $x$ be a continuous function on $\T$ with
modulus of continuity $\omega$. Then (see \cite{J30})
\begin{displaymath}
  \|Q_Nx\|_0\leq C_0\omega\left(\frac{2\pi}{N}\right)\log N
\end{displaymath}
where $C_0$ is a constant that does not depend on $x$ or $N$.  A
function $\omega:[0,\infty)\mapsto[0,\infty)$ is called a modulus of
continuity of a continuous function $x$ if
\begin{displaymath}
  |x(t)-x(s)|\leq \omega(|t-s|)\mbox{.}
\end{displaymath}
holds for all $s$ and $t\in\T$.  For a function $x\in C^0(\T;\R^n)$ the
anti-derivative
\begin{displaymath}
  [Lx](t)=\int_0^tx(s)-R_0x\,\d s
\end{displaymath}
has the Lipschitz constant $\|x\|_0=\max\{|x(t)|:t\in\T\}$ such that a
modulus of continuity for $Lx$ is $\omega(h)=\|x\|_0
h$. Consequently,
\begin{equation}\label{eq:qnlf0}
  \|Q_NLx\|_0\leq C_0\frac{2\pi\|x\|_0}{N}\log N\mbox{,}
\end{equation}
where $C_0$ does not depend on $x$ or $N$. This proves the claim of the
lemma for $j=0$. For $x\in C^j(\T;\R^n)$ we notice that all
derivatives of $x$ up to order $j$ are continuous. Applying estimate
\eqref{eq:qnlf0} to each of the derivatives of $x$ we get
\begin{displaymath}
  \|Q_NLx^{(l)}\|_0\leq \frac{2\pi C_0}{N}\log N\,\|x^{(l)}\|_0
  \mbox{\quad for $l=0,\ldots j$.}  
\end{displaymath}
Consequently, the maximum of the left-hand sides over all
$l\in\{0,\ldots,j\}$ must be less than the maximum of the right-hand
sides:
\begin{displaymath}
  \|Q_NLx\|_j=\max_{l=0,\ldots,j}\|Q_NLx^{(l)}\|_0\leq   
  \frac{2\pi C_0}{N}\log N\,\max_{l=0,\ldots,j}\|x^{(l)}\|_0=
  2\pi C_0\frac{\log N}{N}\|x\|_j\mbox{,}
\end{displaymath}
which implies the desired estimate for $\|Q_NL\|_j$ using the constant
$C=2\pi C_0$.  

The estimate of $Q_NL$ in the Lipschitz norm is a continuity
argument. The operator $Q_NL$ is bounded (and, thus, continuous) on
$C^{0,1}(\T;\R^n)$. For every element $y$ of $C^1(\T;\R^n)$ (which is
a dense subspace of $C^{0,1}(\T;\R^n)$) the Lipschitz constant is
identical to $\|y'\|_0=\max_{t\in\T}|y'(t)|$, and, thus,
$\|y\|_1=\|y\|_{0,1}$. Let $x_n\in C^1(\T;\R^n)$ be a sequence of
continuously differentiable functions that converges to $x\in
C^{0,1}(\T;\R^n)$ in the $\|\cdot\|_{0,1}$-norm: $\|x_n-x\|_{0,1}\to0$
for $n\to\infty$. Then
\begin{displaymath}
  \|Q_NLx_n\|_{0,1}=\|Q_NLx_n\|_1\leq C\frac{\log N}{N}\|x_n\|_1=
  C\frac{\log N}{N}\|x_n\|_{0,1}\mbox{.}
\end{displaymath}
On both sides of the inequality the limit for $n\to\infty$ exists,
resulting in the desired estimate for $\|Q_NL\|_{0,1}$.  \eop

A direct consequence of Lemma~\ref{thm:qnl} is that the Lipschitz norm
of $Q_Nx$, $\|Q_Nx\|_{0,1}$, goes to zero for $N\to\infty$ for
elements of $C^{1,1}(\T;\R^n)$, so, for example, for a solution $x$ of
a periodic BVP:
\begin{equation}\label{eq:qn}
  \|Q_Nx\|_{0,1}=\|Q_NLx'\|_{0,1}\leq 
  C\frac{\log N}{N}\|x'\|_{0,1}\leq C\frac{\log N}{N}\|x\|_{1,1}\mbox{.}
\end{equation}

\subsection{Proof of Splitting Lemma~\ref{thm:split}}
\label{sec:splitproof}
For any given integer $N\geq0$ we have to show that the pair $(x,p)\in
C^0(\T;\R^n)\times\R^{n\times(2\,N+1)}$ satisfies
\begin{align}\label{eq:inteq:app}
  x(t)&=x(0)+\int_0^tF(x)(s)\d s\mbox{\quad for all $t\in\T$, and}\\
  p&=R_Nx\mbox{,\quad (or, equivalently, $E_Np=P_Nx$)}\label{eq:peqpnx:app}\\
  \intertext{ if and only if it satisfies the system}
  x&=E_Np+Q_NLF(x)\mbox{,}\label{eq:fixp:app}\\
  0&=R_N\left[P_0F(x)+Q_0\left(E_Np-P_NLF(x)\right)\right]
  \label{eq:lowmodes:app}\mbox{,}  
\end{align}
``$\Rightarrow$'': Assume that $x\in C^0(\T;\R^n)$ satisfies
\eqref{eq:inteq:app}, and let $p=R_Nx$.  Subtracting
equation~\eqref{eq:inteq:app} for $t=-\pi$ from \eqref{eq:inteq:app}
for $t=\pi$ implies that the average of $F(x)$ is zero. Thus,
$R_0F(x)=0$ and $P_0F(x)=0$.  Since $Ly=\int_0^ty(s)-R_0y\d s$, the
identity \eqref{eq:inteq:app} implies (in combination with
$R_0F(x)=0$)
\begin{equation}\label{eq:intident:app}
  x(t)=x(0)+[LF(x)](t)\mbox{.}
\end{equation}
Applying projection $Q_N$ to this identity we obtain
$Q_Nx=Q_NLF(x)$. Adding \eqref{eq:peqpnx:app} to this we obtain
equation~\eqref{eq:fixp:app}. Applying projection $P_NQ_0$ (which is
the same as $Q_0P_N$) to \eqref{eq:intident:app} we obtain
$Q_0P_Nx=Q_0P_NLF(x)$. Inserting $E_Np$ for $P_Nx$ into this identity
leads to $Q_0[E_Np-P_NLF(x)]=0$. Since $P_0F(x)=0$, this in turn
implies \eqref{eq:lowmodes:app}.

``$\Leftarrow$'': Applying $P_N$ to \eqref{eq:fixp:app} implies
$P_Nx=E_Np$ (and $p=R_Nx$) immediately.  The expression inside the
parentheses of $R_N$ in equation \eqref{eq:lowmodes:app} is a sum of
two parts that each have to be zero (since they are both in the image
of $P_N$ on which $R_N$ in injective). The projection $Q_0$ subtracts
the average from its argument. Hence, $(x,p)\in
C^0(\T;\R^n)\times\R^{n\times(2\,N+1)}$ satisfies
\eqref{eq:fixp:app}--\eqref{eq:lowmodes:app} if and only if there
exists a constant $c\in\R^n$ such that the triple $(c,x,p)$ satisfies
the system of equations consisting of \eqref{eq:fixp:app} and 
\begin{align}
  0&=R_0F(x)\label{eq:avg:app}\\
  E_0c&=E_Np-P_NLF(x)\label{eq:othermodes:app}\mbox{.}
\end{align}
Note that $E_0$ maps the constant $c\in\R^n$ to a function that equals
this constant for all $t\in\T$. In this system, \eqref{eq:avg:app}
ensures that the average of $F(x)$ is
zero. Equation~\eqref{eq:othermodes:app} is an equation in the
finite-dimensional space $\rg P_N$.  Subtracting
\eqref{eq:othermodes:app} from \eqref{eq:fixp:app} gives
\begin{displaymath}
  x=E_0c+LF(x)\mbox{.}
\end{displaymath}
This equals \eqref{eq:inteq:app}, keeping in mind that
$[Ly](t)=\int_0^ty(s)-R_0y\d s$ ($[Ly](0)=0$ for all $y\in
C^0(\T;\R^n)$, hence, $x(0)=c$), and using $R_0F(x)=0$ (see
equation~\eqref{eq:avg:app}).\eop

\subsection{Unique solvability of the fixed point problem \eqref{eq:fixp:intro}}
\label{sec:fixproof}
Let $x_0$ be an element of $C^{1,1}(\T;\R^n)$, for example, a solution
of the periodic boundary value problem $\dot
x(t)=f(\Delta_tx)=F(x)(t)$. Consider a closed ball
$B_\delta^{0,1}(x_0)$ of radius $\delta$ around $x_0$ in the Lipschitz
norm:
\begin{displaymath}
  B_\delta^{0,1}(x_0)=\{x\in C^{0,1}(\T;\R^n):
\|x-x_0\|_{0,1}\leq \delta\}\mbox{.}
\end{displaymath}
The superscript ``$0,1$'' indicates which norm is used to measure the
distance from $x_0$ and that only elements of $C^{0,1}(\T;\R^n)$ are
included.  

Lemma~\ref{thm:loclip} implies that $F$ is Lipschitz continuous with
respect to the $\|\cdot\|_0$-norm in $B_\delta^{0,1}(x_0)$ if we
choose $\delta$ sufficiently small (thus, $F$ is also called $EC$
Lipschitz continuous in $B_\delta^{0,1}(x_0)$):
\begin{equation}\label{eq:flip0}
  \|F(x)-F(y)\|_0\leq K\|x-y\|_0
\end{equation}
for all $x$ and $y$ in $B_\delta^{0,1}(x_0)$ and a fixed $K>0$. In any
ball $B_\delta^{0,1}$, in which $F$ is $EC$ Lipschitz continuous, $F$
is also bounded in the Lipschitz norm:
\begin{equation}\label{eq:Fboundlip}
  \|F(x)\|_{0,1}\leq R\mbox{\quad  for all $x\in B_\delta^{0,1}(x_0)$.}    
\end{equation}
See Lemma~\ref{thm:Flipbound} in Appendix~\ref{sec:basicprop} for the
proof.

We can now formulate a lemma about the unique solvability of the fixed
point problem
\begin{displaymath}
  x=E_Np+Q_NLF(x)\mbox{.}
\end{displaymath}
This unique solvability and the Splitting Lemma~\ref{thm:split} allow
us to reduce the periodic BVP $\dot x(t)=f(\Delta_tx)$ to a system of
algebraic equations. Remember that $E_Np$ takes a vector $p$ of $2N+1$
Fourier coefficients and maps it to the periodic function having these
Fourier coefficients, $R_Nx$ extracts the first $2N+1$ Fourier
coefficients from a periodic function $x$, $P_Nx$ projects the
periodic function $x$ onto the space spanned by the basis
$b_{-N},\ldots,b_N$ and $Q_N=\id-P_N$ sets the first Fourier modes of
a function to zero. ($P_N$ and $Q_N$ are projections in the function
space, and $R_N$ and $E_N$ map between the finite-dimensional subspace
$\rg P_N$ and $\R^{n\times(2N+1)}$.)
\begin{lemma}[Unique solvability of fixed point problem]\label{thm:fixp}
  Let $x_0$ be in $C^{1,1}(\T;\R^n)$, and let $\delta>0$ be such that
  \begin{equation}\label{eq:fbounds}
  \|F(x)\|_{0,1}\leq R\mbox{\quad and\quad} \|F(x)-F(y)\|_0\leq K\|x-y\|_0  
  \end{equation}
  for all $x$ and $y\in B_{6\delta}^{0,1}(x_0)$ and for some constants
  $K>0$ and $R>0$ depending on $\delta$. Then for any sufficiently
  large $N$ the fixed point problem
  \begin{equation}\label{eq:fixp}
    x=E_Np+Q_NLF(x)
  \end{equation}
  has a unique solution $x\in B_{6\delta}^{0,1}(x_0)$ for all vectors
  $p\in\R^{n\times(2N+1)}$ in the neighborhood $U$ of $R_Nx_0$ given
  by
  \begin{equation}
    \label{eq:pclose}
        U=\left\{p\in\R^{n\times(2\,N+1)}: 
      \left\|E_N\left[p-R_Nx_0\right]\right\|_{0,1}<2\delta\right\}\mbox{.}
  \end{equation}
  Moreover, if $x\in B_\delta^{0,1}(x_0)$ is continuously
  differentiable and satisfies $x'=F(x)$ then its projection $p=R_Nx$ is
  in the neighborhood $U$, and $x$ and $p$ satisfy \eqref{eq:fixp}.
\end{lemma}
Note that $U$ is an open set of $\R^{n\times(2N+1)}$ since $E_N$ is an
isomorphism between $\rg P_N$, equipped with the
$\|\cdot\|_{0,1}$-norm, and $\R^{n\times(2N+1)}$. We have to prove the
unique solvability of the fixed-point problem in a slightly larger
ball (radius $6\delta$) and for a slightly larger range of parameters
$p$ (note the $2\delta$ in \eqref{eq:pclose}) in order to establish
one-to-one correspondence in the ball of radius $\delta$.


\pstart{Proof}
The idea is, of course, that the function
\begin{align*}
  M_N(\cdot,p): x\mapsto E_Np+Q_NLF(x)
\end{align*}
maps the closed ball $B_{6\delta}^{0,1}(x_0)$ back into itself and is
uniformly contracting for suitably large $N$ and vectors $p\in U$.

First, any closed ball $B_r^{0,1}(x_0)$ is closed (and, thus, forms a
complete metric space) with respect to the $\|\cdot\|_0$-norm. This
completeness is a simple continuity argument: let $y_n=x_0+z_n$ be a
Cauchy sequence in $B_r^{0,1}(x_0)$ with respect to the
$\|\cdot\|_0$-norm. Then $z_n$ converges to a continuous function $z$,
and, since $\|z_n\|_0\leq\|z_n\|_{0,1}\leq r$, for all $n$, the
maximum norm of $z$ is also bounded by $r$: $\|z\|_0\leq r$. We only
have to show that the Lipschitz constant of $z$ is bounded by
$r$, too. Let $\epsilon>0$ be arbitrary and let $t\neq s$ be
arbitrary in $\T$. We select some $n$ such that $\|z-z_n\|_0<
\epsilon|t-s|/2$. Then
\begin{align*}
  |z(t)-z(s)|&\leq|z(t)-z_n(t)|+|z_n(t)-z_n(s)|+|z_n(s)-z(s)|\\
  &< \epsilon|t-s|+r|t-s|\leq (r+\epsilon)|t-s|\mbox{.}
\end{align*}
Thus, the Lipschitz constant of $z$ is less than $r+\epsilon$ for
arbitrary $\epsilon>0$. Hence, $\|z\|_{0,1}\leq r$, completing the
argument for completeness of $B_r^{0,1}(x_0)$ with respect to the
$\|\cdot\|_0$-norm. This completeness implies that we can apply
Banach's contraction mapping principle in a ball $B_r^{0,1}(x_0)$, a
ball of Lipschitz continuous functions, using the weaker maximum norm
in the following.

We choose the radius $r$ of the ball equal to $6\delta$ ($\delta$ was
chosen in the lemma such that the estimates \eqref{eq:fbounds} are
true for the constants $K$ and $R$), Thus, $B_{6\delta}^{0,1}(x_0)$ is
the set to which we want to apply Banach's contraction mapping
principle. To ensure that the map $M_N(\cdot,p)$ maps into
$B_{6\delta}^{0,1}(x_0)$ for $p\in U$, and that $M_N(\cdot,p)$ is a
contraction we pick $N$ large enough. Specifically, we pick $N$ such
that
\begin{equation}\label{eq:nchoice}
  \begin{aligned}
    \|Q_Nx_0\|_{0,1}&\leq 2\delta\mbox{,}&
    \|Q_NL\|_{0,1}&\leq \frac{2\delta}{R}\mbox{,}\\
    \|Q_NL\|_0&\leq \frac{1}{2K}\mbox{,}& C\frac{\log
      N}{N}&< 1/\max\left\{1,\left(R+\|x_0\|_{1,1}\right)/\delta\right\}\mbox{,}
  \end{aligned}
\end{equation}
where $R$ and $K$ are the bounds on $F$ given in the conditions of the
lemma, in Equation~\eqref{eq:fbounds}. We know that these bounds exist
due to Lemma~\ref{thm:loclip} (see Equation~\eqref{eq:flip0}) and
Lemma~\ref{thm:Flipbound} (see Equation~\eqref{eq:Fboundlip}).  We
know that choosing $N$ according to \eqref{eq:nchoice} is possible
from Lemma~\ref{thm:qnl} and estimate \eqref{eq:qn} following
Lemma~\ref{thm:qnl}.

Let us check first that $x\mapsto E_Np+Q_NLF(x)$ maps the closed ball
$B_{6\delta}^{0,1}(x_0)$ back into itself:
\begin{align*}
  \lefteqn{\left\|E_Np+Q_NLF(x)-x_0\right\|_{0,1}\leq}&\\
  &\leq\left\|E_N\left[p-R_Nx_0\right]-Q_Nx_0+Q_NLF(x)\right\|_{0,1}\\
  &\leq \left\|E_N\left[p-R_Nx_0\right]\right\|_{0,1}+
  \left\|Q_Nx_0\right\|_{0,1}+
  \left\|Q_NL\right\|_{0,1}\left\|F(x)\right\|_{0,1}\\
  &<
  2\delta+2\delta+\frac{2\delta}{R}R=6\delta\mbox{.}
\end{align*}
Here we used the bounds \eqref{eq:nchoice} implied by our choice of
$N$ and the definition \eqref{eq:pclose} of the set $U$ of permitted
$p$, and the bound on $\|F(x)\|_{0,1}$, which is determined in
\eqref{eq:fbounds} by our choice of $\delta$.

Second, let us check that $x\mapsto E_Np+Q_NLF(x)$ is a uniform
contraction in $B_{6\delta}^{0,1}$ with respect to the
$\|\cdot\|_0$-norm:
\begin{align*}
  \left\|Q_NL\left[F(x)-F(y)\right]\right\|_0\leq
  \left\|Q_NL\right\|_0\left\|F(x)-F(y)\right\|_0\leq
  \frac{1}{2K}K\|x-y\|_0 \leq \frac{1}{2}\|x-y\|_0\mbox{.}
\end{align*}
Again, we exploited the bounds \eqref{eq:nchoice}, implied by our
choice of $N$, and the Lipschitz constant $K$ of $F$ determined in
\eqref{eq:fbounds} by our choice of $\delta$.

Since $B_{6\delta}^{0,1}(x_0)$ is complete with respect to the
$\|\cdot\|_0$-norm Banach's contraction mapping principle implies that
the fixed point problem \eqref{eq:fixp} has a unique solution $x\in
B_\delta^{0,1}(x_0)$ for $p\in U$.

Finally, let us check that for $x\in B_\delta^{0,1}(x_0)\cap
C^1(\T;\R^n)$ satisfying the periodic BVP $x'=F(x)$ the projection
$p=R_Nx$ is in $U$. For this we have to prove that if
$\|x-x_0\|_{0,1}\leq \delta$ and $x'=F(x)$ then
$\|P_N(x-x_0)\|_{0,1}< 2\delta$. We can estimate $\|P_N(x-x_0)\|_{0,1}$ via
\begin{align}
  \|P_N(x-x_0)\|_{0,1} &\leq \|(I-Q_N)(x-x_0)\|_{0,1}
  \leq \|x-x_0\|_{0,1}+\|Q_N(x-x_0)\|_{0,1}
  \label{eq:pnxx0:tri}\\
  &\leq \delta +C\frac{\log N}{N}\|x-x_0\|_{1,1}
  \label{eq:pnxx0:estqnl}\\
  &\leq\delta +C\frac{\log N}{N}\max\{|x-x_0\|_{0,1},\|x'-x_0'\|_{0,1}\}
  \label{eq:pnxx0:defn11}\\
  &\leq \delta +C\frac{\log N}{N}\max\{\delta,\|x'\|_{0,1}+\|x_0'\|_{0,1}\}
  \label{eq:pnxx0:balltri}\\
  &=  \delta +C\frac{\log N}{N}\max\{\delta,\|F(x)\|_{0,1}+\|x_0\|_{1,1}\}
  \label{eq:pnxx0:xpFx}\\
  &\leq \delta+C\frac{\log N}{N}\max\{\delta,R+\|x_0\|_{1,1}\}<2\delta\mbox{.}
  \label{eq:pnxx0:fbound}
\end{align}
The inequality \eqref{eq:pnxx0:tri} follows from the definition of
$P_N$ and $Q_N$ and the triangular inequality for the
$\|\cdot\|_{0,1}$-norm. The step from \eqref{eq:pnxx0:tri} to
\eqref{eq:pnxx0:estqnl} uses the estimate~\eqref{eq:qn} for the norm
$\|Q_Ny\|_{0,1}$ for elements $y$ of $C^{1,1}(\T;\R^n)$. It also
bounds $\|x-x_0\|_{0,1}$ by the radius $\delta$ of the ball. Step
\eqref{eq:pnxx0:defn11} splits up the $\|\cdot\|_{1,1}$ norm into its
two parts which are estimated separately in the following steps. One
part, $\|x-x_0\|_{0,1}$ is bounded by $\delta$ (the radius of the
ball), the difference of the derivatives is bounded by a triangular
inequality for its parts, $\|x'\|_{0,1}$ and $\|x_0'\|_{0,1}$ in
\eqref{eq:pnxx0:balltri}. To get to \eqref{eq:pnxx0:xpFx} we use that
$x$ satisfies the BVP $x'=F(x)$. We also bound the norm of $x_0'$ by
$\|x_0\|_{1,1}$. Finally, in \eqref{eq:pnxx0:fbound} we estimate the
Lipschitz norm of $F(x)$, $\|F(x)\|_{0,1}$ by the bound $R$ determined
in \eqref{eq:fbounds} by our choice of $\delta$. The right-hand side
of \eqref{eq:pnxx0:fbound} is (strictly) less than $2\delta$ by our
choice of $N$, see \eqref{eq:nchoice}.  \eop
\subsection{Lipschitz continuity of the algebraic system}
\label{sec:algred}
The Splitting Lemma~\ref{thm:split} guarantees in combination with the
unique existence of the fixed point of $M_N(\cdot,p)$, proven in
Lemma~\ref{thm:fixp}, the equivalence between the periodic BVP $\dot
x(t)=f(\Delta_tx)$ and the algebraic equation $g(p)=0$ for $x$ inside
the ball $B_\delta^{0,1}(x_0)$, where $g$ is given in
\eqref{eq:lowmodes:intro} by
\begin{align}
  \label{eq:lowmodes}
  g&:p\in U\mapsto 
  R_N\left[P_0F(X(p))+Q_0\left(E_Np-P_NLF(X(p))\right)\right]
  \in\R^{n\times(2\,N+1)}\mbox{, where}\\
  \label{eq:xpdef}
  X&:p\in U\mapsto C^0(\T;\R^n) \mbox{,\quad and $X(p)$ is the fixed point
    of $M_N(\cdot,p)$ in $B_{6\delta}^{0,1}(x_0)$}\mbox{.}
\end{align}
The relation between $p\in U$ and $x\in B_\delta^{0,1}(x_0)$ is given
via $p=R_Nx$ and $x=X(p)$: if $x$ satisfies the periodic BVP then
$p=R_Nx$ satisfies $g(p)=0$, and, vice versa, if $p\in U$ satisfies
$g(p)=0$ then $x=X(p)$ satisfies the periodic BVP. The domain of
definition, $D(X)=U$ is an open set, however the map $X$ (and, thus,
$g$) can be extended continuously to the boundary of $U$:
$M_N(\cdot,p)$ maps into the closed ball $B_{6\delta}^{0,1}$ back into
itself also for $p$ on the boundary of $U$ and it still has
contraction rate $1/2$ with respect to the $\|\cdot\|_0$-norm.




The remainder of the section addresses the remaining open claim of the
Equivalence Theorem~\ref{thm:main}, namely the regularity of the maps
$X$ and $g$.  Using only local $EC$ Lipschitz continuity
(Definition~\ref{def:loclip}) we can prove the Lipschitz continuity of
$g$ and $X$:
\begin{lemma}[Regularity of $X$ and algebraic
  system]\label{thm:lipalg}
    \begin{enumerate}
    \item\label{thm:Xsmooth1} For all $p$ in the neighborhood $U=D(X)$,
      defined in \eqref{eq:pclose}, the image $X(p)$ is in
      $C^{1,1}(\T;\R^n)$ (that is, $X(p)\in C^1(\T;\R^n)$ and its time
      derivative is Lipschitz continuous),
    \item\label{thm:Xlip1} $X$ is Lipschitz continuous with respect to
      the $\|\cdot\|_1$-norm for its images\textnormal{:} there exists
      a constant $C_N$ such that
      \begin{displaymath}
        \|X(p)-X(q)\|_1\leq C_N|p-q|\mbox{ for all $p$ and $q$ in $U$,}
      \end{displaymath}
    \item\label{thm:algnonlip} the map $p\in U\mapsto
      \left[R_0F(X(p)),P_NLF(X(p))\right]\in\R^n\times\R^{n\times(2\,N+1)}$
      is Lipschitz continuous in $U$.
    \end{enumerate}
\end{lemma}

\pstart{Proof} For a function $y\in\rg P_N$, differentiation is a
bounded operator: $y'=D_Ny$. The vector $R_Ny$ of the first $2N+1$
Fourier coefficients of a function $y$ and the vector $R_N[y']$
satisfy $R_N[y']=\tilde D_NR_Ny$ where $\tilde D_N$ is a matrix
(independent of $y$). Hence, $y'=E_N\tilde D_NR_Ny$ for all $y\in\rg
P_N$ such that we can define $D_N=E_N\tilde D_NR_N$.  Denote $X(p)$ as
$x$. By definition of the map $X$, $x=E_Np+Q_NLF(x)$. The right-hand
side of this fixed-point equation is differentiable with respect to
time, giving
\begin{equation}\label{eq:xpform}
  x'=D_NE_Np+Q_0F(x)-D_NP_NLF(x)\mbox{.}
\end{equation}
This guarantees that $x\in
C^1(\T;\R^n)$. Equation~\eqref{eq:Fboundlip} ensures that
$\|F(x)\|_{0,1}\leq R$, which implies that the right-hand side of
\eqref{eq:xpform} is Lipschitz continuous in time. This in turn
implies that $x'$ is Lipschitz continuous in time (thus, $x\in
C^{1,1}(\T;\R^n)$), and
\begin{displaymath}
  \|x'\|_{0,1}
  \leq\|D_NE_N\|_{0,1}|p|+\|Q_0\|_{0,1}R+\|D_NP_NL\|_{0,1}R\mbox{.}  
\end{displaymath}

Representation \eqref{eq:xpform} also implies point~\ref{thm:Xlip1}:
let $x=X(p)$ and $y=X(q)$ be two functions in the image of $X$:
\begin{equation}\label{eq:xlip1norm}
  \|x'-y'\|_0\leq\|D_NE_N\|_0|p-q|+(\|Q_0\|_0+\|D_NP_NL\|_0)K\|x-y\|_0\mbox{,}
\end{equation}
where $K$ was the $EC$ Lipschitz constant of $F$ in
$B_{6\delta}^{0,1}(x_0)$. The difference $x-y$ in the
$\|\cdot\|_0$-norm is bounded due to the contractivity of the
right-hand side in fixed point problem \eqref{eq:fixp} defining $X$
(the $\|\cdot\|_0$-norm was the metric used to apply the contraction
mapping principle):
\begin{align*}
  \|x-y\|_0&\leq \|E_N\|_0|p-q|+\|Q_NL[F(x)-F(y)\|_0\leq  
  \|E_N\|_0|p-q|+\frac{1}{2}\|x-y\|_0\mbox{.}
  \intertext{Thus,}
  \|x-y\|_0&\leq 2\|E_N\|_0|p-q|\mbox{,}
\end{align*}
which, combined with \eqref{eq:xlip1norm}, gives Lipschitz continuity of
$X$ as a map from $U$ into $C^1(\T;\R^n)$:
\begin{equation}\label{eq:xlip10norm}
  \|x'-y'\|_0\leq\left[\|D_NE_N\|_0+(\|Q_0\|_0+\|D_NP_NL\|_0)2K\|E_N\|_0\right]
  |p-q|=:C_N|p-q|\mbox{.}
\end{equation}

Point~\ref{thm:algnonlip} is a direct consequence of the Lipschitz
continuity of $F$ with respect to $\|\cdot\|_0$-norm in
$B_{6\delta}^{0,1}(x_0)$, the Lipschitz continuity of $X$ on $U$ in
the $\|\cdot\|_0$-norm, and the fact that $X$ maps into
$B_{6\delta}^{0,1}(x_0)$. \eop

\subsection{First-order differentiability 
  of the algebraic system}
\label{sec:algdiff1}
Until now we have only used the $EC$ Lipschitz continuity (in the
sense of Definition~\ref{def:loclip}) of the right-hand side $F$ in
the ball $B_{6\delta}^{0,1}(x_0)$ with respect to the
$\|\cdot\|_0$-norm. We can expect that the right-hand side $g$ of the
algebraic system, defined in \eqref{eq:lowmodes}, is smooth only if we
require more smoothness of the right-hand side $f$ (which enters $F$
in the algebraic system).

We first discuss first-order differentiability of the map $X$ and the
right-hand side $g$, defined in \eqref{eq:lowmodes} and
\eqref{eq:xpdef}. For this we assume $EC^1$ smoothness of $f$ as
defined in Definition~\ref{def:extdiff}. For $x\in C^1(\T;\R^n)\cap
B_{6\delta}^{0,1}(x_0)$ the norm of $\partial^1f(x,\cdot)$ as an
element of $L(C^0(\T;\R^n);\R^n)$ (the space of linear functionals
mapping $C^0(\T;\R^n)$ into $\R^n$) is less than or equal to $K$, the
$EC$ Lipschitz constant of $F$ (and, hence, $f$) in
$B_{6\delta}^{0,1}(x_0)$ assumed to exist in the conditions of
Lemma~\ref{thm:fixp}.

Let us define the map
\begin{align*}
  \partial^1F:&C^1(\T;\R^n)\times C^0(\T;\R^n)\mapsto
  C^0(\T;\R^n)\mbox{,} &
  \left[\partial^1F(v,w)\right](t)&= \partial^1f(\Delta_t
  v,\Delta_tw)\mbox{.}
\end{align*}
If $v\in C^1(\T;\R^n)$ and $w\in C^{0,1}(\T;\R^n)$ then the map
$\partial^1F$ defined above is indeed the derivative of $F$ in
$v$ with respect to the deviation $w$ (see Lemma~\ref{thm:Fdiff} in
Appendix~\ref{sec:basicprop}):
\begin{equation}\label{eq:Fdiff1}
  \lim_{
    \begin{subarray}{c}
      w\in C^{0,1}(\T;\R^n)\\[0.2ex]
      \|w\|_{0,1}\to0
    \end{subarray}
  }
  \frac{\|F(v+w)-F(v)-\partial^1F(v,w)\|_0}{\|w\|_{0,1}}=0\mbox{.}
\end{equation}
Part of the definition of $EC^1$ smoothness for $f$ is that the map
$\partial^1f$ is continuous in both arguments, $v\in C^1(\T;\R^n)$ and
$w\in C^0(\T;\R^n)$. One can then apply Lemma~\ref{thm:Fcont} to
$\partial^1f$ to conclude that the map $\partial^1F$ (a composition of
$\Delta_t$ and $\partial^1f$) is continuous with respect to the
$\|\cdot\|_0$-norm in its image space as a map of both arguments (in
their respective norm),
too. 
For $v\in B_{6\delta}^{0,1}(x_0)$ the norm of the linear map
$\partial^1F(v,\cdot)$ as an element of
$L(C^0(\T;\R^n);C^0(\T;\R^n))$, the space of continuous linear
functionals from $C^0(\T;\R^n)$ back to itself, is bounded by the $EC$
Lipschitz constant $K$ of $F$ in $B_{6\delta}^{0,1}(x_0)$.

The additional regularity assumption on $f$ and its implications for
$F$ permit us to improve our statements about regularity of $X$ and
the algebraic system:
\begin{lemma}[Continuous differentiability of $X$ and the algebraic
  system]\label{thm:Xdiff}
  Assume that the right-hand side $f$ is $EC^1$ smooth in the
  sense of Definition~\ref{def:extdiff}. Then the regularity
  statements about the map $X$, defined in \eqref{eq:xpdef}, and the
  right-hand side of the algebraic system, defined in
  \eqref{eq:lowmodes}, can be extended\textnormal{:}
  \begin{enumerate}
  \item $X(p)$ is in $C^2(\T;\R^n)$ for all $p\in U=D(X)$, the domain
    of definition of $X$, and $p\mapsto X(p)$ is continuous with
    respect to the $\|\cdot\|_2$-norm for its images.
  \item The map $X$, which maps $U$ into $C^1(\T;\R^n)$ according to
    Lemma~\ref{thm:lipalg}, is continuously differentiable with
    respect to its argument $p$ using the $\|\cdot\|_1$-norm for its images.
  \item The map $p\in U\mapsto
    \left[R_0F(X(p)),P_NLF(X(p))\right]\in\R^n\times\R^{n\times(2\,N+1)}$
    is continuously differentiable with respect to $p$.
  \end{enumerate}
\end{lemma}
\pstart{Proof} Let $p\in U=D(X)\subset \R^{n\times(2\,N+1)}$, where
$U$ is defined in \eqref{eq:pclose}, and let us denote $X(p)$ by
$x$. Lemma~\ref{thm:lipalg} ensures already that $x$ is in
$C^{1,1}(\T;\R^n)$. Lemma~\ref{thm:Fimage} in
Appendix~\ref{sec:basicprop} proves that $F(x)\in C^1(\T;\R^n)$ for
$x\in C^1(\T;\R^n)$ (choosing $D=C^0(\T;\R^n)$ and $k=0$ in the
assumptions of Lemma~\ref{thm:Fimage}). This implies the first
statement, that $X(p)\in C^2(\T;\R^n)$: since
\begin{equation}\label{eq:Xdiff:proof:xpc2}
  X(p)=E_Np+Q_NLF(X(p))
\end{equation}
and $X(p)\in C^{1,1}(\T;\R^n)$ (see Lemma~\ref{thm:lipalg}), $F(X(p))$
is in $C^1(\T;\R^n)$, and, thus, $LF(X(p))$ is in
$C^2(\T;\R^n)$. Hence, $X(p)$ is an element of $C^2(\T;\R^n)$,
too. Furthermore, Lemma~\ref{thm:Fimage} states that $F$ is continuous
as a map from $C^1(\T;\R^n)$ into $C^1(\T;\R^n)$. Since $X$ is
continuous as a map from $U$ into $C^1(\T;\R^n)$ (in fact, it is
Lipschitz continuous, see Lemma~\ref{thm:lipalg}), the right-hand side
of \eqref{eq:Xdiff:proof:xpc2} in $p$ is continuous with respect to the
$\|\cdot\|_1$-norm. This proves the first point.

Concerning the second statement: again, let $p_0$ be in $U=D(X)$, and
choose a small open neighborhood $U(p_0)$ which has a positive
distance to the boundary of $U$. We will prove point two for all $p\in
U(p_0)$. Let us choose an initial $\epsilon_0$ sufficiently small such
that $p+hq$ is still in $U$ for $h\in(-\epsilon_0,\epsilon_0)$, all
$p\in U(p_0)$, and all $q$ with $|q|\leq 1$. Let us introduce the
difference quotient for $h\in(-\epsilon_0,\epsilon_0)\setminus\{0\}$:
\begin{equation}\label{eq:Xdiffproof:zdef}
    z(h,q,p)=\frac{1}{h}\left[X(p+hq)-X(p)\right] \mbox{.}
\end{equation}
The maps $z$ maps
$\left[(-\epsilon_0,\epsilon_0)\setminus\{0\}\right]\times B_1(0)\times
U(p_0)\subset \R\times \R^{n\times(2\,N+1)}\times
\R^{n\times(2\,N+1)}$ into $C^1(\T;\R^n)$.  We first prove that $z$
has a limit for $h\to0$ in $C^1(\T;\R^n)$, and that this limit is
achieved uniformly for all $p\in U(p_0)$ and $|q|\leq 1$. By
definition of $X$, $z$ satisfies the fixed point equation (dropping
all arguments from $z$)
\begin{equation}\label{eq:gfixpd0}
  z=E_Nq+Q_NL\frac{1}{h}\left[F(X(p)+hz)-F(X(p))\right]
\end{equation}
for $h\in(-\epsilon_0,\epsilon_0)\setminus\{0\}$. Let us introduce
\begin{equation}\label{eq:Xdiffproof:tA1def}
  \tilde A_1(p,z,h)=
  \begin{cases}
    \frac{1}{h}\left[F(X(p)+hz)-F(X(p))\right] & \mbox{if $h\neq0$}\\
    \partial^1F(X(p),z) & \mbox{if $h=0$,}
  \end{cases}
\end{equation}
which maps $U(p_0)\times C^{0,1}(\T;\R^n)\times\R$ into
$C^0(\T;\R^n)$. The limit \eqref{eq:Fdiff1} implies that $\tilde A_1$
is continuous in all arguments (insert $v=x$, $w=hz$ into
\eqref{eq:Fdiff1}).  Using $\tilde A_1$ we extend the fixed point
problem \eqref{eq:gfixpd0} to $h=0$:
\begin{equation}\label{eq:gfixpd1}
  z=E_Nq+Q_NL\tilde A_1(X(p),z,h)\mbox{.}
\end{equation}
The following intermediate lemma proves that the fixed point problem
\eqref{eq:gfixpd1} has a unique solution:
\begin{lemma}[Fixed point problem for linearization]\label{thm:fixplin}
  There exists an $\epsilon>0$ and constants $C_0>0$ and $C_1>0$ such
  that the map
  \begin{displaymath}
    \gamma:z\in C^{0,1}(\T;\R^n)\mapsto 
    E_Nq+Q_NL\tilde A_1(X(p),z,h)\in C^{0,1}(\T;\R^n)\mbox{,}
  \end{displaymath}
  which depends on the additional parameters $p$, $q$ and $h$, has a
  unique fixed point $z_*$ in
  \begin{displaymath}
    B=\{z\in C^{0,1}(\T;\R^n):\|z\|_0\leq C_0 \mbox{ and }
    \|z\|_{0,1}\leq C_1\}
  \end{displaymath}
  for all $h\in(-\epsilon,\epsilon)$, all $p\in U(p_0)\subset
  U=D(X)\subset\R^{n\times(2\,N+1)}$ and all $q\in\R^{n(2N+1)}$ with
  $|q|<1$. The fixed point $z_*$ is an element of $C^1(\T;\R^n)$ and
  depends continuously on $h$, $p$ and $q$ with respect to the
  $\|\cdot\|_1$-norm.
\end{lemma}
Note that the $\epsilon$ we have to choose in Lemma~\ref{thm:fixplin}
is smaller than the initial $\epsilon_0$ for which the difference
quotient $z$ is defined. 

\pstart{Proof of Lemma~\ref{thm:fixplin}} First of all, since $\tilde
A_1$ is continuous in all arguments, the map $\gamma$ is
continuous. Moreover, since $x'=(X(p))'$ and $x=X(p)$ depend
continuously on $p$ (see Lemma~\ref{thm:lipalg} and expression
\eqref{eq:xpform}), the map $\gamma$ also depends continuously on the
parameters $p$, $q$ and $h$ (that is, the expression $E_nq+Q_NL\tilde
A_1(X(p),z,h)$, defining $\gamma$, depends continuously on $z$, $p$, $q$
and $h$ with respect to the $\|\cdot\|_{0,1}$-norm).  We choose the
constants $C_0>0$ and $C_1>0$ such that
\begin{align}
  \label{eq:c0choice}
  C_0&\geq2\|E_N\|_0\\
  \label{eq:c1choice}
  C_1&\geq \|D_NE_N\|_0+
  \left(\|Q_0\|_0+\|D_NP_NL\|_0\right)KC_0\mbox{,}
\end{align}
where $K$ is the Lipschitz constant of $F$ with respect to the
$\|\cdot\|_0$-norm in $B_{6\delta}^{0,1}$.

We choose $\epsilon\leq\epsilon_0$ such that for all $z$ satisfying
$\|z\|_{0,1}\leq C_1$ and all $p\in U(p_0)$ the function $X(p)+hz$ is
in $B_{6\delta}^{0,1}(x_0)$ for all $h\in(-\epsilon,\epsilon)$. This
implies that for any $z_1$ and $z_2$ satisfying $\|z_1\|_{0,1}\leq
C_1$ and $\|z_2\|_{0,1}\leq C_1$ we have
\begin{equation}
\begin{split}
  \frac{1}{h}\|F(X(p)+hz_1)-F(X(p)+hz_2)\|_0&\leq K\|z_1-z_2\|_0\mbox{,}\\
  \|\partial^1F(X(p),z_1-z_2)\|_0&\leq K\|z_1-z_2\|_0\mbox{,}
\end{split}\label{eq:g0lip}
\end{equation}
where $K$ was the Lipschitz constant for $F$ in
$B_{6\delta}^{0,1}(x_0)$, and, thus,
\begin{align}
  \|\gamma(z_1)-\gamma(z_2)\|_0&\leq \frac{1}{2}\|z_1-z_2\|_0\label{eq:g0contract}
\end{align}
for all $h\in(-\epsilon,\epsilon)$ by choice of $N$ ($N$ was such that
$\|Q_NL\|_0\leq (2K)^{-1}$). This estimate for $\gamma$ implies
\begin{equation}\label{eq:g0bound}
  \|\gamma(z)\|_0\leq \|E_N\|_0+\frac{1}{2}\|z\|_0
  \mbox{\qquad if $\|z\|_{0,1}\leq C_1$,}
\end{equation}
since $\gamma(0)=E_Nq$ and $|q|\leq1$. Moreover, the two
inequalities~\eqref{eq:g0lip} imply that for
$h\in(-\epsilon,\epsilon)$, $\|z\|_{0,1}\leq C_1$ and $p\in U(p_0)$
the maximum norm of $\tilde A_1(p,z,h)$ is bounded by $K\|z\|_0$:
\begin{equation}
  \label{eq:tA1bound}
  \|\tilde A_1(p,z,h)\|_0\leq K\|z\|_0
\end{equation}
The time derivative of $\gamma(z)$ exists and its $\|\cdot\|_0$-norm can be
estimated by differentiating the expression $E_nq+Q_NL\tilde
A_1(X(p),z,h)$, defining $\gamma$, with respect to time in the same manner
as we obtained \eqref{eq:xlip1norm} (we insert \eqref{eq:tA1bound} to
bound $\|\tilde A_1(p,z,h)\|_0$):
\begin{equation}\label{eq:g1bound}
  \left\|\frac{\d}{\d t}\gamma(z)\right\|_0\leq
  \|D_NE_N\|_0+(\|Q_0\|_0+\|D_NP_NL\|_0)K\|z\|_0\mbox{.}
\end{equation}
The combination of the bounds \eqref{eq:g0bound} and
\eqref{eq:g1bound} and the definition of the constants $C_0$ and $C_1$
guarantee that $\gamma(z)$ maps the set
\begin{displaymath}
  B=\{z\in C^{0,1}(\T;\R^n):\|z\|_0\leq C_0 \mbox{ and }
  \|z\|_{0,1}\leq C_1\}
\end{displaymath}
back into itself. The contraction estimate \eqref{eq:g0contract} for
the $\|\cdot\|_0$-norm and the completeness of $B$ with respect to the
$\|\cdot\|_0$-norm make the contraction mapping principle applicable
with a uniform contraction rate for all $p\in U(p_0)$, all $|q|\leq1$
and $h\in(-\epsilon,\epsilon)$. This ensures that the fixed point
$z_*$ depends continuously on $p$, $q\in\R^{n(2N+1)}$ and
$h\in(-\epsilon,\epsilon)$ with respect to the $\|\cdot\|_0$-norm
(since $\gamma$ is continuous with respect to $z$, $h$, $q$ and $p$).

The time derivative $z_*'$ of $z_*$ also exists and is continuous in
$p$, $q$ and $h$: we differentiate the fixed point equation
\eqref{eq:gfixpd1} for $z_*$ with respect to time (in the same way as
done in \eqref{eq:xpform}) to get
\begin{align}\label{eq:zpform}
  z_*'=&D_NE_Nq+Q_0\tilde A_1(X(p),z_*,h)- D_NP_NL\tilde
  A_1(X(p),z_*,h)\mbox{,}
\end{align}
which is a continuous function in $p$, $q$ and $h$ with respect to the
$\|\cdot\|_0$-norm (note that $z_*$ depends on $p\in U(p_0)$, $q$ and
$h$). Thus, the fixed point $z_*$ is in $C^1(\T;\R^n)$ and depends
continuously on $p$, $q$ and $h$ with respect to the
$\|\cdot\|_1$-norm.\eop

\pstart{Proof of Lemma~\ref{thm:Xdiff} continued} As a consequence of
Lemma~\ref{thm:fixplin} we may write the fixed point $z_*$ of $\gamma$ as a
function of $h$, $q$ and $p$: $z_*(h,q,p)$ maps
$h\in(-\epsilon,\epsilon)$, $q$ in the unit ball of
$\R^{n\times(2\,N+1)}$ and $p\in U(p_0)$ continuously into
$C^1(\T;\R^n)$. It is also identical to $z(h,q,p)$, defined in
\eqref{eq:Xdiffproof:zdef} as the directional difference quotient of
$X$. Thus, the directional difference quotient $z(h,q,p)$ has a limit
for $h\to0$ in the $\|\cdot\|_1$-norm, and this limit equals
$z_*(0,q,p)$. Moreover, this limit $z_*(0,q,p)$ depends continuously
on $p$ and $q$ in the $\|\cdot\|_1$-norm (as proved in
Lemma~\ref{thm:fixplin}), and it is linear in $q$ (since $\tilde
A_1(p,z,0)$ is linear in $z$).  Thus, $z_*(0,q,p)$ is the Frech{\'e}t
derivative:
\begin{equation}\label{eq:Xdiff1proof:Frechet0}
  \lim_{q\to0}\frac{\|X(p+q)-X(p)-z_*(0,q,p)\|_1}{|q|}=0\mbox{.}
\end{equation}
 Consequently, the map
$(p,q)\mapsto z_*(0,q,p)=\partial^1X(p)\,q$ is continuous in the
$\|\cdot\|_1$-norm as claimed in the lemma.

The third statement of Lemma~\ref{thm:Xdiff} is a consequence of the
second statement and the fact that the difference quotient of $F$ has
a limit in the $\|\cdot\|_0$-norm if it is taken between arguments in
$C^1(\T;\R^n)$ (see \eqref{eq:Fdiff1}). We split the difference
quotients into two parts:
\begin{align}\label{eq:Fd1expr1}
  \frac{F(X(p+hq))-F(X(p))}{h}=&\ \frac{F(X(p)+h\partial^1
    X(p)q)-F(X(p))}{h} +\\
  &\ +\frac{F(X(p+hq))-F(X(p)+h\partial^1 X(p)q)}{h}\label{eq:Fd1expr2}
\end{align}
The right-hand side in \eqref{eq:Fd1expr1} converges in the
$\|\cdot\|_0$-norm to $\partial^1F(X(p),\partial^1X(p)q)$ for $h\to0$,
since $X(p)$ and $\partial^1X(p)q$ are in $C^1(\T;\R^n)$ because $F$
is $EC^1$ continuous (see the second point of the lemma for the
regularity of $\partial^1X(p)q$ and Lemma~\ref{thm:lipalg} for the
regularity of $X(p)$). For the term in \eqref{eq:Fd1expr2} we can
apply the local $EC$ Lipschitz continuity (all arguments are in
$B_{6\delta}^{0,1}(x_0)$ for $p\in U(p_0)$, $|q|\leq 1$ and
$h\in(\epsilon,\epsilon)$) such that we get
\begin{displaymath}
  \left\|\frac{F(X(p+hq))-F(X(p)+h\partial X(p)q)}{h}\right\|_0\leq
  K\left\|\frac{X(p+hq)-X(p)}{h}- 
    \partial X(p)q\right\|_0\mbox{,}
\end{displaymath}
which converges to $0$ for $h\to0$ due to the second statement of the
lemma ($K$ is the $EC$ Lipschitz constant of $F$ in
$B_{6\delta}^{0,1}(x_0)$).  Consequently, we obtain from the limit of
\eqref{eq:Fd1expr1} for $h\to0$ that the directional derivative of
$F(X(p))$ in direction $q$ is equal to $\partial^1F(X(p),\partial^1
X(p)q)$, which is continuous with respect to $p$ and $q$ and linear in
$q$. Thus,
\begin{equation}\label{eq:Fd1expr}
  \left[\frac{\partial}{\partial p}F(X(p))\right]q
  =\partial^1F(X(p),\partial X(p)q)\mbox{,}  
\end{equation}
and $p\mapsto F(X(p))$ is continuously differentiable with respect to
$p$ in the $\|\cdot\|_0$-norm. Note that we use the notation not
enclosing $q$ in the bracket in \eqref{eq:Fd1expr} to highlight that
this is a classical derivative with respect to a finite-dimensional
variable. The linear operators $R_0$ and $P_NL$ preserve the
continuity (and the linearity in $q$) of \eqref{eq:Fd1expr}. \eop

\begin{lemma}[Spectral radius of the linearized fixed point problem]
  \label{thm:specradius} For $x=X(p)$ (where $p\in U=D(X)$) consider
  the linear map
  \begin{displaymath}
    M:z\mapsto Q_NL\partial^1F(x,z)\mbox{.}
  \end{displaymath}
  The spectral radius of $M$ as a map from $C^0(\T;\R^n)$ back into
  itself, or as a map from $C^1(\T;\R^n)$ back into itself, is less or equal $1/2$.
\end{lemma}

\pstart{Proof} Since $M$ is compact as an element of
$L(C^k(\T;\R^n);C^k(\T;\R^n))$ (the space of linear functionals from
$C^k(\T;\R^n)$ back to itself) for $k=0$ and $k=1$, the spectral radius
is identical to the modulus of the maximal (in modulus) eigenvalue,
which is of finite algebraic multiplicity if it is different from
zero. An eigenvector $z$ corresponding to this maximal eigenvalue is
an element of $C^1(\T;\R^n)$ such that the spectral radius of $M$ is
the same for $k=0$ and $k=1$.

Since $x$ and $z$ are both in  $C^1(\T;\R^n)$ we have that
\begin{equation}\label{eq:sradproof:d1f}
\partial^1F(x)z=\lim_{h\to0}\frac{1}{h}\left[F(x+hz)-F(x)\right]
\end{equation}
For $x=X(p)$ where $p\in U$, and $h$ sufficiently small the arguments
of $F$, $x+hz$ and $x$, both lie inside $B_{6\delta}^{0,1}$ such that
the $EC$ Lipschitz constant $K$ applies to the difference:
\begin{displaymath}\label{eq:sradproof:flip}
  \frac{1}{h}\|F(x+hz)-F(x)\|_0\leq K\|z\|_0\mbox{.}
\end{displaymath}
Since $\|Q_NL\|_0\leq 1/(2K)$,
\eqref{eq:sradproof:d1f} and \eqref{eq:sradproof:flip} combine to
\begin{displaymath}
  \|Mz\|_0\leq \frac{1}{2}\|z\|_0\mbox{.}
\end{displaymath}
As $z$ is an eigenvector corresponding to the largest eigenvalue, the
spectral radius of $M$ is less or equal $1/2$.\eop

Thus, the derivative $z=\tpartial X(p)\,q$ of $X$ in $p$ is the unique
solution of the contractive linear fixed point problem in
$C^1(\T;\R^n)$
\begin{equation}\label{eq:dxfixp}
  z=E_Nq+Q_NL\partial^1F(X(p),z)\mbox{.}
\end{equation}

\subsection{Higher degrees of smoothness}
\label{sec:smooth}
We observe that $(x,y)=(X(p),\partial^1X(p)\,q)$ satisfies the
system of equations 
\begin{equation}\label{eq:fixp:ext}
  \begin{split}
    x&=E_Np+Q_NLF(x)\\
    y&=E_Nq+Q_NL\partial^1F(x,y)\mbox{.}
  \end{split}
\end{equation}
This has a similar structure to the original fixed point problem
\eqref{eq:fixp} but in dimension $n_1=2n$ with the variables $(x,y)$
and parameters $(p,q)$. Thus, we aim to apply a linear version of the
arguments of Section~\ref{sec:algdiff1} recursively, assuming that $f$
is $EC^k$ smooth as recursively defined in
Definition~\ref{def:extdiff}. Throughout this section we assume that
$f$ is $EC^k$ smooth for all degrees up to order $j_{\max}$.

For higher-order derivatives, we introduce the spaces $D_j$ and the
operators $\partial^jF$ for $j\geq0$ recursively:
\begin{align*}
  D_0&=C^0(\T;\R^n) & D_j&=D_{j-1}^1\times D_{j-1}\\
  \partial^jF:&D_j\mapsto C^0(\T;\R^n)\mbox{,} &
  [\partial^jF(x)](t)&=\partial^jf(\Delta_tx)\mbox{.}
\end{align*}
The spaces $D_j$ are products of the type \eqref{eq:Dspacedef}, and
the argument $x$ of $\partial^jF$ and $\partial^jf$ is in $D_j$, a
product of $2^j$ spaces. We also recall that the notion of subspaces
$D_j^k$ of higher-oder ($k\geq0$) differentiability for product spaces
such as $D_j$ was introduced in Section~\ref{sec:theorem}.
For example,
\begin{align*}
  D_0^k&=C^k(\T;\R^n)\mbox{,}\\
  D_1^k&=D_0^{k+1}\times D_0^k=C^{k+1}(\T;\R^n)\times C^k(\T;\R^n)\mbox{,}\\
  D_2^k&=D_1^{k+1}\times D_1^k=C^{k+2}(\T;\R^n)\times
  C^{k+1}(\T;\R^n)\times C^{k+1}(\T;\R^n)\times C^k(\T;\R^n)\mbox{, etc.,}
\end{align*}
all with their natural maximum norms. The maps $\partial^jF$ are all
continuous and map indeed into $C^0(\T;\R^n)$ due to the continuity of
$\partial^jf$ and $\Delta_t$ (applying Lemma~\ref{thm:Fcont} to $D_j$,
$\partial^jF$ and $\partial^jf$). It is also clear from the definition
that $\partial^{j+k}F=\partial^j[\partial^kF]$ if $j+k\leq
j_{\max}$. We will also use the notation $L(D_j^k;D_i^\ell)$ for the
space of linear bounded functionals mapping from $D_k^k$ into
$D_i^\ell$.

The following lemma is a consequence of the $EC^k$ smoothness of
$f$. 
\begin{lemma}\label{thm:rhs:smooth}
  For $j+k\leq j_{\max}$ the operator $\partial^jF$ is a continuous map from
  $D_j^k$ into $C^k(\T;\R^n)$.
\end{lemma}
\pstart{Proof of Lemma} We have to apply Lemma~\ref{thm:Fimage} from
Appendix~\ref{sec:basicprop} inductively over the order of
differentiability ($k$). To start the induction for $k=0$ we can apply
Lemma~\ref{thm:Fcont} to $D_j$, $\partial^jF$ and $\partial^jf$.  For
the inductive step let us assume that for $k$ we know that
$\partial^jF:D_j^k\mapsto C^k(\T;\R^n)$ is continuous for all $j\leq
j_{\max}-k$. Let us fix a $j\leq j_{\max}-k-1$. We have to show that
$\partial^jF$ maps $D_j^{k+1}$ continuously into
$C^{k+1}(\T;\R^n)$. We know (by inductive assumption) that
$\partial^jF$ maps $D_j^k$ continuously into $C^k(\T;\R^n)$ and that
$\partial^{j+1}F$ maps $D_{j+1}^k=D_j^{k+1}\times D_j^k$ continuously
into $C^k(\T;\R^n)$. Thus, we can apply Lemma~\ref{thm:Fimage} to
$\partial^jF$ (this takes the place of the operator $F$ in
Lemma~\ref{thm:Fimage}) and $D=D_j^k$, obtaining that
$\partial^jF:D_j^{k+1}\mapsto C^{k+1}(\T;\R^n)$ is continuous. \eop

An immediate consequence of Lemma~\ref{thm:rhs:smooth} is that $X(p)$
and $\partial X(p)\,q$, as constructed in Section~\ref{sec:algdiff1},
are as smooth as the right-hand-side:
\begin{lemma}[Smoothness of $X$ and $\partial X$ in
  time]\label{thm:xtsmooth}
  Let $f$ be $EC^{j_{\max}}$ smooth. For every $p\in U=D(X)$ and every
  $q\in R^{n(2N+1)}$ the functions $X(p)$ and $\partial X(p)\,q$
  satisfy $X(p)\in C^{j_{\max}+1}(\T;\R^n)$ and $\partial X(p)\,q\in
  C^{j_{\max}}(\T;\R^n)$. Moreover, the maps
  \begin{align*}
    p&\mapsto X(p)\in C^{j_{\max}+1}(\T;\R^n)\mbox{\quad and\quad} 
    [p,q]\mapsto \partial X(p)\,q\in C^{j_{\max}}(\T;\R^n)
  \end{align*}
 are continuous.
\end{lemma}
\pstart{Proof} The function $x=X(p)$ satisfies
$x=E_Np+Q_NLF(x)$. Since $F$ maps $D_0^k=C^k(\T;\R^n)$ back into
itself continuously for all $k\leq j_{\max}$, $Q_NL$ maps $D_0^k$ into
$D_0^{k+1}$ continuously for all $k$, and $E_Np\in C^\infty(\T;\R^n)$,
the fixed point equation implies the following: if $x\in D_0^k$ then
$F(x)\in D_0^k$, thus, $x=E_Np+Q_NLF(x)\in D_0^{k+1}$(for all $k\leq
j_{\max}$). Similarly, $z=E_Nq+Q_NL\partial^1F(x)\,z$, and
$\partial^1F$ maps $D_1^k$ into $D_0^k$ for all $k\leq
j_{\max}-1$. Thus, the fixed point equation implies: if $z\in D_0^k$
and $x\in D_0^{k+1}$ then $(x,z)\in D_1^k$, thus, $\partial^1F(x,z)\in
D_0^k$, thus, $z=E_Nq+Q_NL\partial^1F(x,z)\in D_0^{k+1}$ for all
$k\leq j_{\max}-1$. All of the above dependencies are continuous such
that the continuous dependence on $p$ and $q$ in the norms of
$D_0^{j_{\max}+1}$ and $D_0^{j_{\max}}$, respectively, follows.\eop

We plan to find the derivatives of the map $X$ inductively through
fixed point equations of the form \eqref{eq:fixp:ext}. In order to set
up the recursion we define inductively the operators $F_j$ by
\begin{align}
  F_0(x)&=F(x) &&\mbox{for $x\in D_0$}\label{eq:F0def}\\
  F_j
  \begin{pmatrix}
    x\\ y
  \end{pmatrix}
&=
\begin{bmatrix}
  F_{j-1}(x)\\
  \partial^1F_{j-1}(x,y)
\end{bmatrix}\mbox{,} &&\mbox{for\ }
\begin{bmatrix}
  x\\ y
\end{bmatrix}\in D_j=D_{j-1}^1\times D_j\mbox{.}\label{eq:Fjdef}
\end{align}
Note that $F_j$ is always linear in its second argument, $y$, since
$\partial^1F_{j-1}$ is linear in its second argument. The operators
$F_j$ are combinations of derivatives of $F$. The plan is to study
fixed-point problems of the type $x=E_Np+Q_NLF_j(x)$ (with $j=1$ we
obtain \eqref{eq:fixp:ext}). Before doing so, we establish which
spaces the operators $F_j$ map into:
\begin{lemma}[Image of right-hand side]\label{thm:rhs:image}
  For $j+l+k\leq j_{\max}$ the operator $\partial^lF_j$ maps
  $D_{j+l}^k$ continuously into $D_j^k$. In particular, $F_j$ maps
  $D_j$ continuously back into itself.
\end{lemma}
\pstart{Proof} The statement of the lemma follows inductively from the
definition of $F_j$ and $D_j^k$. We apply Lemma~\ref{thm:rhs:smooth}
to start our induction over $j$ (for $j=0$ the statement is identical
to Lemma~\ref{thm:rhs:smooth}). For the inductive step let us assume
that we know that $\partial^lF_{j-1}$ maps $D_{j+l-1}^k$ continuously
into $D_{j-1}^k$ for all $k$ and $l$ satisfying $l+k\leq
j_{\max}-j+1$. By definition \eqref{eq:Fjdef} of $F_j$ the derivative
$\partial^lF_j$ for $l\leq j_{\max}-j$ is
\begin{align*}
  \partial^lF_j
  \begin{pmatrix}
    x\\ y
  \end{pmatrix}=
  \begin{bmatrix}
    \partial^lF_{j-1}(x)\\
    \partial^{l+1}F_{j-1}(x,y)
  \end{bmatrix}
  &&\mbox{for\ }
  \begin{bmatrix}
    x\\ y
  \end{bmatrix}\in D_{l+j}=D_{l+j-1}^1\times D_{l+j-1}\mbox{.}
\end{align*}
The first component, $\partial^lF_{j-1}$ maps $D_{l+j-1}^{k+1}$
continuously into $D_{j-1}^{k+1}$ for all $k$ from $0$ to
$j_{\max}-l-j$ (this is the assumption of the inductive step when one
shifts the index $k$ down by $1$). Similarly, $\partial^{l+1}F_{j-1}$
maps $D_{j+l-1}^{1+k}\times D_{j+l-1}^k=D_{j+l}^k$ continuously into
$D_{j-1}^k$ for all $k$ from $0$ to $j_{\max}-l-j$, again due to the
assumption of the inductive step. Consequently, $\partial^lF_j$ maps
$D_{j+l-1}^{k+1}\times D_{j+l-1}^k=D_{j+l}^k$ continuously into
$D_{j}^k$ for all $k$ from $0$ to $j_{\max}-l-j$, which is the
statement we had to prove for the inductive step.\eop

Even though the map $x\in D_j^1\mapsto \partial^1F_j(x,\cdot)\in
L(D_j;D_j)$ is in general not continuous, the following map is:
\begin{lemma}[Continuity in operator norm]\label{thm:speccont}
  For $j<j_{\max}$ the map $x\in D_j^1\mapsto
  Q_NL\partial^1F_j(x,\cdot)\in L(D_j^1;D_j^1)$ is continuous with
  respect to $x\in D_j^1$.
\end{lemma}
\pstart{Proof of Lemma~\ref{thm:speccont}} The $EC^k$ smoothness of
$f$ (for $k\leq j_{\max}$) implies that $F_j$ is continuously
differentiable (in the classical sense) as a map from $D_j^1$ into
$D_j$. Thus, the map $x\mapsto\partial^1F_j(x,\cdot)$ as a map from
$D_j^1$ into $L(D_j^1;D_j)$ is continuous. Recall that the operator
$L$ involves taking the anti-derivative of its argument such that
$L:D_j\mapsto D_j^1$.  Since $Q_NL$ maps $D_j$ continuously into
$D_j^1$, the map $x\mapsto Q_NL\partial^1F_j(x,\cdot)$ is continuous
as a map from $D_j^1$ into $L(D_j^1;D_j^1)$.  \eop

The following theorem provides continuous differentiability of order
$j_{\max}$ for $X$ and the map $p\mapsto F(X(p))$ if the right-hand
side is $EC^k$ smooth in the sense of Definition~\ref{def:extdiff} for
$k\leq j_{\max}$:
\begin{theorem}[Smoothness of algebraic system and $X$]\label{thm:smooth}
  Define $n_0=n(2N+1)$ and $n_j=2^jn_0$, and the maps
  \begin{align*}
    X_0&: p\in U=D(X)\subseteq \R^{n_0}\mapsto X(p)\in D_0\mbox{\ and}\\
    Y_0&:p\in U=D(X)\subseteq \R^{n_0}\mapsto F(X(p))\in D_0\mbox{,}
  \end{align*}
  and assume that $f:D_0=C^0(\T;\R^n)\mapsto\R^n$ is $EC^{j_{\max}}$
  smooth. Then the following maps exist and are continuous for all $j$
  up to $j_{\max}$:
  \begin{align*}
    X_j&:[p,q]\in D(X_j):=D(X_{j-1})\times\R^{n_{j-1}}\subseteq\R^{n_j}
    \mapsto [X_{j-1}(p),\partial X_{j-1}(p)\,q]\in D_j\mbox{,}\\
    Y_j&:[p,q]\in D(X_j)
    \phantom{\ :=D(X_{j-1})\times\R^{n_{j-1}}\subseteq\R^{n_j}}
    \mapsto [Y_{j-1}(p),\partial Y_{j-1}(p)\,q]\in D_j\mbox{.}
  \end{align*}
\end{theorem}
The proof of Theorem~\ref{thm:smooth} does not require the application
of the contraction mapping principle for nonlinear maps. It uses only
Lemma~\ref{thm:specradius}, Lemma~\ref{thm:rhs:image} and
Lemma~\ref{thm:speccont} inductively.

\pstart{Proof of Theorem~\ref{thm:smooth}} The main work is the proof
of the existence and continuity of $X_j$, which we will do
first.
The assumption of the inductive step is
comprised of the following two statements. We assume for $j$:
\begin{enumerate}
\item\label{p:ass:exfix0} The map $(p_1,p_2)\in D(X_{j-1})\times
  \R^{n_{j-1}}\mapsto X_j(p_1,p_2)\in D_j$ exists and is
  continuous. Moreover, the pair $(x_1,x_2)=X_j(p_1,p_2)$ satisfies
  \begin{align}
    x_1&=E_Np_1+Q_NLF_{j-1}(x_1)\label{eq:ivf}\\
    x_2&=E_Np_2+Q_NL\partial^1F_{j-1}(x_1,x_2)\mbox{.}\label{eq:ivderiv}
  \end{align}
\item\label{p:ass:specradius} The linear map $z\mapsto
  Q_NL\partial^1F_{j-1}(x_1,z)$ maps $D_{j-1}^1$ back into itself and
  has spectral radius less or equal $1/2$.
\end{enumerate}
Both statements of the assumption of the inductive step have been
proven for $j=1$ in Lemma~\ref{thm:Xdiff} and
Lemma~\ref{thm:specradius} . Let $j$ be smaller than $j_{\max}$.

\subsubsection*{Regularity of  $X_j(p)$}
Let us first establish that the map
\begin{displaymath}
p=
\begin{bmatrix}
 p_1\\ p_2 
\end{bmatrix}\mapsto x=
\begin{bmatrix}
 x_1\\ x_2 
\end{bmatrix}
=X_j
  \begin{pmatrix}
    p_1\\ p_2    
  \end{pmatrix}
  =
  \begin{bmatrix}
    X_{j-1}(p_1)\\
    \partial^1X_{j-1}(p_1)\,p_2      
  \end{bmatrix}
\end{displaymath}
does not only map continuously into $D_j$ but into $D_j^k$ for all
$k\leq j_{\max}-j+1$. 

The argument is the same as in the proof of
Lemma~\ref{thm:xtsmooth}: the map $F_j$ maps $D_j^k$ continuously
back into $D_j^k$ for all $k\leq j_{\max}-j$. If $x\in D_j^k$ then
$F_j(x)\in D_j^k$, thus, $x=E_Np+Q_NLF(x)\in D_j^{k+1}$ for all $k\leq
j_{\max}-j$ (and the dependence on $p$ is continuous because all
dependencies are continuous).

\subsubsection*{Proof of existence and continuity of $\partial^1X_j(p)\,q$}
Let us use the notation $p=(p_1,p_2)$ and $x=(x_1,x_2)=X_j(p)$.  Let
$p_0\in D(X_j)$ be arbitrary. We first show that $\partial^1X_j(p)\,q$
exists for all $p$ in a neighborhood $U(p_0)$ with positive distance
to the boundary of $D(X_j)$. We can choose $\epsilon>0$ sufficiently
small such that $p+hq\in D(X_j)$ for all $h\in(-\epsilon,\epsilon)$,
all $q=(q_1,q_2)$ with $|q|<1$ and all $p\in U(p_0)$. Consider the
difference quotient
\begin{displaymath}
\frac{X_j(p+hq)-X_j(p)}{h}=
\frac{1}{h}
\begin{bmatrix}
X_{j-1}(p_1+hq_1)-X_{j-1}(p_1)\\
\partial^1X_{j-1}(p_1+hq_1)\,[p_2+hq_2]-\partial^1X_{j-1}(p_1)\,p_2
\end{bmatrix}=:
\begin{bmatrix}
  z_1\\ z_2
\end{bmatrix}
\mbox{.}
\end{displaymath}
By assumption of the inductive step, $X_{j-1}$ is continuously
differentiable such that the first row of this difference quotient
has the form
\begin{equation}\label{eq:xjm1p1q1}
  z_1(h,p_1,q_1)=\frac{1}{h}\left[X_{j-1}(p_1+hq_1)-X_{j-1}(p_1)\right]=
  \int_0^1\partial^1X_{j-1}(p_1+hsq_1)\,q_1\d s
\end{equation}
for $h\neq0$.  As established above
$(p_1,q_1)\mapsto \partial^1X_{j-1}(p_1)\,q_1 \in D_{j-1}^k$ is
continuous for all $k\leq j_{\max}-j+1$ such that
\begin{displaymath}
z_1(h,p_1,q_1)\in D_{j-1}^{j_{\max}-j+1} \subseteq D^2_{j-1} 
\end{displaymath}
($j_{\max}-j+1\geq2$ since $j<j_{\max}$), and $z(h,p_1,q_1)$ depends
continuously on its arguments, also when $h=0$. Let us use the
abbreviations
\begin{align*}
  x_1(p_1)&=X_{j-1}(p_1)\mbox{,} \\
  x_2(p_1,p_2)&=\partial^1X_{j-1}(p_1)p_2\mbox{,}\\
  z_1(h,p_1,q_1)&=\frac{1}{h}\left[X_{j-1}(p_1+hq_1)-X_{j-1}(p_1)\right]=
  \int_0^1\partial^1X_{j-1}(p_1+hsq_1)\,q_1\d s\\
  z_2(h,p_1,p_2,q_1,q_2)&=\frac{1}{h}\left[\partial^1X_{j-1}(p_1+hq_1)\,[p_2+hq_2]
    -\partial^1X_{j-1}(p_1)\,p_2\right]\mbox{\quad for $h\neq0$.}
\end{align*}
With these notations we have $X_{j-1}(p_1+hq_1)=x_1+hz_1$ and, for
non-zero $h$, $\partial^1X_{j-1}(p_1+hq_1)[p_2+hq_2]=x_2+hz_2$. The
fixed-point equations \eqref{eq:ivf} and \eqref{eq:ivderiv} imply a
fixed-point equation for the difference quotient $z_2$ for non-zero
$h$:
\begin{align}
  z_2&=E_Nq_2+\frac{1}{h}Q_NL\left[\partial^1F_{j-1}(x_1+hz_1,x_2+hz_2)-
    \partial^1F_{j-1}(x_1,x_2)\right]\nonumber\\
  &=E_Nq_2+\tilde
  z(p_1,p_2,q_1,h)+Q_NL\partial^1F_{j-1}(x_1+hz_1,z_2)
  \mbox{\quad where}\label{eq:z2fixp}\\
  \tilde z(p_1,p_2,q_1,h)&=Q_NL\frac{\partial^1F_{j-1}(x_1+hz_1,x_2)-
        \partial^1F_{j-1}(x_1,x_2)}{h}\nonumber
\end{align}
The regularity of $x_1$, $x_2$ and $z_1$ is: 
\begin{equation}\label{eq:x1x2z1reg}
  \begin{split}
    x_1&\in D_{j-1}^{j_{\max}-j+2}\subseteq D_{j-1}^3\mbox{,} \\
    x_2&\in D_{j-1}^{j_{\max}-j+1}\subseteq D_{j-1}^2\mbox{\quad and} \\
    z_1&\in D_{j-1}^{j_{\max}-j+1}\subseteq D_{j-1}^2\mbox{.}
  \end{split}
\end{equation}
We can apply the mean value theorem to the difference quotient
appearing in $\tilde z$ since $x_1$ and $z_1$ are at least in $D_{j-1}^2$
and $x_2$ is at least in $D_{j-1}^1$ (see
Lemma~\ref{thm:rhs:smooth}, and Lemma~\ref{thm:fmeanval} and
Lemma~\ref{thm:Fdiff} in Appendix~\ref{sec:basicprop}):
\begin{align*}
  \tilde z(p_1,p_2,q_1,h)&=
  Q_NL\int_0^1\partial^2F_{j-1}(x_1+shz_1,x_2,z_1,0)\d s\mbox{.}
\end{align*}
The map
$(x_1,x_2,z_1,h)\mapsto \partial^2F_{j-1}(x_1+shz_1,x_2,z_1,0)$ maps
$x_1$, $x_2$, $z_1$ and $h$ continuously into the space
$D_{j-1}^{j_{\max}-j-1}$ (we see this by applying
Lemma~\ref{thm:rhs:smooth} to $\partial^2F_{j-1}$, setting $k$ in
Lemma~\ref{thm:rhs:smooth} to $j_{\max}-j-1$). Thus, the quantity
$\tilde z(p_1,p_2,q_1,h)$ is in $D_{j-1}^{j_{\max}-j}\subseteq
D_{j-1}^1$ (since $j\leq j_{\max}-1$). It depends continuously on
$p_1$, $p_2$, $q_1$ and $h$ in this space, and can be extended to
$h=0$ continuously (such that $\tilde z(p_1,p_2,q_1,0)\in
D_{j-1}^{j_{\max}-j}$, too).

Hence, \eqref{eq:z2fixp} is a linear fixed-point problem for $z_2$
where the inhomogeneity is in $D_{j-1}^{j_{\max}-j}$ and depends
continuously on $(p,q,h)$. The linear map $M(h):z_2\mapsto
Q_NL\partial^1F_{j-1}(x_1+hz_1)\,z_2$ in front of $z_2$ on the
right-hand side of \eqref{eq:z2fixp} depends continuously on $h$ as an
element of $L(D_{j-1}^1;D_{j-1}^1)$ (see Lemma~\ref{thm:speccont} and
note that $x_1$ and $z_1$ are in $D_{j-1}^1$). Since the spectral
radius of the map $M(0)$ (for $h=0$) is less or equal than $1/2$ by
assumption of our inductive step, the spectral radius of $M(h)$ is
less than unity if we choose $h$ sufficiently small. Thus, for all
$p\in D(X_j)$ and $q\in \R^{n_j}$ and sufficiently small $h$, $z_2$
satisfies a contractive linear fixed point equation with an
inhomogeneity in $D_{j-1}^1$ and a contractive linear map that maps
into $D_{j-1}^1$ where all coefficients depend continuously on
$(h,p,q)$. Consequently, $z_2$ has a limit in $D_{j-1}^1$ for $h\to0$
that depends continuously on $(p,q)$. For $h=0$ the fixed point
equation for $(z_1,z_2)$ simplifies to
\begin{equation}
  \label{eq:isz}
  \begin{split}
    z_1&=E_Nq_1+Q_NL\partial^1F_{j-1}(x_1,z_1)\\
    z_2&=E_Nq_2+Q_NL\left[\partial^2F_{j-1}(x_1,x_2,z_1,0)+
      \partial^1F_{j-1}(x_1,z_2)\right]\mbox{.}    
  \end{split}
\end{equation}
Both equations are linear in $q$ and $z=(z_1,z_2)$. Consequently,
$z(0,p,q)$, which is by definition the directional derivative of $X_j$
in $p$ in direction $q$, depends linearly on $q$ and continuously on
$p$ and $q$. Consequently,
\begin{displaymath}
z(0,p,q)=\left[\frac{\partial}{\partial p}X_{j}(p)\right]q
\end{displaymath}
is the Frech{\'e}t derivative of $X_j$. 

\subsubsection*{Collection to finish proof of statement 1 of inductive step}
The functions $x=X_j(p)$ and $z=\partial^1X_j(p)q$ satisfy
\begin{equation}\label{eq:ivfjp1}
  \begin{aligned}
    x=&E_Np+Q_NLF_j(x) &&\mbox{by inductive assumption
      \eqref{eq:ivf}--\eqref{eq:ivderiv},}\\
    z=&E_Nq+Q_NL\partial^1F_j(x,z) &&\mbox{by \eqref{eq:isz} and
      definition of $F_j$.}
  \end{aligned}
\end{equation}
The variable $x=X_j(p)$ depends continuously on $p$ with respect to
the norm of $D_j^1$ by the assumption of the inductive step and the
step ``Regularity of $X_j(p)$''. The variable $z=\partial^1X(p)\,q$
depends continuously on $p$ and $q$ as shown in the previous step,
``Existence and continuity of $\partial^1X(p)\,q$''. Thus
$(x,z)=(X_j(p),\partial^1X(p)\,q)=X_{j+1}(p,q)\in D_j^1\times
D_j=D_{j+1}$ depends continuously on $(p,q)$, and satisfies
\eqref{eq:ivf}--\eqref{eq:ivderiv} for $j+1$ (which is identical to
system \eqref{eq:ivfjp1}). This completes the proof of
statement~\ref{p:ass:exfix0} of the inductive assumption for $j+1$.

\subsubsection*{Spectral radius of map $z\mapsto Q_NL\partial^1F_j(x,z)$}
The map $\partial^1F_j$ maps $D_{j+1}$ continuously into $D_j$ (by
Lemma~\ref{thm:rhs:smooth}). Thus, for fixed $x$ the linear map
$z\mapsto \partial^1F_j(x,z)$ maps $D_j$ continuously into $D_j$, and,
hence, the map $M_j(x): z\mapsto Q_NL\partial^1F_j(x,z)$ maps $D_j$
continuously into $D_j^1$, making $M_j(x)$ a compact linear
operator. Thus, the spectral radius of $M_j(x)$ is determined by its
largest eigenvalue (which has finite modulus and algebraic
multiplicity if it is non-zero). Splitting $M_j(x)$ into its two
components we get
\begin{displaymath}
  M_j(x):\begin{bmatrix}
    z_1\\ z_2
  \end{bmatrix}\mapsto
  \begin{bmatrix}
    \begin{aligned}
      &Q_NL\phantom{[\ }\partial^1F_{j-1}(x_1,z_1)\\
      &Q_NL\left[\partial^2F_{j-1}(x_1,x_2,z_1,0)+
        \partial^1F_{j-1}(x_1,z_2)\right]
    \end{aligned}
  \end{bmatrix}
\end{displaymath}
If $(\lambda,(z_1,z_2))$ is an eigenpair of $M_j(x)$ then the first
row of the definition of $M_j(x)$ implies that, either $(\lambda,z_1)$
is an eigenpair of $z_1\mapsto Q_NL\partial^1F_{j-1}(x_1,z_1)$, or
$z_1=0$. If $(\lambda,z_1)$ is an eigenpair of $z_1\mapsto
Q_NL\partial^1F_{j-1}(x_1,z_1)$ then, by inductive assumption,
$|\lambda|\leq 1/2$. If $z_1=0$ then the term
$\partial^2F_{j-1}(x_1,x_2,z_1,0)$ vanishes in the second row, such
that $(\lambda,z_2)$ is an eigenpair of $z_2\mapsto
Q_NL\partial^1F_{j-1}(x_1,z_2)$. Thus, by inductive assumption,
$|\lambda|\leq 1/2$ in this case, too. Consequently, the spectral
radius of $M_j(x)$ is also less or equal to $1/2$, which proves
statement~\ref{p:ass:specradius} of the inductive assumption for
$j+1$.

\subparagraph*{Existence of $Y_j$} We show inductively that
$Y_j(p)=F_j(X_j(p))$. For $j=1$ this statement was proven in
Lemma~\ref{thm:Xdiff}. Let $j<j_{\max}$ and assume that
$Y_j=F_j(X_j(p))$ for $p\in D(X_j)$. Since
\begin{displaymath}
  X_j(p)=E_Np+Q_NLF_j(X_j(p))
\end{displaymath}
and $F_j$ maps $D_j^1$ into $D_j^1$, $X_j$ is an element of
$D_j^1$. Let $q\in \R^{n_j}$ be arbitrary, and let us denote
$(x_1,x_2)=(X_j(p),\partial^1X_j(p)\,q)=X_{j+1}(p,q)$. The component $x_2$ satisfies
\begin{displaymath}
  x_2=E_Nq+Q_NL\partial^1F_j(x_1,x_2)
\end{displaymath}
such that $x_2$ is in $D_j^1$, too. Consequently,
\begin{align}
  \frac{Y_j(p+hq)-Y_j(p)}{h}&=
  \frac{F_j(X_j(p+hq))-F_j(X_j(p))}{h}\nonumber\\
  &=\frac{F_j(x_1+hx_2)-F_j(x_1)}{h}+
  \frac{F_j(X_j(p+hq))-F_j(x_1+hx_2)}{h}\mbox{.}\label{eq:Yjsplit}
\end{align}
Since $F_j$ is continuously differentiable for $x_1\in D_j^1$ and
deviations $hx_2\in D_j^1$ the first quotient in the
expression~\eqref{eq:Yjsplit} converges to
$\partial^1F_j(x_1,x_2)$. Since $F_j$ as a map from $D_j^1$ into $D_j$
is locally Lipschitz continuous the second term in \eqref{eq:Yjsplit}
can be bounded by
\begin{eqnarray*}
  \lefteqn{\left\| \frac{F_j(X_j(p+hq))-
        F_j(X_j(p)+h\partial^1X_j(p)q)}{h}\right\|_{D_j}\leq}\\
   &\leq& 
  K_1\left\|\frac{X_j(p+hq)-X_j(p)}{h}-\partial^1X_j(p)q\right\|_{D_j}\mbox{,}
\end{eqnarray*}
with some constant $K_1$, for sufficiently small $h$, which converges
to zero for $h\to0$ because $X_j$ is differentiable. Consequently, the
directional derivative of $Y_j$ in $p$ in direction $q$ is
$\partial^1F_j(X_j(p))[\partial X_j(p)\,q]$, which is continuous in
$p$ and $q$ and linear in $q$. Therefore, the Frech{\'e}t derivative
of $Y_j$ exists and
\begin{displaymath}
  \left[\frac{\partial}{\partial p}Y_j(p)\right]q=
  \partial^1F_j(X_j(p),\partial X_j(p)\,q)\mbox{,}
\end{displaymath}
which implies by definition of $F_j$ and $X_j$ that $Y_{j+1}=F_{j+1}(X_{j+1}(p,q))$.
\eop

We can refine the statement of Theorem~\ref{thm:smooth} slightly by
noting that $X_j:D(X_j)\mapsto D_j^1$ is continuous for all $j\leq
j_{\max}$ (instead of $X_j:D(X_j)\mapsto D_j$). This follows from
the continuity of $Y_j=F_j(X_j(p))$ as a map into $D_j$ and the
relation
\begin{displaymath}
  X_j(p)=E_Np+Q_NLY_j(p)\mbox{.}
\end{displaymath}
Theorem~\ref{thm:smooth} completes the proof of the Equivalence
Theorem~\ref{thm:main}. The refinement (that $X_j$ maps into $D_j^1$)
ensures that the image $X(p)$ is in $C^{j_{\max}+1}(\T;\R^n)$, as
claimed in Theorem~\ref{thm:main}

\section{Proof of Hopf Bifurcation Theorem}
\label{sec:hopf:proof}
First, we note that $x\mapsto S(x,\omega)^{-1}=x(\omega^{-1}\cdot)$
maps $C^k(\T;\R^n)$ into a closed subspace of $C^k([-\tau,0];\R^n)$,
if we extend functions $x$ on $\T$ to the whole real line by setting
$x(t)=x(t_{\mod[-\pi,\pi)})$. This implies that, if the functional
$f:C^0([-\tau,0];\R^n)\times\R\mapsto\R^n$ is $EC^k$ smooth then the
functional
\begin{displaymath}
  (x,\mu,\omega)\in C^0(\T;\R^n)\times\R^2\mapsto\frac{1}{\omega}
  f(S(x,\omega),\mu)\in\R^n
\end{displaymath}
is $EC^k$ smooth, too, such that we can reduce the problem of finding
periodic orbits of frequency $\omega$ to the algebraic system
\eqref{eq:lowmodes:par}. The right-hand side $F_y$ in
\eqref{eq:lowmodes:par} is defined by
\begin{displaymath}
  \left[F_y(x,\omega,\mu)\right](t)= 
  \frac{1}{\omega}f(S(\Delta_tx,\omega),\mu)\mbox{.}
\end{displaymath}
Let us choose the periodic orbit $x_0=(x,\omega,\mu)$ with $x=0$,
$\omega=\omega_0$, $\mu=0$ as the solution in the neighborhood of
which we construct the equivalent algebraic system. We choose the
number $N$ of Fourier modes and the size $\delta$ of the neighborhood
$B_\delta^{0,1}(x_0)$ in $C^{0,1}(\T;\R^{n+2})$ such that the
conditions of Theorem~\ref{thm:main} are satisfied in
$B_\delta^{0,1}(x_0)$. The full algebraic system
\eqref{eq:lowmodes:par} then reads (after multiplication with $\omega$
and mapping it onto the space $\rg P_N$ from $\R^{n(2N+1)}$ by
applying $R_N^{-1}$)
\begin{equation}
  \begin{aligned}
  0=&P_0F_y(X_y(p,\omega,\mu),\omega,\mu)+
  \omega Q_0P_NE_Np- Q_0P_NL F_y(X_y(p,\omega,\mu),\omega,\mu)    
  \end{aligned}
  \label{eq:hopf:modes}
\end{equation}
The variables are $p\in\R^{n(2N+1)}$ (which was called $p_y$ in
\eqref{eq:lowmodes:par}), $\mu$ and $\omega$. We know from
Theorem~\ref{thm:main} that
\begin{align*}
  Y_y:&(p,\omega,\mu)\in\R^{n(2N+1)}\times\R\times\R
  \mapsto F(X_y(p,\omega,\mu),\omega),\mu)\in C^0(\T;\R^n)\mbox{,}\\
  X_y:&(p,\mu,\omega)\in\R^{n(2N+1)}\times\R\times\R\mapsto X_y(p,\omega,\mu)
  \in C^0(\T;\R^n)
\end{align*}
are $k$ times differentiable, and note that
\begin{equation}\label{eq:hopf:fzero}
  F_y(X_y(0,\omega,\mu),\omega),\mu)=0  
\end{equation}
for all $\omega\approx\omega_0$ and $\mu\approx0$ (because
$x_0=(0,\omega,\mu)$ is a solution). The derivative of the right-hand
side $F_y$ in $x=0$, $\omega\approx\omega_0$ and $\mu\approx0$ with
respect to $x$ is $A(\omega,\mu)x$, defined by
\begin{align*}
  \left[A(\omega,\mu)x\right](t)=a(\mu)\left[x(t+\omega\cdot)\right]\mbox{,}
\end{align*}
where $a(\mu)$ is the same linear functional as used in the definition
of the characteristic matrix $K(\lambda,\mu)$ in \eqref{eq:charmdef}
(the derivatives of $F$ with respect to $\omega$ and $\mu$ are zero
due to \eqref{eq:hopf:fzero}). We observe that $A(\omega,\mu)$
commutes with $P_j$ and $Q_j$ for all $j\geq0$.

Let us now determine the linearization of $X_y(p,\omega,\mu)$ in
$(p,\omega,\mu)=(0,\omega,\mu)$. Due to \eqref{eq:hopf:fzero}
$X_y(0,\omega,\mu)$ is equal to zero for all $\omega\approx\omega_0$
and $\mu\approx0$: since $0$ is a solution to the periodic BVP and
$P_N0=0$, the zero solution must also be equal to
$X_y(0,\omega,\mu)$. Thus, we have
\begin{align*}
  0&=\left.\frac{\partial}{\partial\omega} X_y(p,\omega,\mu)
  \right\vert_{\textstyle p=0}\mbox{\quad and} &
  0&=\left.\frac{\partial}{\partial\mu}
    X_y(p,\omega,\mu)\right\vert_{\textstyle p=0}\mbox{.}
\end{align*}
Moreover, the fixed point equation \eqref{eq:dxfixp} defining
$z=[\partial X_y/\partial p] (p,\omega,\mu)\,q$, evaluated in
$(p,\omega,\mu)=(0,\omega,\mu)$ reads
\begin{equation}\label{eq:hopf:qnz}
  z=E_Nq+Q_NLA(\mu,\omega)\,z=
  E_Nq+Q_NLA(\mu,\omega)\,Q_Nz\mbox{,}
\end{equation}
exploiting that $Q_NL=Q_NLQ_N$ and
$Q_NA(\omega,\mu)=A(\omega,\mu)\,Q_N$. In the neighborhood
$B_\delta^{0,1}(x_0)$ the spectral radius of $Q_NLA(\mu,\omega)$ is
less than unity (see Lemma~\ref{thm:specradius}). Application of $Q_N$
to \eqref{eq:hopf:qnz} gives $Q_Nz=Q_NLA(\mu,\omega)Q_Nz$. Since the
spectral radius of $Q_NLA(\mu,\omega)$ is less than unity this implies
that $Q_Nz=0$, and, thus
\begin{displaymath}
  \left[\left.\frac{\partial}{\partial p}
    X_y(p,\omega,\mu)\right\vert_{\textstyle p=0}\right]q= E_Nq\mbox{.}
\end{displaymath}

Consequently, the linearization of the algebraic system
\eqref{eq:hopf:modes} in $(0,\omega,\mu)$ with respect
to the first variable is
\begin{equation}\label{eq:hopf:lin}
  0=P_0A(\omega,\mu)E_Np+\omega Q_0P_NE_Np-Q_0P_NLA(\omega,\mu)E_Np
\end{equation}
for all $\omega\approx \omega_0$ and $\mu\approx 0$ (also using $p$
for the argument of the linearization in \eqref{eq:hopf:lin}). We
observe that the linear system \eqref{eq:hopf:lin} decouples into
equations for
\begin{align*}
  y_0&=P_oE_Np=E_0p=p_0 &&\mbox{(the average of $E_Np$),}\\
  y_1&=Q_0E_1p=p_{-1}\sin t+p_1\cos t &&\mbox{(the first Fourier component of $E_Np$),}\\
  y_j&=Q_{j-1}E_jp=p_{-j}\sin(jt)+p_j\cos(jt) &&\mbox{(the $j$-th
    Fourier component of $E_Np$,}\\
  &&&\mbox{$2\leq j\leq N$),}
\end{align*}
where we denote the components of $p$ by $p_j\in \R^n$ ($j=-N\ldots N$).
This decoupling is achieved by pre-multiplication of
\eqref{eq:hopf:lin} with $P_0$ and $Q_{j-1}P_j$ for $j=1\ldots N$:
\begin{align}
    P_0\cdot\mbox{\eqref{eq:hopf:lin}}:&&
    0&=A(\omega,\mu)y_0=a(\mu)\,p_0
      \label{eq:hopf:lin:dec0}\\
    Q_0P_1\cdot\mbox{\eqref{eq:hopf:lin}}:&&
    0&=\omega y_1-Q_0LA(\omega,\mu)\,y_1    
    \label{eq:hopf:lin:dec1}\\
    Q_{j-1}P_j\cdot\mbox{\eqref{eq:hopf:lin}}:&&
    0&=\omega y_j-Q_0LA(\omega,\mu)\,y_j &&
    \mbox{\ for $j=2\ldots N$.}
    \label{eq:hopf:lin:decj}
\end{align}
Inserting the definition of $y_j$ into the equations~\eqref{eq:hopf:lin:dec1}
and \eqref{eq:hopf:lin:decj} gives for $j\geq 1$
\begin{align*}
  0=&\omega\left[p_{-j}\sin(jt)+p_j\cos(jt)\right]-
  Q_0\int_0^ta(\mu)[p_j\sin(js+j\omega\cdot)+p_j\cos(js+j\omega\cdot)]\d s\\
  =&\omega\left[p_{-j}\sin(jt)+p_j\cos(jt)\right]-
  \frac{1}{j}\sin(jt)a(\mu)[p_{-j}\sin(j\omega\cdot)+p_j\cos(j\omega\cdot)]\\
  &\phantom{x}-\frac{1}{j}\cos(jt)a(\mu)
  [-p_{-j}\cos(j\omega\cdot)+p_j\sin(j\omega\cdot)]\mbox{.}
\end{align*}
These equations are satisfied if and only if the coefficients in front
of $\sin(jt)$ and $\cos(jt)$ are zero. The resulting system of
equations reads in complex notation (splitting up again into the cases
$j=1$ and $j>1$)
\begin{align}
  \label{eq:hopf:red1}
  i\omega u_1-a(\mu)\left[u_1\exp(i\omega s)\right]&=K(i\omega,\mu)\,u_1=0\mbox{,}\\
  \label{eq:hopf:redj}
  ij\omega u_j-a(\mu)\left[u_j\exp(ij\omega s)\right]&=K(ij\omega,\mu)\,u_j=0
  \mbox{\quad ($2\leq j\leq N$),}
\end{align}
that is, $u_j=p_{-j}+ip_j\in\C^n$ is a solution of
\eqref{eq:hopf:red1} (or \eqref{eq:hopf:redj}, respectively) if and
only if $y_j=p_{-j}\sin(jt)+p_j\cos(jt)$ is a solution of
\eqref{eq:hopf:lin:dec1} (or \eqref{eq:hopf:lin:decj}, respectively).

The non-resonance assumption of the theorem guarantees that equation
\eqref{eq:hopf:lin:dec0} is a regular linear system for $p_0$, and
that \eqref{eq:hopf:redj} is a regular linear algebraic system for
$p_{-j}$ and $p_j$ ($j\geq 2$) at $\mu=0$ and $\omega=\omega_0$ (and,
hence, for all $\omega$ and $\mu$ near-by).  The condition on the
simplicity of the eigenvalue $i\omega_0$ of $K$ ensures that equation
\eqref{eq:hopf:red1} (and, thus, \eqref{eq:hopf:lin:dec1}) has a
one-dimensional (in complex notation) subspace of solutions for
$\omega=\omega_0$ and $\mu=0$, spanned by the nullvector $v_1$ of
$K(i\omega,0)$.  Let us denote the adjoint nullvector of
$K(i\omega_0,0)$ by $w_1\in\C^n$ (again, using complex notation,
$w_1^HK(i\omega_0,0)=0$). Since $i\omega_0$ is simple, the
relationship
\begin{displaymath}
  w_1^H\frac{\partial K}{\partial\lambda}(i\omega,0)\,v_1\neq0
\end{displaymath}
holds\mbox{,} which implies that we can choose $w_1\in\C^n$ without
loss of generality such that
\begin{displaymath}
  w_1^H\frac{\partial K}{\partial\lambda}(i\omega,0)\,v_1=1\mbox{.}  
\end{displaymath}
With this convention we observe that
\begin{align}
  w_1^H\frac{\partial K}{\partial\mu}(i\omega,0)\,v_1&=
  -\left.\frac{\partial\lambda}{\partial\mu}\right\vert_{\textstyle\mu=0}
  =:c_\mu\in\C\mbox{, and}&
  w_1^H\frac{\partial}{\partial\omega}K(i\omega,0)\,v_1&=
  i\in\C\label{eq:hopf:derivs}
\end{align}
where $\Re c_\mu\neq0$ by the transversal crossing assumption of the
theorem.  In complex notation any scalar multiple of the nullvector
$v_1=v_r+iv_i$ is also a nullvector. Thus, the complex scalar factor
$\alpha+i\beta$ in front of $v_1$ makes up two components of the
variable $p$ (in real notation): in short, $p$ solves the linearized
algebraic system \eqref{eq:hopf:lin} if and only if all $p_j$ with
$|j|\neq1$ are zero and $p_{-1}\sin t+p_1\cos
t=\Re\left[(\alpha+i\beta)v_1\exp(it)\right]$ for some
$\alpha,\beta\in\R$, that is,
\begin{equation}\label{eq:hopf:abintro}
  \begin{bmatrix}
   p_{-1}\\ p_{1\phantom{-}} 
  \end{bmatrix}=
  \alpha 
  \begin{bmatrix}
    -v_i\\ \phantom{-}v_r
  \end{bmatrix}+\beta
  \begin{bmatrix}
    -v_r\\ -v_i
  \end{bmatrix}=:\alpha b_r+\beta b_i\mbox{.}
\end{equation}
Let us collect the statements so far and introduce coordinates.  We
collect all components $p_j$ with $|j|\neq1$ and the orthogonal
complement in $\R^{2n}$ of the space spanned by $\{b_1,b_2\}$ into a
single variable $q$ (of real dimension $n_q=n(2N-1)+2(n-1)$). Then a
set of coordinates for $p$ are the variables
\begin{align*}
  (\alpha,\beta)&=:r\in\R^2\mbox{,\quad and\quad}  
  q\in\R^{n_q}\mbox{.}
\end{align*}
We split up the full algebraic system of equations
\eqref{eq:hopf:modes} in the same way as we did for the linearized
problem, by pre-multiplication with $P_0$ and $Q_{j-1}P_j$ for
$j=1\ldots N$:
\begin{align}
    P_0\cdot\mbox{\eqref{eq:hopf:modes}}:&&
     0&=P_0F(X_y(p,\omega,\mu),\omega,\mu)
      \label{eq:hopf:nlin0}\\
    Q_0P_1\cdot\mbox{\eqref{eq:hopf:modes}}:&&
    0&=\omega Q_0E_1p-Q_0P_1LF(X_y(p,\omega,\mu),\omega,\mu)
    \label{eq:hopf:nlin1}\\
    Q_{j-1}P_j\cdot\mbox{\eqref{eq:hopf:modes}}:&& 0&=\omega
    Q_{j-1}E_jp-Q_{j-1}P_jLF(X_y(p,\omega,\mu),\omega,\mu)\mbox{.}
    \label{eq:hopf:nlinj}
\end{align}
We split equation~\eqref{eq:hopf:nlin1} further using $w_1^H$ and its
orthogonal complement, the projection $w_1^\perp=\id-w_1w_1^H/(w_1^Hw_1)$. This gives
rise to a splitting into two real equations
($w_1^H\cdot$\eqref{eq:hopf:nlin1}) and $2(n-1)$ real equations
($w_1^\perp\cdot$\eqref{eq:hopf:nlin1}). Collecting
$w_1^\perp\cdot$\eqref{eq:hopf:nlin1} and the equations
\eqref{eq:hopf:nlin0} and \eqref{eq:hopf:nlinj} into a subsystem of
$n(2N-1)+2(n-1)=n_q$ equations the full algebraic system
\eqref{eq:hopf:modes} in the coordinates $(r,q)$ has the form
\begin{equation}\label{eq:hopf:nlin:matrixform}
  0=\begin{bmatrix}
    M_{rr}(r,q,\omega,\mu) & M_{rq}(r,q,\omega,\mu)\\
    M_{qr}(r,q,\omega,\mu) & M_{qq}(r,q,\omega,\mu)
  \end{bmatrix}
  \begin{bmatrix}
    r\\ q
  \end{bmatrix}\mbox{.}
\end{equation}
By our choice of coordinates the matrices $M_{rr}\in\R^{2\times2}$,
$M_{rq}\in\R^{2\times n_q}$ and $M_{qr}\in\R^{n_q\times2}$ are
identically zero in $r=0$, $q=0$, $\mu=0$, $\omega=i\omega_0$ such
that the system matrix has the form
\begin{displaymath}
  \begin{bmatrix}
    \begin{bmatrix}
      0&0\\ 0&0
    \end{bmatrix} &
    \begin{bmatrix}
      0 & \ldots & 0\\
      0 & \ldots & 0
    \end{bmatrix}\\
    \begin{bmatrix}
      0&0\\
      \vdots&\vdots\\
      0&0
    \end{bmatrix} &
    \begin{matrix}
    M_{qq}(0,0,i\omega,0)\\
    \mbox{(regular)}      
    \end{matrix}
  \end{bmatrix}
\end{displaymath}
$(r,q,\mu,\omega)=(0,0,0,\omega_0)$. Thus, we can perform a
Lyapunov-Schmidt reduction: we eliminate $q$ by solving the $n_q$
lower equations for $q$, obtaining a graph $q(r,\omega,\mu)\,r$
locally in a neighborhood of $(r,q,\mu,\omega)=(0,0,0,\omega_0)$. This
graph respects rotational invariance:
$q(\Delta_sr,\omega,\mu)\,\Delta_sr=\Delta_s[q(r,\omega,\mu)\,r]$.
Note that the application of $\Delta_s$ to $r=(\alpha,\beta)$
corresponds to the rotation of $r$ by angle $s$ (the same as the
multiplication $\exp(is)(\alpha+i\beta)$).  The Lyapunov-Schmidt
reduction of \eqref{eq:hopf:nlin:matrixform}, replacing $q$ by the
graph $q(r,\omega,\mu)\,r$, then reads
\begin{equation}
  \label{eq:hopf:nlin:reduced}
  0=M_{rr}(r,q(r,\omega,\mu)\,r,\omega,\mu)\,r=:M_r(r,\omega,\mu)\,r\mbox{,}
\end{equation}
where $M_r$ is still rotationally symmetric in $r$:
$M_r(\Delta_sr,\omega,\mu)\,\Delta_sr=\Delta_sM_r(r,\omega,\mu)\,r$.
Equation~\eqref{eq:hopf:derivs} in real notation implies that
\begin{align*}
  \frac{\partial M_r}{\partial\omega}(0,\omega_0,0)&=
  \frac{\partial M_{rr}}{\partial\omega}(0,0,\omega_0,0)=
  \begin{bmatrix}
    0&{-1}\\ 1 &\phantom{-}0
  \end{bmatrix}\mbox{,} \\
  \frac{\partial M_r}{\partial\mu}(0,\omega_0,0)&=
  \frac{\partial M_{rr}}{\partial\mu}(0,0,\omega_0,0)=
  \begin{bmatrix}
    \Re c_\mu&-\Im c_\mu\\ \Im c_\mu &\phantom{-}\Re c_\mu
  \end{bmatrix}\mbox{.}
\end{align*}
Equation~\eqref{eq:hopf:nlin:reduced} is a system of two equations
with four unknowns ($r=(\alpha,\beta)$, $\omega$ and $\mu$). We now
fix one of the unknowns setting
\begin{displaymath}
  \alpha=0
\end{displaymath}
such that we can expect one-parametric families of solutions
$(\beta,\omega,\mu)$. 
Introducing $M_\beta$ as the second column of $M_r$ and dropping the
dependence on $\alpha$ (which is zero), the first derivative of
$M_\beta(\beta,\omega,\mu)$ in $(0,\omega_0,0)$ with respect to the
pair $\omega$ and $\mu$ is:
\begin{align*}
  \begin{bmatrix}
  {\displaystyle\frac{\partial M_\beta}{\partial \omega}}&
  {\displaystyle\frac{\partial M_\beta}{\partial \mu}}
  \end{bmatrix}
  (0,\omega_0,0)=    
  \begin{bmatrix}
    -1 &-\Im c_\mu\\ \phantom{-}0& \phantom{-}\Re c_\mu
  \end{bmatrix}\mbox{,}
\end{align*}
which is regular (as established in \eqref{eq:hopf:derivs}, since $\Re
c_\mu\neq0$ due to the assumption that the eigenvalue crosses the
imaginary axis transversally).  Note that $M_\beta$ itself is a
projection of the first derivative of the original right-hand side of
the full algebraic system \eqref{eq:hopf:modes}. Thus,
$M_\beta$ is $k-1$ times continuously differentiable, and we end up
with a system of two equations for three scalar variables
$(\beta,\omega,\mu)$:
\begin{align*}
  0=M_\beta(\beta,\omega,\mu)\,\beta\mbox{.}
\end{align*}
Hence, either $\beta=0$, which corresponds to the trivial solution or
(after division by $\beta$)
\begin{align}
  0=&M_\beta(\beta,\omega,\mu)\mbox{,}\label{eq:hopf:nlin:factor}
\end{align}
where $M_\beta(0,\omega_0,0)=(0,0)$ and the derivative with respect to
the pair $(\omega,\mu)$ is regular in $(0,\omega_0,0)$.  Thus, we can
apply the Implicit Function Theorem to \eqref{eq:hopf:nlin:factor} to
obtain a unique graph $(\omega(\beta),\mu(\beta))$ solving
\eqref{eq:hopf:nlin:factor}. The graph satisfies
$(\omega(0),\mu(0))=(\omega_0,0)$, and, thus, branches off from the
trivial solution (which has $\beta=0$ and $\omega$ and $\mu$
arbitrary). The rotational symmetry of $M_r$ implies reflection
symmetry of $M_\beta$ in $\beta$ such that
$M_\beta(-\beta,\omega,\mu)=M_\beta(\beta,\omega,\mu)$ for all
$\beta$, $\omega$ and $\mu$. Hence, the solution graph is reflection
symmetric, too: $\omega(-\beta)=\omega(\beta)$ and
$\mu(-\beta)=\mu(\beta)$. Thus, for small $\beta$ there is a unique
non-trivial solution of the full algebraic system of the form
$r=(0,\beta)$, $q=q(r,\omega(\beta),\mu(\beta))\,r$. As
Equation~\eqref{eq:hopf:abintro} shows, we can extract the coordinates
$\alpha$ (which is zero) and $\beta$ from the full solution $x\in
C^k(\T;\R^n)$ by applying the projections
\begin{align*}
  \frac{1}{\pi}\int_{-\pi}^\pi\cos(t)v_r^Tx(t)-\sin(t)v_i^Tx(t)\d t&=
  \frac{1}{\pi}\int_{-\pi}^\pi\Re\left[v_1\exp(it)\right]^Tx(t)\d t=\alpha\mbox{,}\\
  \frac{1}{\pi}\int_{-\pi}^\pi\sin(t)v_r^Tx(t)+\cos(t)v_i^Tx(t)\d t&=
  \frac{1}{\pi}\int_{-\pi}^\pi\Im\left[v_1\exp(it)\right]^Tx(t)\d t=-\beta\mbox{,}
\end{align*}
which determines the First Fourier coefficients of
$x$ as claimed in \eqref{eq:hopfphase} in
Theorem~\ref{thm:hopf}. (Recall that the vector $v_1=v_r+v_i$ was
scaled to have unit length and that the decomposition was
orthogonal.) 
\eop

\bibliographystyle{plain}
\bibliography{delay}

\appendix


\section{Basic differentiability properties of the right-hand side}
\label{sec:basicprop}
Let $J$ be a compact interval or $\T$. Let $(D,\|\cdot\|_D)$ be a
Banach space of the form
\begin{displaymath}
D=C^{k_1}(J;\R^{m_1})\times\ldots\times C^{k_\ell}(J;\R^{m_\ell})  
\end{displaymath}
where $\ell\geq 1$, the integers $k_j$ are non-negative and the integers
$m_j$ are positive. We use the natural maximum norm on the product $D$:
\begin{displaymath}
\|x\|_D=\|(x_1,\ldots,x_\ell)\|_D=\max_{j\in\{1,\ldots,\ell\}}\|x_j\|_{k_j}\mbox{,}
\end{displaymath} and use the notation
\begin{align*}
  D^k&=C^{k_1+k}(J;\R^{m_1})\times\ldots\times C^{k_\ell+k}(J;\R^{m_\ell})\mbox{,}
   & \|x\|_{D,k}&=\max_{0\leq j\leq
    k} \|x^{(j)}\|_D\mbox{,}\\
    D^{0,1}&=\left\{x\in D: L(x)<\infty\right\}\mbox{, with the norm}&
    \|x\|_{D,L}&=\max\left\{\|x\|_D,L(x)\right\}\mbox{,}
    \intertext{where $x^{(j)}$ is the component-wise $j$th derivative
      and the Lipschitz constant $L(x)$ is defined as} L(x)&=\sup_{
      \begin{subarray}{c}
        t\neq s
      \end{subarray}
    }\ \max_{j=1\ldots,\ell}\ \frac{|x_j^{(k_j)}(t)-x_j^{(k_j)}(s)|}{|t-s|}
    \mbox{,}
\end{align*}
where $t$ and $s$ in the index of $\sup$ are taken from $J$, if $J$ is
a compact interval, and from $\R$ if $J=\T$. Balls that are closed and
bounded in $D^{0,1}$ are complete with respect to the norm of $D$.  

\subsection{Basic properties of $f$}
\label{sec:fprop}
This section proves three properties that $EC^1$ smooth
functionals $f$ have: first that the derivative limit
\eqref{eq:ass:contdiff:j} exists also for Lipschitz continuous
deviations, second a weaker form of the mean value theorem, and third,
local $EC$ Lipschitz continuity.
\begin{lemma}[Extension of derivative to 
  deviations in $D^{0,1}$]\label{thm:weakcontdiff}
  Let $f:D\mapsto\R^n$ be $EC^1$ smooth in the sense of
  Definition~\ref{def:extdiff}.  Then the limit required to exist in
  Definition\ref{def:extdiff} exists also in the
  $\|\cdot\|_{D,L}$-norm: for all $x\in D^1$
  \begin{align}\allowdisplaybreaks
     \label{eq:ass:contdifflip} 
      \lim_{
        \begin{subarray}{c}
          y\in D^{0,1}\\[0.2ex]
          \|y\|_{D,L}\to0
        \end{subarray}
      }&
      \frac{|f(x+y)-f(x)-\partial^1f(x,y)|}{\|y\|_{D,L}}=0\mbox{.}
  \end{align}
\end{lemma}
Note that in \eqref{eq:ass:contdifflip} the norm in which $y$ goes to
zero is $\|\cdot\|_{D,L}$ instead of $\|\cdot\|_{D,1}$.

\pstart{Proof} This is a simple continuity argument. Let $\epsilon>0$
be arbitrary. We pick $\delta>0$ such that
\begin{equation}\label{eq:ydsmall}
  |f(x+\tilde y)-f(x)-\partial^1f(x,\tilde y)|<\epsilon\|\tilde y\|_{D,1}
\end{equation}
for all $\tilde y\in D^1$ satisfying $\|\tilde y\|_{D,1}<\delta$. Let
$y\in D^{0,1}$ be such that $\|y\|_{D,L}<\delta$. For every $\rho>0$
we can find a $\tilde y\in D^1$ such that $\|\tilde y\|_{D,1}\leq
\|y\|_{D,L}$ and $\|\tilde y-y\|_D\leq \rho$ (for example $\tilde
y=F_n*y$ with a Fej{\'e}r kernel with sufficiently large $n$ will meet
this criterion).  We can choose our $\rho$ sufficiently small such
that the $\tilde y\in D^1$ thus constructed satisfies
\begin{align}
  \|\tilde y\|_{D,1}&\leq \|y\|_{D,L}\label{eq:yydsmall}\\
  |f(x+y)-f(x+\tilde y)|&<\epsilon\|y\|_{D,L}\label{eq:fyydsmall}\\
  |\partial^1f(x,y-\tilde y)|&<\epsilon\|y\|_{D,L}\label{eq:ayydsmall}\mbox{.}
\end{align}
Condition \eqref{eq:fyydsmall} follows from the continuity of $f$ and
the density of $D^{0,1}$ in $D^1$ (both with respect to the $D$-norm),
and \eqref{eq:ayydsmall} follows from the continuity of $\partial^1f$
as a map on $D^1\times D$, and the density of $D^{0,1}$ in $D^1$ with
respect to $D$-norm. Combining estimate \eqref{eq:ydsmall} with
\eqref{eq:yydsmall}--\eqref{eq:ayydsmall} we obtain
\begin{displaymath}
  |f(x+y)-f(x)-\partial^1f(x,y)|<3\epsilon\|y\|_{D,L}\mbox{.}
\end{displaymath}
\eop

\begin{lemma}[Existence of mean value]\label{thm:fmeanval}
  There exists a continuous function
  \begin{displaymath}
    \tilde a:D^1\times D^1\times D\mapsto \R^n
  \end{displaymath}
  which is linear in its third argument and satisfies for all $x,y\in
  D^1$
  \begin{equation}\label{eq:meanval}
    f(x+y)-f(x)=\tilde a(x,y,y)\mbox{.}
  \end{equation}
  Moreover, $\tilde a(x,0,y)=\partial^1f(x,y)$ for all $x\in D^1$ and
  $y\in D$.
\end{lemma}

\pstart{Proof}
The argument for the existence of a mean value follows exactly the
proof of the general mean value theorem \cite{HKWW06}: the candidate
for $\tilde a(u,v,w)$ is
\begin{equation}\label{eq:meandiff}
  \tilde a(u,v,w)=\int_0^1\partial^1f(u+sv,w)\d s\mbox{.}
\end{equation}
Since $\partial^1f$ is assumed to be continuous in its arguments the
integral is well defined and continuous in its arguments $u\in D^1$,
$v\in D^1$, $w\in D$. All one has to show is that the $\tilde a$
defined in \eqref{eq:meandiff} satisfies \eqref{eq:meanval}: let $x,y
\in D^1$ and $\epsilon>0$ be arbitrary, and choose $m$ such that
uniformly for all $s\in[0,1]$
\begin{align*}
  \left|\int_0^1\partial^1f(x+sy,y)\d s- \frac{1}{m}\sum_{k=0}^{m-1}
    \partial^1f\left(x+\frac{k}{m}y,y\right)\right|&<\epsilon\mbox{,}\\
  \left|f\left(x+sy+\frac{y}{m}\right)
    -f(x+sy)-\partial^1f\left(x+sy,\frac{y}{m}\right)\right|&<\frac{\epsilon}{m}\mbox{.}
\end{align*}
Then it follows that
\begin{displaymath}
  \left|f(x+y)-f(x)-\int_0^1\partial^1f(x+sy,y)\d s\right|<2\epsilon\mbox{.}
\end{displaymath}
Since $\epsilon>0$ was arbitrary the left-hand side must be zero. \eop

\begin{lemma}[Local $EC$ Lipschitz continuity]\label{thm:flip}  
  For all $x\in D^1$ there exists a neighborhood $U(x)\subseteq
  D^{0,1}$ and a constant $K_x>0$ such that for all $y_1$ and $y_2\in
  U(x)$ the following Lipschitz estimate holds:
  \begin{displaymath}
    |f(y_1)-f(y_2)|\leq K_x\|y_1-y_2\|_D\mbox{.}
  \end{displaymath}
\end{lemma}
Note that the upper bound depends only on the $\|\cdot\|_D$-norm, not
on the $\|\cdot\|_{D,L}$-norm, which would be a weaker statement.

\pstart{Proof}
Let $x$ be an element of $D^1$. We prove the Lipschitz continuity first for $y_1$ and $y_2$ from a
ball $B_{\delta/2}^1(x)\subseteq D^1$ around $x\in
D^1$ with sufficiently small $\delta$. 

 Since the mean value $\tilde a$ is
continuous in $(x,0,0)$, and $\tilde a(x,0,0)=0$, we have a $\delta>0$
such that for all $u,v\in D^1$ and $w\in D$ satisfying
$\|u\|_{D,1}<\delta$, $\|v\|_{D,1}<\delta$ and $\|w\|_D<\delta$
\begin{displaymath}
  |\tilde a(x+u,v,w)|<\epsilon\mbox{.}
\end{displaymath}
This implies that $|\tilde a(x+u,v,w)|<[\epsilon/\delta]\|w\|_D$ for
$u$ and $v$ with $\max\{\|u\|_{D,1},\|v\|_{D,1}\}<\delta$ and $w\in D$
(since $\tilde a$ is linear in its third argument). Thus, $\|\tilde
a(x+u,v,\cdot)\|_D\leq\epsilon/\delta$ for $\tilde
a(x+u,v,\cdot)$ as an element of $L(D;D)$ in the operator norm
corresponding to $D$. Consequently, if $\|y_1-x\|_{D,1}<\delta/2$ and
$\|y_2-x\|_{D,1}<\delta/2$
\begin{displaymath}
  |f(y_1)-f(y_2)|=\left|\int_0^1\tilde a(y_2,y_1-y_2,y_1-y_2)\d s\right|\leq 
  \frac{\epsilon}{\delta}\|y_1-y_2\|_D\mbox{,}
\end{displaymath}
such that we can choose $K_x=\epsilon/\delta$.  The extension of the
statement to $D^{0,1}$ follows from the continuity of $f$ in $D$:
$B_{\delta/2}^1(x)\subset D^1$ is dense in
$B_{\delta/2}^{0,1}(x)\subset D^{0,1}$ when using the
$\|\cdot\|_D$-norm. Pick two sequences $y_n$ and $z_n$ in
$B_{\delta/2}^1(x)$ that converge to $y$ and $z$ in
$B_{\delta/2}^{0,1}(x)$ in the $D$-norm, that is, $\|y_n-y\|_D\to0$,
$\|z_n-z\|_D\to 0$ for $n\to\infty$. Then $f(y_n)\to f(y)$ and
$f(z_n)\to f(z)$ since $f$ is continuous in $D$. Moreover,
$\|y_n-z_n\|_D\to\|y-z\|_D$ for $n\to\infty$. Since
\begin{equation} |f(y_n)-f(z_n)|\leq
  K_x\|y_n-z_n\|_D
\end{equation}
for all $n$ the inequality also holds for the limit for $n\to\infty$.
\eop

\subsection{Basic properties of $F$}
\label{sec:Fprop}
In this section we restrict ourselves to the periodic case:
$J=\T$. Let $F:D\mapsto C^0(\T;\R^n)$ be defined as
$F(x)(t)=f(\Delta_tx)$.
\begin{lemma}[Continuity of $F$]\label{thm:Fcont}
  Let $f:D\mapsto\R^n$ be continuous. Then $F:D\mapsto C^0(\T;\R^n)$
  is also continuous.
\end{lemma}
\pstart{Proof} This is a simple consequence of the continuity of $f$,
the continuity of $(t,x)\mapsto \Delta_tx$ with respect to both
arguments ($t$ and $x$) in the $\|\cdot\|_0$-norm, and the
compactness of $\T$. Let $\epsilon>0$ and $x\in D$ be arbitrary. We
want to prove continuity of $F$ in $x$. So, we have to find a
$\delta>0$ such that
\begin{equation}\label{eq:Fcont:epsdelta}
  \left|f(\Delta_sx+h)-f(\Delta_sx)\right|<\epsilon 
  \mbox{\quad for all $s\in\T$ and $h\in D$, satisfying $\|h\|_D<\delta$.}
\end{equation}
(Since $\|\Delta_sh\|_D=\|h\|_D$ we can replace $\Delta_sh$ by $h$.)
The continuity of $f$ implies that for every $r>0$ and every $t\in\T$
we find a $\delta_x(t,r)$ such that
\begin{equation}
  |f(\Delta_tx+h)-f(\Delta_tx)|<r
  \mbox{\quad whenever $\|h\|_D<\delta_x(t,r)$.}\label{eq:Fcont:fdelta}
\end{equation}
For every $t\in\T$ there exists an
open neighborhood $U(t)\subset \T$ such that
\begin{displaymath}
  \|\Delta_sx-\Delta_tx\|_D<\delta_x(t,\epsilon/2)/2
  \mbox{\quad for all $s\in U(t)$,}
\end{displaymath}
because the function $t\in\T\mapsto \Delta_tx$ is continuous in $t$.
These neighborhoods $U(t)$ are an open cover of the compact set $\T$,
so there exist finitely many $t_1,\ldots,t_m\in\T$ such that the union
of the neighborhoods $U(t_j)$ contains all points $s\in\T$. We choose
\begin{displaymath}
  \delta=\min_{j=1,\ldots,m}\delta_x(t_j,\epsilon/2)/2\mbox{,}
\end{displaymath}
which is a positive quantity. Let $s\in\T$ be arbitrary and let $h\in
D$ satisfy $\|h\|_D<\delta$. We have to check the inequality
\eqref{eq:Fcont:epsdelta}. The point $s$ is in one of the
neighborhoods $U(t_j)$, say without loss of generality, $s\in
U(t_1)$. Thus,
$\|\Delta_sx-\Delta_{t_1}x\|_D<\delta_x(t_1,\epsilon/2)/2$, and,
consequently,
$\|\Delta_sx-\Delta_{t_1}x+h\|_D<\delta_x(t_1,\epsilon/2)$ (because
also $\|h\|_D<\delta\leq \delta_x(t_1,/\epsilon/2)/2$). Therefore, we
can split up the difference $|f(\Delta_sx+h)-f(\Delta_sx)|$:
\begin{align*}
  |f(\Delta_sx+h)-f(\Delta_sx)|\leq&\
  \left|\left[f\left(\Delta_{t_1}x+(\Delta_sx-\Delta_{t_1}x+h)\right)
      -f(\Delta_{t_1}x)\right]\right|\\
   &\ +\left|\left[f\left(\Delta_{t_1}x+(\Delta_sx-\Delta_{t_1}x)\right)
      -f(\Delta_{t_1}x)\right]\right|\\
  <&\ \epsilon/2+\epsilon/2=\epsilon
\end{align*}
Note that the deviations from $\Delta_{t_1}x$ in the arguments of $f$
in both terms of the sum are less than or equal to
$\delta_x(t_1,\epsilon/2)$ such that we can apply
\eqref{eq:Fcont:fdelta} for $t=t_1$, $r=\epsilon/2$. \eop

The following lemma lists properties that $F$ has if $f$ satisfies
local $EC$ Lipschitz continuity in the sense of
Definition~\ref{def:loclip}. That is, we do \emph{not} assume that $f$
is $EC^1$ smooth in the sense of Definition~\ref{def:extdiff} for
Lemma~\ref{thm:Flipbound}. Since Lemma~\ref{thm:flip} was proved using
only the assumption of $EC^1$ smoothness of $f$, local $EC$ Lipschitz
continuity is a weaker condition.
\begin{lemma}[$EC$ Lipschitz continuity of $F$]\label{thm:Flipbound}
  Assume that $f:D\mapsto\R^n$ is locally $EC$ Lipschitz continuous in the sense of
  Definition~\ref{def:loclip}. Then $F$ has the following properties:
  \begin{enumerate}
  \item \label{thm:Flip} for all $x\in D^1$ there exists a neighborhood
    $U(x)\subseteq D^{0,1}$ and a constant $K_x>0$ such that for all
    $y_1$ and $y_2\in U(x)$ 
    \begin{displaymath}
      \|F(y_1)-F(y_2)\|_0\leq K_x\|y_1-y_2\|_D\mbox{.}
    \end{displaymath}
  \item\label{thm:Fd01bound} $F$ maps elements of $D^{0,1}$ into
    $C^{0,1}(\T;\R^n)$. Moreover, for every $x\in D^{0,1}$, any
    bounded neighborhood $U(x)\subseteq D^{0,1}$ for which the
    Lipschitz constant $K_x$ exists has a bounded image under $F$:
    there exists a bound $R>0$ such that $\|F(y)\|_{0,1}\leq R$ for
    all $y\in U(x)$ ($R$ depends on $U(x)$).
  \end{enumerate}
\end{lemma}

\pstart{Proof} Statement~\ref{thm:Flip} is a consequence of the local
$EC$ Lipschitz continuity of $f$ and the compactness of
$\{\Delta_tx:t\in\T\}$ in $D^{0,1}$ (which allows one to choose a
uniform Lipschitz bound $K_x$ for all $t\in\T$). Note that the Lemma
requires $x\in D^1$ to ensure continuity of $t\to\Delta_tx$ with
respect to the $D^1$-norm (and, hence, with respect to the $D^{0,1}$-norm).

Concerning statement~\ref{thm:Fd01bound}: let $x\in D^{0,1}$ be
arbitrary, and let the neighborhood $U(x)$ be bounded (say,
$\|y-x\|_{D,L}\leq \delta$) such that $F$ has a Lipschitz constant
$K_x$ in $U(x)$.  Then we have for all $y,z\in U(x)$ and $t,s\in\T$
the estimate
\begin{displaymath}
  |f(\Delta_ty)-f(\Delta_sz)|\leq K_x\|\Delta_ty-\Delta_sz\|_D=
  K_x\|\Delta_{t-s}y-z\|_D\mbox{.}
\end{displaymath}
Initially setting $z=x$ and $s=t$ we get a bound on $\|F(y)\|_0$:
$\|F(y)\|_0\leq \|F(x)\|_0+K_x\delta=:R_0$ for all $y\in U(x)$.  It
remains to be shown that the Lipschitz constant of $F(y)$ is bounded
for $y\in U(x)$:
\begin{align*}
  |F(y)(t)-F(y)(s)|=|f(\Delta_ty)-f(\Delta_sy)|
  \leq K_x\|\Delta_ty-\Delta_sy\|_D
  \leq K_x\|y\|_{D,L}|t-s|\mbox{.}
\end{align*}
Since $\|y-x\|_{D,L}\leq \delta$ for $y\in U(x)$,
choosing
\begin{displaymath}
  R=\max\left\{R_0,K_x\left(\|x\|_{D,L}+\delta\right)\right\}  
\end{displaymath}
ensures that $\|F(y)\|_{0,1}\leq R$.  \eop

Define the maps 
\begin{align*}
  \partial^1F(u,v)(t)&=\partial^1f(\Delta_tu,\Delta_tv) 
  &&\mbox{for $u\in D^1$, $v\in D$,}\\
  \tilde A(u,v,w)(t)&=\tilde a(\Delta_tu,\Delta_tv,\Delta_tw)
  &&\mbox{for $u\in D^1$, $v\in D^1$, $w\in D$.}
\end{align*}
The following Lemma~\ref{thm:Fdiff}, and Lemma~\ref{thm:Fimage}
assume that $f$ is $EC^1$ smooth in $D$ in the sense of
Definition~\ref{def:extdiff}.
\begin{lemma}[Differentiability of $F$]\label{thm:Fdiff}
  Let $f:D\mapsto\R^n$ be $EC^1$ smooth. Then $F$, $\partial^1F$ and
  $\tilde A$ have the following properties.
  \begin{enumerate}
  \item\label{thm:Fcontdiff} The map $(u,v)\mapsto\partial^1F(u,v)$ is
    continuous in both arguments (and linear in its second argument)
    as a map from $D^1\times D$ into $C^0(\T;\R^n)$. It satisfies for
    all $x\in D^1$
    \begin{equation}\label{eq:Fcontdiff}
      \lim_{
        \begin{subarray}{c}
          y\in D^{0,1}\\[0.2ex]
          \|y\|_{D,L}\to0
        \end{subarray}
      }\frac{\|F(x+y)-F(x)-
        \partial^1F(x,y)\|_0}{\|y\|_{D,L}}=0\mbox{.}
    \end{equation}
  \item\label{thm:Fmeanval} The map $\tilde A(u,v,w)$ is continuous in
    all three arguments (and linear in its third argument) as a map
    from $D^1\times D^1\times D$ into $C^0(\T;\R^n)$. It satisfies for
    all $x,y\in D^1$
    \begin{displaymath}
      F(x+y)-F(x)=\tilde A_1(x,y,y)\mbox{.}
    \end{displaymath}
    Moreover, $\tilde A(x,0,y)=\partial^1F(x,y)$ for all $x\in D^1$
    and $y\in D$.
  \end{enumerate}
\end{lemma}
Note that in the limit \eqref{eq:Fcontdiff} we allow for deviations
$y\in D^{0,1}$.

\pstart{Proof} The continuity of $\partial^1F$ follows from the
continuity of $\partial^1f$ by applying Lemma~\ref{thm:Fcont} to
$\partial^1f:D^1\times D\mapsto\R^n$ instead of $f$. The linearity of
$\partial^1F$ in its second argument follows from the linearity of
$\partial^1f$ in its second argument.

The limit \eqref{eq:Fcontdiff} also follows from the corresponding
limit \eqref{eq:ass:contdifflip}: let $x\in D^1$ and $\epsilon>0$ be
arbitrary. For every fixed $t$ there exists a $\delta(t)>0$ such that
\begin{equation}\label{eq:Fcontdiffproof:ineq}
  \frac{|f(\Delta_tx+\Delta_ty)-f(\Delta_tx)-
    \partial^1f(\Delta_tx,\Delta_ty)|}{\|y\|_{D,L}}<\epsilon
\end{equation}
for all $y$ with $\|y\|_{D,L}<\delta(t)$.  As $f$ and $\partial^1f$
are continuous in their arguments $x\in D^1$ and $y\in D^{0,1}$, the
inequality also holds for all $s$ in a sufficiently small open
neighborhood of $t$, $U(t)$. The set of neighborhoods $U(t)$ for all
$t\in\T$ are a cover of the compact set $\T$ by open sets. Choosing a
finite subcover from this cover, and labeling the times
$t_1,\ldots,t_m$, we can choose
\begin{displaymath}
  \delta=\min_{k=1,\ldots,m}\delta(t_k)
\end{displaymath}
such that \eqref{eq:Fcontdiffproof:ineq} holds for all uniformly
$t\in\T$. This proves statement~\ref{thm:Fcontdiff} of the lemma.

Concerning statement \ref{thm:Fmeanval}: for the continuity of $\tilde
A$ we invoke again Lemma~\ref{thm:Fcont}, this time for $\tilde a$ on
$D^1\times D^1\times D$. The linearity of $\tilde A$ in its third
argument follows from the linearity of $\tilde a$ in its third
argument. The relations $F(x+y)(t)-F(x)(t)=\tilde A(x,y,y)(t)$ and
$\tilde A(x,0,y)(t)=\partial^1F(x,y)(t)$ in every $t\in\T$
follow from the corresponding relations for $f$ and $\tilde a$, as
stated in Lemma~\ref{thm:fmeanval}. \eop

\begin{lemma}[Differentiability of images of $F$]\label{thm:Fimage}
  Let $f:D\mapsto\R^n$ be $EC^1$ smooth and let $k\geq0$ be some
  integer. We assume that $F:D\mapsto C^k(\T;\R^n)$ and
  $\partial^1F:D^1\times D\mapsto C^k(\T;\R^n)$ are continuous maps.
  Then $F$ maps elements of $D^1$ into $C^{k+1}(\T;\R^n)$, and $F$ is
  continuous as a map from $D^1$ to $C^{k+1}(\T;\R^n)$.
\end{lemma}

\pstart{Proof}
Let $x$ be in $D^1$, that is, $x'\in D$. If $\partial^1F:D^1\times
D\mapsto C^k(\T;\R^n)$ is continuous then $\tilde A:D^1\times
D^1\times D\mapsto C^k(\T;\R^n)$, which is given by $\tilde
A(u,v,w)=\int_0^1\partial^1F(u+sv,w)\d s$, is continuous, too. Using
statement~\ref{thm:Fmeanval} of Lemma~\ref{thm:Fdiff} we have
\begin{align}
  \frac{F(\Delta_hx)-F(x)}{h}&=\tilde
  A\left(x,\Delta_hx-x,\frac{\Delta_hx-x}{h}\right)\mbox{.}
  \label{eq:Fimage:meanval}
\end{align}
On the right side $\|\Delta_hx-x\|_{D,1}$ converges to $0$ for
$h\to0$. Also,
\begin{displaymath}
  \left\|\frac{\Delta_hx-x}{h}-x'\right\|_D\to0\mbox{\quad for $h\to0$,}
\end{displaymath}
because $x\in D^1$. Since $\tilde A$ is continuous in its arguments
the right side converges to $\tilde A(x,0,x')=\partial^1F(x,x')\in
C^k(\T;\R^n)$ for $h\to0$. Thus, the limit of the left-hand side in
\eqref{eq:Fimage:meanval} for $h\to0$ exists, too, such that $F(x)\in
C^{k+1}(\T;\R^n)$ and the time derivative $(F(x))'$ equals
$\partial^1F(x,x')$.  Since $(v,w)\in D^1\times
D\mapsto \partial^1F(v,w)\in C^k(\T;\R^n)$ is continuous in $(u,v)$,
the time derivative of $F(x)$, $(F(x))'=\partial^1F(x,x')$ is also
continuous in $x$ if we use the norm $\|\cdot\|_{D,1}$ for the
argument and $\|\cdot\|_k$ for the image. \eop
\end{document}